\documentclass[12pt]{smfart}

\usepackage{smfenum}
\usepackage[francais]{babel}
\usepackage{smfthm}
\usepackage{amssymb}
\usepackage{amscd}
\usepackage{mathrsfs}
\usepackage[all]{xypic}

\newcommand{\Qp}{\mathbf{Q}_p}
\newcommand{\Zp}{\mathbf{Z}_p}
\newcommand{\Cp}{\mathbf{C}_p}
\newcommand{\ZZ}{\mathbf{Z}}
\newcommand{\Qpbar}{\overline{\mathbf{Q}}_p}
\newcommand{\Ind}{\mathrm{Ind}}
\newcommand{\Res}{\mathrm{Res}}
\newcommand{\zn}{\zeta_{p^n}}
\newcommand{\RR}{\mathbf{R}}
\newcommand{\Knr}{K^\mathrm{nr}}
\newcommand{\Fr}{\mathrm{Fr}}
\newcommand{\OO}{\mathcal{O}}
\newcommand{\bcris}{\mathbf{B}_{\mathrm{cris}}}
\newcommand{\bst}{\mathbf{B}_{\mathrm{st}}}
\newcommand{\bdR}{\mathbf{B}_{\mathrm{dR}}}
\newcommand{\Fil}{\mathrm{Fil}}
\newcommand{\dcris}{\mathbf{D}_{\mathrm{cris}}}
\newcommand{\dpst}{\mathbf{D}_{\mathrm{pst}}}
\newcommand{\ddR}{\mathbf{D}_{\mathrm{dR}}}
\newcommand{\eps}{\varepsilon}
\newcommand{\iso}{\overset{\sim}{\to}}
\newcommand{\Tam}{\mathrm{Tam}}
\newcommand{\calD}{\mathcal{D}}
\newcommand{\calE}{\mathcal{E}}
\newcommand{\calF}{\mathcal{F}}
\newcommand{\calH}{\mathcal{H}}
\newcommand{\calK}{\mathcal{K}}
\newcommand{\Exp}{\mathrm{Exp}}
\newcommand{\pscal}[1]{\langle #1 \rangle}
\newcommand{\et}{\widetilde{\mathbf{E}}}
\newcommand{\afont}{\mathbf{A}}
\newcommand{\bfont}{\mathbf{B}}
\newcommand{\dfont}{\mathbf{D}}
\newcommand{\aplus}{\mathbf{A}^+}
\newcommand{\bplus}{\mathbf{B}^+}
\newcommand{\nwach}{\mathbf{N}}
\newcommand{\bhol}{\mathbf{B}^+_{\mathrm{rig},K}}

\renewcommand{\geq}{\geqslant}
\renewcommand{\leq}{\leqslant} 
\renewcommand{\projlim}{\varprojlim}
\renewcommand{\injlim}{\varinjlim}
\renewcommand{\tilde}{\widetilde}
\renewcommand{\hat}{\widehat}
\renewcommand{\det}{\operatorname{det}}

\author[D. Benois]{Denis Benois}
\address{Facult\'e des Sciences et Techniques \\
Universit\'e de Besan\c{c}on \\ 
16 route de Gray \\
25030 Besan\c{c}on cedex \\ 
France}
\email{denis.benois@math.univ-fcomte.fr}

\author[L. Berger]{Laurent Berger}
\address{CNRS \& IHES \\
Le Bois-Marie\\
35 route de Chartres\\
91440 Bures-sur-Yvette \\ 
France}
\email{laurent.berger@ihes.fr}
\urladdr{www.ihes.fr/\~{}lberger/}

\title{Th\'eorie d'Iwasawa des repr\'esentations cristallines II} 

\date{Septembre 2005}

\NumberTheoremsIn{subsection}

\begin{document}

\begin{abstract}
Soit $K$ une extension finie non-ramifi\'ee de $\Qp$ et $V$ une repr\'esentation cristalline de $\mathrm{Gal}(\Qpbar/K)$. Dans cet article, on montre la conjecture $C_{\mathrm{EP}}(L,V)$ pour $L \subset \Qp^{\mathrm{ab}}$ et sa version \'equivariante $C_{\mathrm{EP}}(L/K,V)$ pour $L \subset \cup_{n=1}^\infty K(\zn)$. Les principaux ingr\'edients sont la conjecture $\delta_{\Zp}(V)$ sur l'int\'egralit\'e de l'exponentielle de Perrin-Riou, que nous d\'emontrons en utilisant la th\'eorie des $(\varphi,\Gamma)$-modules, et des techniques de descente en th\'eorie d'Iwasawa pour montrer que $\delta_{\Zp}(V)$ implique $C_{\mathrm{EP}}(L/K,V)$.
\end{abstract}

\begin{altabstract}
Let $K$ be a finite unramified extension of $\Qp$ and let $V$ be a crystalline representation of $\mathrm{Gal}(\Qpbar/K)$.
In this article, we give a proof of the $C_{\mathrm{EP}}(L,V)$ conjecture for $L \subset \Qp^{\mathrm{ab}}$ as well as a proof of its equivariant version $C_{\mathrm{EP}}(L/K,V)$ for $L \subset \cup_{n=1}^\infty K(\zn)$. The main ingredients are the $\delta_{\Zp}(V)$ conjecture about the integrality of Perrin-Riou's exponential, which we prove using the theory of $(\varphi,\Gamma)$-modules, and Iwasawa-theoretic descent techniques used to show that $\delta_{\Zp}(V)$ implies $C_{\mathrm{EP}}(L/K,V)$.
\end{altabstract}

\maketitle

\setcounter{tocdepth}{2}

\tableofcontents

\setlength{\baselineskip}{20pt}

\section*{Introduction}

Soient $p$ un nombre premier impair, $K$ une extension finie de $\Qp$ et $V$ une repr\'esentation potentiellement semi-stable de $G_K = \mathrm{Gal}(\Qpbar/K)$. Fontaine et Perrin-Riou ont formul\'e une conjecture qu'ils ont appel\'ee $C_{\mathrm{EP}}(K,V)$, conjecture qui entra\^{\i}ne la compatibilit\'e de la conjecture de Bloch et Kato sur les valeurs sp\'eciales des fonctions $L$ avec l'\'equation fonctionnelle. L'objet de ce texte est de montrer la conjecture $C_{\mathrm{EP}}(L,V)$ pour toute extension finie $L$ de $K$ telle que $L \subset \Qp^{\mathrm{ab}}$, 
quand $K$ est non-ramifi\'e sur $\Qp$ et $V$ est une repr\'esentation cristalline de $G_K$ ainsi que, sous les m\^emes hypoth\`eses, la version \'equivariante $C_{\mathrm{EP}}(L/K,V)$ de cette conjecture pour toute extension finie $L$ de $K$ contenue dans $K_\infty = \cup_{n=1}^\infty K(\zn)$. Comme ingr\'edient de la d\'emonstration, on montre aussi la conjecture $\delta_{\Zp}(V)$ de Perrin-Riou, que nous appelons $C_{\mathrm{Iw}}(K_\infty/K,V)$ en raison de son lien avec la th\'eorie d'Iwasawa de $V$.

Rappelons tout d'abord la conjecture $C_{\mathrm{EP}}(L/K,V)$. Pour cela, on se donne une extension ab\'elienne finie $L/K$ de groupe de Galois $G=\mathrm{Gal}(L/K)$, une repr\'esentation potentiellement semi-stable $V$ de $G_K$ et un r\'eseau $T$ de $V$ stable sous l'action de $G_K$. On d\'efinit la droite d'Euler-Poincar\'e de $V$ en posant :
\[ \Delta_{\mathrm{EP}}(L/K,V) = \det_{\Qp[G]}\RR\Gamma(L,V) \otimes \det_{\Qp[G]}(\Ind_{L/\Qp}V). \]
On sait que $\RR\Gamma(L,T)$ est un complexe parfait de $\Zp[G]$-modules et que l'image de $\Delta_{\mathrm{EP}}(L/K,T) = \det_{\Zp[G]}\RR\Gamma(L,T) \otimes
\det_{\Zp[G]}(\Ind_{L/\Qp}T)$ dans $\Delta_{\mathrm{EP}}(L/K,V)$ ne d\'epend pas du choix de $T$. 

On note $\dcris^L(V)$, $\dpst(V)$ et $\ddR^L(V)$ les modules associ\'es \`a la restriction de $V$ \`a $G_L$ par la th\'eorie de Fontaine, et $t_V(L) = \ddR^L(V) /
\Fil^0 \ddR^L(V)$ l'espace tangent de $V$ sur $L$. La suite exacte :
\begin{multline*}
0 \to H^0(L,V) \to \dcris^L(V) \to \dcris^L(V) \oplus t_V(L) \to H^1(L,V) \to \\
\to \dcris^L(V^*(1))^* \oplus t^*_{V^*(1)}(L) \to \dcris^L(V^*(1))^* \to H^2(L,V) \to 0,
\end{multline*}
qui provient de la suite exacte fondamentale (cf \S\ref{kato}) et l'isomorphisme $t_{V^*(1)}^*(L) \simeq \Fil^0 \ddR^L(V)$ donnent un isomorphisme canonique :
\[ \det_{\Qp[G]} \RR\Gamma(L,V) \iso \det^{-1}_{\Qp[G]} \ddR^L(V). \]

La th\'eorie des constantes locales permet d'autre part de d\'efinir un \'el\'ement $\eps(L/K,V) \in \Qp(\zeta_{p^\infty})[G]$ associ\'e \`a l'action de $G_K$ sur $\dpst(V)$. L'isomorphisme de comparaison :
\[ \bdR \otimes_{\Qp} \Ind_{L/\Qp} (V) \simeq \bdR \otimes_{\Qp} \ddR^L(V), \]
normalis\'e par $\eps(L/K,V)$ et par le facteur $\Gamma$ habituel $\Gamma^*(V)$,
fournit un isomorphisme :
\[ \det_{\Qp[G]}^{-1} \ddR^L(V) \otimes \det_{\Qp[G]} \Ind_{L/\Qp} (V)
\simeq \Qp[G]_{V,L/K}, \]
o\`u $\Qp[G]_{V,L/K}$ est un certain $\Qp[G]$-module libre de rang $1$ qui contient un sous-$\Zp[G]$-module inversible canonique $\Zp[G]_{V,L/K}$ (cf d\'efinition \ref{defzpvlk}). En composant ces isomorphismes, on obtient une trivialisation canonique de la droite d'Euler-Poincar\'e :
\[ \delta_{V,L/K} : \Delta_{\mathrm{EP}}(L/K,V) \simeq \Qp[G]_{V,L/K}. \]
Dans son manuscrit non-publi\'e \cite{K2}, Kato a propos\'e la conjecture suivante (qu'il appelle {\og local $\eps$-conjecture \fg}) : 
\begin{enonce*}{Conjecture $C_{\mathrm{EP}}(L/K,V)$}
Si $V$ est une repr\'esentation potentiellement semi-stable et si $L/K$ est une
extension ab\'elienne finie, alors l'application $\delta_{V,L/K}$ envoie $\Delta_{\mathrm{EP}}(L/K,T)$ sur $\Zp[G]_{V,L/K}$.
\end{enonce*}

C'est la conjecture \ref{new244} de cet article. Si $L=K$, alors on retrouve la conjecture 
$C_{\mathrm{EP}}(K,V)$ de Fontaine et Perrin-Riou que l'on peut d'ailleurs reformuler en termes de nombres de Tamagawa (cf conjecture \ref{new245}).

Rappelons \`a pr\'esent la conjecture $C_{\mathrm{Iw}}(K_\infty/K,V)$. On suppose pour cela que $K$ est non-ramifi\'e, on fixe une suite compatible de racines primitives $p^n$-i\`emes de l'unit\'e $\eps=(\zn)_{n \geq 0}$ et pour $n \geq 1$, on pose $K_n=K(\zn)$ ainsi que $K_\infty = \cup_{n \geq 1} K_n$. Soient $H_K = \mathrm{Gal}(\Qpbar / K_\infty)$, $\Gamma=\mathrm{Gal}(K_\infty/K)$ et $\Gamma_n=\mathrm{Gal}(K_\infty/K_n)$ ce qui fait que $\Gamma  = \Delta_K \times \Gamma_1$ o\`u $\Delta_K$ est le sous-groupe de torsion de $\Gamma$. Soit $\calH$ l'alg\`ebre des s\'eries formelles $f(X) \in \Qp[[X]]$ qui convergent sur le disque unit\'e ouvert et $\calH(\Gamma_1) = \{ f(\gamma_1-1) \mid \gamma_1 \in \Gamma_1$ et $f \in \calH \}$. On pose $\Lambda = \Zp[[\Gamma]]$, $\calH(\Gamma) = \Qp[\Delta_K] \otimes_{\Qp} \calH(\Gamma_1)$ et $\calK(\Gamma)$ est l'anneau total des fractions de $\calH(\Gamma)$. On d\'efinit la cohomologie d'Iwasawa d'une repr\'esentation $V$ en posant : 
\[ H^1_{\mathrm{Iw}}(K,T)=\projlim_{\mathrm{cor}_{K_n / K_{n-1}}} H^1(K_n,T), \] 
et $H^1_{\mathrm{Iw}} (K,V) = \Qp \otimes_{\Zp} H^1_{\mathrm{Iw}} (K,T)$. 

Supposons \`a pr\'esent que $V$ est cristalline. Dans \cite{PR1}, Perrin-Riou a construit une famille d'applications :
\[ \Exp_{V,h}^\eps : \calD(V)^{\Delta =0} \to \calH(\Gamma) \otimes_{\Lambda} H^1_{\mathrm{Iw}}(K,V) / V^{H_K}, \]
qui interpolent les exponentielles de Bloch et Kato. Plus pr\'ecisement, pour tout  $h \geq 1$ v\'erifiant $\Fil^{-h}\ddR(V) = \ddR(V)$, on a un diagramme commutatif :  
\[ \begin{CD}
\calD(V)^{\Delta=0} @>\Exp^{\eps}_{V,h}>> \calH(\Gamma) \otimes_{\Lambda} H^1_{\mathrm{Iw}}(K,V) / V^{H_K} \\
@V\Xi^{\eps}_{V,n}VV         @VV\mathrm{pr}_{T,n}V \\
\ddR^{K_n}(V) @>(h-1)!\exp_{V,K_n}>> 
H^1(K_n, V) / H^1(\Gamma_n, V^{H_K}).
\end{CD} \]

Ici, $\calD(V)$ est isomorphe \`a $\Lambda \otimes_{\Zp} \dcris(V)$ et les applications $\Delta$ et $\Xi^{\eps}_{V,n}$ sont explicites, mais leur d\'efinition est un peu technique pour cette introduction (cf paragraphe \ref{benois}). Cette construction joue un r\^ole important dans la th\'eorie des fonctions $L$ $p$-adiques (cf \cite{PR2} et \cite{C3}). Posons maintenant :
\[ \Delta_{\mathrm{Iw}}(K_{\infty}/K,T) = \det_{\Lambda} \RR\Gamma_{\mathrm{Iw}}(K,T)\otimes\det_{\Lambda}(\Ind_{K_\infty/\Qp} T), \]
et $\Delta_{\mathrm{Iw}}(K_{\infty}/K,V) = \Qp \otimes_{\Zp} \Delta_{\mathrm{Iw}}(K_{\infty}K,T)$. On pose $\ell_j=j- \log \gamma_1 / \log \chi (\gamma_1)$ et on d\'efinit un facteur $\Gamma$ par la formule : 
\[ \mathbf{\Gamma}_h(V) = \prod_{j > -h} 
(\ell_{-j})^{\dim_{\Qp}\Fil^j\dcris (V)}. \]

Le d\'eterminant de $\Exp_{V,h}^\eps$ normalis\'e par $\mathbf{\Gamma}_h(V)^{-1}$ ne d\'epend alors pas de $h$, et la loi de r\'eciprocit\'e de Perrin-Riou entra\^{\i}ne qu'il induit un isomorphisme canonique :
\[ \delta_{V,K_\infty/K} : \Delta_{\mathrm{Iw}}(K_{\infty}/K,V)
\to \Qp \otimes_{\Zp} \Lambda_{V,K_\infty/K}, \]
o\`u $\Lambda_{V,K_\infty/K}$ est un certain $\Lambda$-module libre de rang $1$ (cf le paragraphe \ref{nobody}). Perrin-Riou a propos\'e la conjecture suivante (appel\'ee $\delta_{\Zp}(V)$ dans \cite{PR1} et \cite{PR2}) relativement au d\'eterminant de $\Exp_{V,h}^\eps$.
\begin{enonce*}{Conjecture $C_{\mathrm{Iw}}(K_\infty/K,V)$} 
Si $V$ est une repr\'esentation cristalline de $G_K$, alors l'application
$\delta_{V,K_\infty/K}$ envoie $\Delta_{\mathrm{Iw}}(K_{\infty}/K,T)$ sur $\Lambda_{V,K_\infty/K}$.
\end{enonce*}

Le r\'esultat principal de cet article est le suivant : 
\begin{enonce*}{Th\'eor\`eme A}
Si $K$ est une extension non-ramifi\'ee de $\Qp$ et si $V$ est une repr\'esen\-tation cristalline de $G_K$, alors :
\begin{enumerate}
\item la conjecture $C_{\mathrm{Iw}}(K_\infty/K,V)$ est vraie;
\item la conjecture $C_{\mathrm{EP}}(L/K,V)$ est vraie pour toute extension 
finie $L$ de $K$ contenue dans $K_\infty$.
\end{enumerate}
\end{enonce*}
En utilisant les propri\'et\'es fonctorielles de la conjecture $C_{\mathrm{EP}}(L/K,V)$, on en d\'eduit le corollaire suivant :
\begin{enonce*}{Corollaire B}
Si $K$ est une extension non-ramifi\'ee de $\Qp$ et si $V$ est une repr\'esen\-tation cristalline de $G_K$, alors :
\begin{enumerate}
\item la conjecture $C_{\mathrm{EP}}(L,V)$ est vraie pour toute extension $L/K$ contenue dans $\Qp^{\mathrm{ab}}$; 
\item la conjecture $C_{\mathrm{EP}}(K,V(\eta))$ est vraie pour tout caract\`ere de Dirichlet $\eta$ de $\Gamma$.
\end{enumerate}
\end{enonce*}

Le th\'eor\`eme A et le corollaire B sont d\'emontr\'es \`a la fin de cet article (cf le th\'eor\`eme \ref{new345} et le corollaire \ref{old436}). Disons quelques mots du plan de l'article. Les chapitres \ref{repsst} et \ref{constloc} sont consacr\'es \`a des rappels, qui aboutissent \`a l'\'enonc\'e de la conjecture $C_{\mathrm{EP}}(L/K,V)$. Les chapitres \ref{expbpr} et \ref{nobalt} sont le coeur technique de l'article. On commence par y rappeler la construction de l'exponentielle de Perrin-Riou, puis on y \'enonce la conjecture $C_{\mathrm{Iw}}(K_\infty/K,V)$. Apr\`es cela on montre dans les paragraphes \S\S\ref{same1}, \ref{same2}, en utilisant des techniques de descente en th\'eorie d'Iwasawa, que la conjecture $C_{\mathrm{Iw}}(K_\infty/K,V)$ est \'equivalente \`a la conjecture $C_{\mathrm{EP}}(K_n/K,V)$ pour tout $n \geq 1$. Enfin dans le \S\ref{main} on d\'emontre la conjecture $C_{\mathrm{Iw}}(K_\infty/K,V)$.

Les m\^emes arguments, avec un peu plus de calculs, permettent de d\'emontrer la conjecture $C_{\mathrm{EP}}(L/K,V)$ pour toute extension $L/K$ contenue dans $\Qp^{\mathrm{ab}}$. Cette petite g\'en\'eralisation est importante pour la version \'equivariante des conjectures de Bloch et Kato; nous en laissons les d\'etails au lecteur.

Pour terminer cette introduction, remarquons que dans le cas o\`u $V$ est ordinaire, ces r\'esultats \'etaient d\'ej\`a connus (voir \cite{PR1, BN, BF3}).

\section{Repr\'esentations potentiellement semi-stables}\label{repsst}

Dans tout cet article, le corps $K$ est une extension finie de $\Qp$ (dans le chapitre \ref{expbpr}, on suppose qu'elle est non-ramifi\'ee). L'anneau des entiers de $K$ est not\'e $\OO_K$ et son corps r\'esiduel $k_K$ est de cardinal $q_K$.
On fixe une fois pour toutes une suite compatible de racines primitives $p^n$-i\`emes de l'unit\'e $\eps=(\zn)_{n \geq 0}$ et pour $n \geq 1$, on pose $K_n=K(\zn)$ ainsi que $K_\infty = \cup_{n \geq 1} K_n$. La notation $K_0$ d\'esigne le sous-corps maximal non-ramifi\'e de $K$.

On pose :
\begin{alignat*}{2}
G_K & = \mathrm{Gal}(\Qpbar/K)  \qquad H_K & = \mathrm{Gal}(\Qpbar/K_\infty) \\
\Gamma_n & = \mathrm{Gal}(K_\infty/K_n)  \qquad G_n & = \mathrm{Gal}(K_n/K)
\end{alignat*}
et $\Lambda=\Zp[[\Gamma]]$ est l'alg\`ebre d'Iwasawa de $\Gamma$. Profitons-en pour remarquer que le caract\`ere cyclotomique $\chi$ envoie $\Gamma$ dans $\Zp^\times$ et que cette application est un isomorphisme si $K$ est non-ramifi\'e sur $\Qp$.

L'objet de ce chapitre est de donner quelques rappels, sur la th\'eorie de Hodge $p$-adique, la th\'eorie des $(\varphi,\Gamma)$-modules, la cohomologie galoisienne et l'exponentielle de Bloch-Kato.

\Subsection{Th\'eorie de Hodge $p$-adique}\label{hodge}

Dans ce paragraphe, on rappelle quelques unes des constructions de Fontaine (voir \cite{F3,F4}) qui sont utilis\'ees dans la suite de cet article. On note $\sigma$ le Frobenius arithm\'etique absolu agissant sur $\Qp^{\mathrm{nr}}$. 

Soient $\bcris$, $\bst$ et $\bdR$ les anneaux de p\'eriodes $p$-adiques construits par Fontaine (voir \cite{F3} par exemple). Le corps $\bdR = \bdR^+[1/t]$ est une $\Qp$-alg\`ebre qui contient $\Qpbar$ et qui est munie d'une action de $G_K$ ainsi que d'une filtration d\'ecroissante exhaustive et s\'epar\'ee par des $\Fil^i \bdR = t^i \bdR^+$. 
Remarquons que l'uniformisante $t = \log [\eps]$ d\'epend du choix de
$\eps=(\zn)_{n \geq 0}$ que l'on a fait ci-dessus.
L'anneau $\bst$ est une $\Qp$-alg\`ebre qui contient $\hat{\mathbf{Q}}_p^{\mathrm{nr}}$ et qui est munie d'une action de $G_K$ ainsi que d'un endomorphisme $\varphi$ commutant \`a l'action de $G_K$ et $\sigma$-semi-lin\'eaire et d'un op\'erateur de monodromie $N : \bst \to \bst$ qui commute \`a l'action de $G_K$ et v\'erifie $N \circ \varphi = p \varphi \circ N$.  Enfin, $\bcris=\bst^{N=0}$. On a donc $\bcris \subset \bst$ et de plus on a une injection $\Qpbar \otimes_{\Qp^{\mathrm{nr}}}Ê \bst \hookrightarrow \bdR$.

Pour toute repr\'esentation $p$-adique $V$ de $G_K$, on pose $\ddR^K(V) = (\bdR \otimes_{\Qp} V)^{G_K}$, ce qui fait que $\ddR^K(V)$ est un $K$-espace vectoriel filtr\'e de dimension finie. S'il n'y a pas de confusion possible quant au corps $K$, on \'ecrit plus simplement $\ddR(V)$. De mani\`ere analogue on pose : 
\[ \dcris^K(V) = (\bcris \otimes_{\Qp} V)^{G_K} \quad\text{et}\quad \dpst(V) = \injlim_{L/K} (\bst \otimes_{\Qp} V)^{G_L}, \] 
o\`u $L$ parcourt l'ensemble des extensions finies de $K$, ce qui fait de $\dcris^K(V)$ un $K_0$-espace vectoriel muni d'une action $\sigma$-semi-lin\'eaire de $\varphi$ et de $\dpst(V)$ un $K_0^{\mathrm{nr}}$-espace vectoriel muni des op\'erateurs $\varphi$ et $N$ v\'erifiants $N \circ \varphi = p \varphi \circ N$. Comme ci-dessus, on \'ecrit $\dcris(V)$ s'il n'y a pas de confusion possible. On a : 
\[ \dim_{K_0} \dcris^K(V) \leq \dim_{K_0^{\mathrm{nr}}} \dpst(V) \leq \dim_K \ddR^K(V) \leq \dim_{\Qp} V. \]
On dit que $V$ est cristalline (resp. potentiellement semi-stable,
resp. de de Rham) si $\dim_{K_0} \dcris^K(V) = \dim_{\Qp} V$ (resp. si $\dim_{K_0^{\mathrm{nr}}} \dpst(V) = \dim_{\Qp} V$, resp. si $\dim_K \ddR^K(V) = \dim_{\Qp} V$). Remarquons que l'on sait maintenant que toute repr\'esentation de de Rham est potentiellement semi-stable (cf \cite{Ber1}).

Si $V$ est une repr\'esentation de de Rham, on pose : 
\[ h_i(V) = \dim_K (\Fil^i \ddR^K(V) / \Fil^{i+1}\ddR^K(V)). \] 
La d\'ecomposition de Hodge-Tate de $V$ s'\'ecrit alors $\Cp \otimes_{\Qp} V \simeq \oplus_{i \in \ZZ} \Cp(-i)^{h_i(V)}$ o\`u $\Cp$ est le compl\'et\'e $p$-adique de $\Qpbar$. Les oppos\'es des entiers $i$ tels que $h_i(V) \neq 0$ sont les poids de Hodge-Tate de $V$. On pose $t_H(V) = \sum_{i \in \ZZ} i h_i(V)$. 

\Subsection{Modules de Wach et $(\varphi,\Gamma)$-modules}\label{wach}

Soit $K/\Qp$ une extension finie que l'on suppose ici non-ramifi\'ee. On note 
$\et^+ = \varprojlim (\OO_{\Qpbar}/ p \OO_{\Qpbar})$ l'anneau construit par
Fontaine (voir \cite{F2} par exemple, cet anneau s'y appelle $\mathcal{R}$), $\et=\mathrm{Frac}(\et^+)$ son corps des fractions et $W(\et)$ l'anneau des vecteurs de Witt \`a coefficients dans $\et$. On pose $X=[\eps]-1$,
avec $\eps=(\zn)_{n \geq 0}$, $\aplus_K = \OO_K[[X]]$ et on note $\afont_K$ le compl\'et\'e
$p$-adique  de $\aplus_K [1/X]$.  Les anneaux $\aplus_K$ et $\afont_K$ sont munis d'un Frobenius $\varphi$ et d'une action de $\Gamma = \mathrm{Gal}(K_\infty/K)$, donn\'es par les formules $\varphi(X) = (1+X)^p-1$ et $\gamma (X)=(1+X)^{\chi (\gamma)}-1$ pour $\gamma \in \Gamma$, o\`u $\chi : \Gamma \to \Zp^\times$ est le caract\`ere cyclotomique. 
Soit $\bfont$ le compl\'et\'e $p$-adique de l'extension maximale non-ramifi\'ee du corps $\bfont_K = \Qp \otimes_{\Zp} \afont_K $ dans $W(\et)$. On pose $\afont = \bfont \cap W(\et)$,  $\bplus = \bfont \cap W(\et^+)[1/p]$ et $\aplus = \afont \cap W(\et^+)$. Tous ces anneaux sont munis d'une action de $G_K$ et d'un Frobenius $\varphi$. Enfin, on a $\afont_K= \afont^{H_K}$.

Un $(\varphi,\Gamma)$-module est un module libre de rang fini sur 
$\afont_K$ muni d'un Frobenius semi-lin\'eaire $\varphi$ et d'une
action continue et semi-lin\'eaire de $\Gamma$ commutant
avec $\varphi$. Dans \cite{F2}, Fontaine a d\'efini un foncteur :
\[ \dfont : T \mapsto \dfont(T) = (\afont \otimes_{\Zp} T)^{H_K} , \]
qui fournit une \'equivalence entre la cat\'egorie des $\Zp$-repr\'esentations de $G_K$ et la cat\'egorie des $(\varphi,\Gamma)$-modules \'etales. Le foncteur :
\[ M \mapsto (\afont \otimes_{\afont_K} M)^{\varphi=1} \]
est un quasi-inverse de $\dfont$. De m\^eme, le foncteur $\dfont : V \mapsto (\bfont \otimes_{\Qp} V)^{H_K}$ donne une \'equivalence entre la cat\'egorie des repr\'esentations $p$-adiques de $G_K$ et la cat\'egorie des  $(\varphi,\Gamma)$-modules \'etales sur $\bfont_K = \Qp \otimes_{\Zp} \afont_K$.

Si $V$ est une repr\'esentation cristalline, alors un r\'esultat de Colmez \cite{C2} dit qu'il existe une base de $\dfont(T)$ dans laquelle les matrices de $\varphi$ et de $\gamma \in \Gamma$ sont \`a coefficients dans $\aplus_K$. Plus pr\'ecisement, on a le r\'esultat suivant, d\'emontr\'e dans \cite{Ber2} :

\begin{prop}\label{old412} 
Si $V$ est une repr\'esentation cristalline dont les oppos\'es des poids de Hodge-Tate sont $0 = r_1 \leq r_2 \leq  \cdots \leq r_d=h$, alors il existe un unique sous $\aplus_K$-module $\nwach(T)$ de $\dfont^+(T) = (\afont^+ \otimes_{\Zp} T)^{H_K}$, stable par $\varphi$ et qui satisfait les conditions suivantes :
\begin{enumerate}
\item $\nwach(T)$ est un $\aplus_K$-module libre de rang $d=\dim (V)$ et contient une base de $\dfont(T)$ sur $\afont_K$;
\item l'action de $\Gamma$ pr\'eserve $\nwach(T)$ et elle est triviale sur $\nwach(T)/X \nwach(T)$;
\item $X^h \dfont^+(T)\subset \nwach(T)$.
\end{enumerate}
De plus, on a $q^h \nwach(T) \subset \varphi^*\nwach(T)$, o\`u $q = \varphi (X) / X$ et $\varphi^*\nwach(T)$ est le $\aplus_K$-module engendr\'e par $\varphi(\nwach(T))$.
\end{prop}

Soit $\bhol$ l'ensemble des s\'eries formelles $f(X) = \sum_{k=0}^\infty a_k X^k$,  
avec $a_k \in K$ et telles que $f(X)$ converge sur le disque unit\'e ouvert
$\{x \in \Cp \mid |x|_p <1\}$. L'anneau $\bhol$ est de B\'ezout \cite{L} et de plus il admet la th\'eorie des diviseurs \'el\'ementaires; il est aussi muni d'actions de $\varphi$
et de $\Gamma$ et on a un plongement $\varphi^{-n} : \bhol \hookrightarrow K_n[[t]]
\subset \bdR^+$ qui envoie $X$ sur $\zn \exp{(t/{p^n})}-1$.

\begin{prop}\label{old414}
Si $V$ est une repr\'esentation cristalline dont les oppos\'es des poids de Hodge-Tate sont
$0 = r_1 \leq r_2\leq \hdots \leq r_d=h$, alors $\dcris(V)\simeq (\bhol \otimes_{\aplus_K} \nwach(T))^\Gamma$ et : 
\[ \left[  \bhol \otimes_{\aplus_K} \nwach(T) :
 \bhol \otimes_K \dcris(V) \right] =
\left [\left (\frac{t}{X} \right)^{r_1} ;\cdots ; 
\left (\frac{t}{X} \right )^{r_h} \right ]. \]
\end{prop}

\begin{proof} 
Voir \cite[prop III.4]{Ber2}. 
\end{proof}

\Subsection{Cohomologie galoisienne}\label{herr}

Rappelons maintenant comment on peut calculer la cohomologie galoisienne des repr\'esentations $p$-adiques \`a partir des $(\varphi,\Gamma)$-modules. 
On suppose toujours que $K/\Qp$ est non-ramifi\'ee, on pose $\Gamma_n = \mathrm{Gal}(K_\infty/K_n)$, on fixe un g\'en\'erateur topologique $\gamma_1$ de $\Gamma_1$ et on pose $\gamma_n=\gamma_1^{p^{n-1}}$. Si $T$ est une repr\'esentation $\Zp$-adique de $G_K$, on note $C_{\varphi,\gamma_n}(K_n,T)$ le complexe :
\[ 0 \to \dfont(T) \overset{f}{\to} \dfont(T) \oplus \dfont(T) \overset{g}{\to} \dfont(T) \to 0, \]
o\`u les applications $f$ et $g$ sont d\'efinies par $f(x)=((\varphi-1)x,(\gamma_n-1)x)$ et 
$g(y,z)=(\gamma_n-1)y -(\varphi-1)z$.

Dans \cite{H1}, Herr a montr\'e que les groupes de cohomologie $H^i(C_{\varphi,\gamma_n}(K_n,T))$ s'identifient canoniquement aux groupes de cohomologie galoisienne $H^i(K_n,T)$ (voir \cite[prop 1.3.2]{Ben} ou bien \cite[prop I.4.1]{CC} ou encore \cite[prop I.8]{Ber3} pour une description explicite de cet isomorphisme quand $i=1$).

Enfin, on peut aussi retrouver la cohomologie d'Iwasawa :
\[ H^i_{\mathrm{Iw}}(K,T) = \projlim_{\mathrm{cor}_{K_n / K_{n-1}}} H^i(K_n,T), \]
en utilisant les $(\varphi,\Gamma)$-modules. Pour cela, on utilise l'op\'erateur $\psi : \bfont \to \bfont$
qui est d\'efini par la formule :  
\[ \psi (x)= \frac{1}{p}\varphi^{-1}(\mathrm{Tr}_{\bfont/\varphi (\bfont)}(x)). \] 
L'op\'erateur $\psi$ commute \`a l'action de $G_K$ et on a $\psi \circ \varphi = \mathrm{id}$. La cohomologie du complexe : 
\[ \dfont(T) \overset{\psi-1}{\to} \dfont(T) \]
s'identifie canoniquement \`a la cohomologie d'Iwasawa de $T$, c'est-\`a-dire que 
$H_{\mathrm{Iw}}^1(K,T) \simeq \dfont(T)^{\psi=1}$ et que $H_{\mathrm{Iw}}^2(K,T) \simeq \dfont(T) / (\psi-1)$ (voir \cite[\S II.3]{CC}). Donnons une description explicite
du premier isomorphisme. Si $\alpha \in \dfont(T)^{\psi=1}$, alors $(\varphi-1)\alpha \in \dfont(T)^{\psi=0}$ et comme $\gamma_n -1$ est inversible sur $\dfont(T)^{\psi=0}$ (cf. \cite{H1} ou \cite[prop I.5.1]{CC}), il existe $x_n \in \dfont(T)$ v\'erifiant $(\gamma_n-1)x_n = (\varphi-1)\alpha$. Les $\mathrm{cl}(x_n,\alpha) \in  H^1(C_{\varphi,\gamma_n}(K_n,T))$ forment alors un syst\`eme compatible
d'\'el\'ements de $H^1(K_n,T)$.

\Subsection{L'exponentielle de Bloch-Kato}\label{kato}

Soit $V$ une repr\'esentation de de Rham. Bloch et Kato ont d\'efini  (voir \cite[\S4]{BK}) la partie exponentielle (resp. parties finie et g\'eom\'etrique) de $H^1(K,V)$ en posant :
\begin{align*}
&H^1_e(K,V)=\ker (H^1(K,V)\to H^1(K,\bcris^{\varphi=1} \otimes_{\Qp} V)),\\
&H^1_f(K,V)=\ker (H^1(K,V)\to H^1(K,\bcris \otimes_{\Qp} V)),\\
&H^1_g(K,V)=\ker (H^1(K,V)\to H^1(K,\bdR \otimes_{\Qp} V)).
\end{align*}

La dualit\'e locale fournit un accouplement $(\cdot,\cdot)_V : H^1(K,V) \times H^1(K,V^*(1)) \to \Qp$ pour lequel l'orthogonal de $H^1_e(K,V)$ est $H^1_g(K,V^*(1))$ et celui de $H^1_f(K,V)$ est $H^1_f(K,V^*(1))$. L'espace tangent de $V$ sur $K$ est par d\'efinition le quotient :
\[ t_V(K) = \ddR^K(V) / \Fil^0 \ddR^K(V). \]

Les anneaux $\bcris$ et $\bdR$ sont reli\'es par l'inclusion $\bcris \subset \bdR$ mais aussi et surtout par les suites exactes fondamentales :
\begin{align*}
0 \to \Qp \to \bcris^{\varphi=1} \overset{\alpha}{\to} \bdR / \Fil^0 \bdR \to 0, \\
0 \to \Qp \to \bcris \overset{\beta}{\to} \bcris \oplus \bdR / \Fil^0\bdR \to 0,
\end{align*}
o\`u $\alpha (x) = x \mod{\Fil^0 \bdR}$ et $\beta(x)=((1-\varphi)x,x \mod{\Fil^0 \bdR})$. En prenant les produits tensoriels de ces suites par $V$ et les invariants sous l'action de $G_K$, on obtient des suites exactes longues de cohomologie qui nous donnent les deux suites exactes :
\begin{align*}
0 \to H^0(K,V) \to \dcris(V)^{\varphi=1} \to t_V(K) \to H_e^1(K,V) \to 0, \tag{eq1} \\
0 \to H^0(K,V) \to \dcris(V) \to \dcris(V) \oplus t_V(K) \to H_f^1(K,V) \to 0. \tag{eq2} 
\end{align*}
L'application de connexion $\exp_{V,K} : t_V(K) \to H^1(K,V)$ dans la premi\`ere suite s'appelle l'exponentielle de Bloch et Kato. L'exponentielle duale $\exp^*_{V,K} : H^1(K,V) \to \Fil^0 \ddR^K (V)$ est d\'efinie par la formule : 
\[ \mathrm{Tr}_{K/\Qp}[\exp^*_{V,K}(x),y]_V=
(x,\exp_{V^*(1),K}(y))_V, \]
o\`u $[\cdot,\cdot]_V : \ddR^K(V) \times \ddR^K(V^*(1)) \to K$ est la dualit\'e
canonique. On v\'erifie facilement que $\ker (\exp^*_{V,K}) = H^1_g(K,V)$.

\begin{lemm}\label{old232}
On a des isomorphismes canoniques :
\begin{align*}
\exp_{V,f/e} & : \frac{\dcris(V)}{(1-\varphi)\dcris(V)} \iso 
\frac{H^1_f(K,V)}{H^1_e(K,V)},\\
\exp^*_{V,g/f} & : \frac{H^1_g(K,V)}{H^1_f(K,V)}
\iso  \dcris(V)^{\varphi=p^{-1}}.
\end{align*}
\end{lemm}

\begin{proof} 
On remarque que la suite (eq1) s'injecte dans (eq2), d'o\`u on obtient le premier isomorphisme. D'autre part, pour tout $K_0$-espace vectoriel $W$ de dimension finie muni d'un op\'erateur $\sigma$-semi-lin\'eaire $\varphi$, on a un isomorphisme :
\begin{align*}
\mathrm{Hom}_{K_0}(W,K_0)^{\varphi=1} & \iso  
\mathrm{Hom}_{\Qp}(W/(1-\varphi)W,\Qp) \\
f & \mapsto \mathrm{Tr}_{K_0/\Qp} f,
\end{align*}
ce qui fait que :
\[ \frac{H^1_g(K,V)}{H^1_f(K,V)}\simeq
\left(\frac{H^1_f(K,V^*(1))}{H^1_e(K,V^*(1))}\right)^* \simeq
\left(\frac{\dcris(V^*(1))}{(1-\varphi)\dcris(V^*(1))}\right)^*
\iso \dcris(V)^{\varphi=p^{-1}}, \]
et le lemme est d\'emontr\'e.
\end{proof}

On pose maintenant $L_f(K,V) = \det_{\Qp} H^0(K,V) \otimes_{\Qp} \det_{\Qp}^{-1}H^1_f(K,V)$. La suite exacte (eq2) fournit alors un isomorphisme canonique $i_V : L_f(K,V) \simeq \det_{\Qp}^{-1}t_V(K)$. Soit $T$ un $\Zp$-r\'eseau de $V$ stable sous l'action de $G_K$ et soit $\omega$ une base de $\det_{\Qp} t_V(K)$. On note $H_f^1(K,T)$ l'image inverse de $H_f^1(K,V)$ dans $H^1(K,T)$ et l'on pose : 
\[ L_f(K,T) = \det_{\Zp} H^0(K,T) \otimes_{\Zp} \det_{\Zp}^{-1} H^1_f(K,T). \]

\begin{defi}\label{deftamag}
On appelle nombre de Tamagawa, et on note $\Tam^0_{K,\omega}(T)$, l'unique puissance de $p$ telle que $i_V(L_f(K,T)) = \Zp \Tam^0_{K,\omega}(T) \omega^{-1}$, o\`u $\omega^{-1}$ est la base duale de $\omega$ (voir \cite{BK,PR2}).
\end{defi}

Ces nombres interviennent dans la formulation de la conjecture $C_{\mathrm{EP}}(K,V)$ (conjecture \ref{new245} ci-dessous).

\section{D\'eterminants et constantes locales}\label{constloc}

L'objet de ce chapitre est d'\'enoncer la conjecture $C_{\mathrm{EP}}(L/K,V)$. On commence par des rappels sur la th\'eorie des d\'eterminants g\'en\'eralis\'es, puis on passe en revue la construction des constantes locales, pour les repr\'esentations de Weil-Deligne tout d'abord, et pour les repr\'esentations potentiellement semi-stables ensuite.

\Subsection{D\'eterminants g\'en\'eralis\'es}\label{deligne}

Dans le reste de cet article, nous avons besoin de la construction de d\'eterminants sur des anneaux tels que $\Zp[G]$ ou $\Qp[G]$, pour un groupe ab\'elien fini $G$, ou encore $\Zp[[X]]$ et $\Qp \otimes_{\Zp} \Zp[[X]]$. Nous commen\c{c}ons donc par quelques rappels, tir\'es de \cite{KM, D2, BF2}, sur le formalisme tr\`es g\'en\'eral des d\'eterminants. 

Soit $A$ un anneau commutatif unitaire. On note $\mathbf{M}(A)$ la cat\'egorie des $A$-modules et $\mathbf{P}(A)$ la sous-cat\'egorie de $\mathbf{M}(A)$ form\'ee des modules projectifs de type fini.

On appelle cat\'egorie de Picard une cat\'egorie $\mathscr{P}$ dont toute 
fl\`eche est un isomorphisme, munie d'un foncteur $\boxtimes : \mathscr{P}
\times \mathscr{P} \to \mathscr{P}$ et d'une contrainte d'associativit\'e
pour $\boxtimes$. On peut d\'eduire de ces axiomes l'existence d'un objet unit\'e 
$\mathbf{1}_{\mathscr{P}}$, unique \`a isomorphisme pr\`es. Tout objet $X$ de 
$\mathscr{P}$ admet un inverse $X^{-1}$ tel que $X \boxtimes X^{-1} \simeq
\mathbf{1}_{\mathscr{P}}$. On dit qu'une cat\'egorie de Picard $\mathscr{P}$ est 
commutative si elle est munie d'une contrainte de commutativit\'e compatible 
\`a la contrainte d'associativit\'e. 

Soit $(\mathbf{P}(A),\mathrm{is})$ la cat\'egorie dont les objets sont ceux 
de $\mathbf{P}(A)$ et dont les fl\`eches sont les isomorphismes. On appelle foncteur d\'eterminant un foncteur $\det : (\mathbf{P}(A),\mathrm{is}) \to 
\mathscr{P}$ v\'erifiant les propri\'et\'es suivantes :
\begin{enumerate}
\item Pour toute suite exacte $0 \to P' \to P \to P'' \to 0$, on a un isomorphisme fonctoriel :
$\det(P) \simeq \det(P') \boxtimes \det(P'')$.
\item Pour toute suite exacte $0 \to P \overset{\alpha}{\to} Q \to 0$, l'application
$\det(\alpha)$ co\"{\i}ncide avec le compos\'e :
\[ \det(P) \simeq \det(0) \boxtimes \det(Q) \simeq \det(Q), \] et de m\^eme 
$\det(\alpha)^{-1}$ co\"{\i}ncide avec le compos\'e :
\[ \det(Q) \simeq \det(P) \boxtimes \det(0) \simeq \det(P). \] 
\item Si $P=P' \oplus P''$ et si :$0 \to P' \to P \to P'' \to 0$ et $0 \to P'' \to P \to P' \to 0$ sont les suites exactes naturelles, alors le diagramme :
\[ \xymatrix{& \det(P) \ar[dl] \ar[dr] & \\
\det(P') \boxtimes \det(P'') \ar[rr] & &  \det(P'') \boxtimes \det(P')} \]
est commutatif. 
\item Pour tout module projectif $P$ muni d'une filtration $P \supset P '\supset  P'' \supset \{0\}$, le diagramme : 
\[ \begin{CD}
\det(P) @>>> \det(P') \boxtimes \det(P/P') \\
@VVV @VVV \\
\det(P'') \boxtimes \det(P/P'') @>>> \det(P'') \boxtimes 
\det(P'/P'') \boxtimes \det(P/P') 
\end{CD} \]
est commutatif.
\end{enumerate}

Soit $\mathbf{K}(A)=\mathbf{K}(\mathbf{M}(A))$ la cat\'egorie 
des complexes de $A$-modules. On dit qu'un morphisme de complexes 
$f : M^\bullet \to N^\bullet$ 
est un quasi-isomorphisme si pour tout $i$, l'application
$H^i(M^{\bullet}) \to H^i(N^{\bullet})$ est un isomorphisme.
La cat\'egorie d\'eriv\'ee $\mathbf{D}(A)=\mathbf{D}(\mathbf{K}(A))$ est la localisation de $\mathbf{K}(A)$ par rapport aux quasi-isomorphismes.

On dit qu'un objet $M^{\bullet}$ de $\mathbf{D}(A)$ est parfait
s'il existe un complexe born\'e de $A$ modules projectifs de type
fini : $P^{\bullet} = ( \cdots \to P_{i+1} \to P_i \to P_{i-1} \to \cdots)$
quasi-isomorphe \`a $M^{\bullet}$. 
Soit $\mathbf{D}^{\mathrm{p}}(A)$ la sous-cat\'egorie 
de $\mathbf{D}(A)$ form\'ee
des objets parfaits. Pour tout objet $M^{\bullet}$ de 
$\mathbf{D}^{\mathrm{p}}(A)$, on fixe un complexe $P^{\bullet}$ v\'erifiant les conditions ci-dessus et l'on pose :
\[ \det(M^{\bullet}) = \boxtimes_{i \in \ZZ} \det(P_i)^{(-1)^i}. \]
On obtient ainsi une extension du foncteur $\det$, unique \`a \'equivalence pr\`es,
\`a un foncteur (encore not\'e $\det$) :
\[ \det : (\mathbf{D}^{\mathrm{p}}(A),\mathrm{qis}) \to \mathscr{P}. \]
Si les modules de cohomologie $H^i(M^{\bullet})$ sont parfaits
en toutes dimensions, on a alors un isomorphisme fonctoriel :
\[ \det(M^{\bullet}) \simeq \boxtimes_{i \in \ZZ} 
\det(H^i(M^{\bullet}))^{(-1)^i}. \]

Dans \cite{D2}, Deligne construit une cat\'egorie de Picard
commutative $\mathscr{V}(A)$ et un foncteur d\'eterminant universel
$[\cdot]_A : (\mathbf{P}(A),\mathrm{is}) \to \mathscr{V}(A)$
tel que tout foncteur d\'eterminant $\det$ s'\'ecrit comme le compos\'e
de $[\cdot]_A$ avec un foncteur additif $\mathscr{V}(A) \to \mathscr{P}$.
On en d\'eduit en particulier un foncteur 
$[\cdot]_A : (\mathbf{D}^{\mathrm{p}}(A),\mathrm{qis}) \to \mathscr{V}(A)$. La proposition ci-dessous rassemble quelques propri\'et\'es du foncteur $[\cdot]_A$.

\begin{prop}\label{detuniv}
Si $f : A \to B$ est un morphisme d'anneaux, alors le foncteur {\og extension des scalaires \fg} $f^*M = B \otimes_A M$ induit des foncteurs  
$\mathbf{L} f^* : \mathbf{D}^{\mathrm{p}}(A) \to \mathbf{D}^{\mathrm{p}}(B)$ et $f^* :   
\mathscr{V}(A) \to \mathscr{V}(B)$, et les foncteurs $[\cdot]_B \circ \mathbf{L} f^*$
et $f^* \circ [\cdot]_A : (\mathbf{D}^{\mathrm{p}}(A),\mathrm{qis}) \to \mathscr{V}(B)$ 
sont quasi-isomorphes. 

Si on suppose de plus que $B$ est projectif de type fini sur $A$, alors la restriction
des scalaires induit des foncteurs $f_* : \mathbf{D}^{\mathrm{p}}(B) \to \mathbf{D}^{\mathrm{p}}(A)$ et $f_* : \mathscr{V}(B) \to \mathscr{V}(A)$, et les foncteurs $f_* \circ [\cdot]_B$ et $[\cdot]_A \circ f_* : (\mathbf{D}^{\mathrm{p}}(B),\mathrm{qis}) \to \mathscr{V}(A)$ sont quasi-isomorphes. 
\end{prop}

\begin{proof}
Voir \cite[section 4.11]{D2}.
\end{proof}

On note $\mathscr{P}(A)$ la cat\'egorie des $A$-modules inversibles
gradu\'es. Un objet de $\mathscr{P}(A)$ s'identifie \`a une paire
$(X,\alpha)$ o\`u $X$ est un $A$-module inversible et $\alpha : \mathrm{Spec} (A) \to \ZZ$ est une fonction localement
constante. Une fl\`eche $f : (X,\alpha) \to (Y,\beta)$ n'existe que si $\alpha = \beta$, auquel cas c'est un isomorphisme. On munit $\mathscr{P}(A)$ d'un produit tensoriel en posant
\[ (X,\alpha) \otimes (Y,\beta) = (X \otimes_A Y, \alpha + \beta). \]
Munie de la contrainte de
commutativit\'e donn\'ee par la r\`egle de Koszul :
\begin{align*} 
\psi : X \otimes_A Y & \to Y \otimes_A X \\
\psi(x \otimes y) & = (-1)^{\alpha \beta} y \otimes x,
\end{align*}
la cat\'egorie $\mathscr{P}(A)$ est alors une cat\'egorie de Picard commutative.
On identifie l'oppos\'e $(X,\alpha)^{-1}$ d'un \'el\'ement
$(X,\alpha)$ \`a $(X^*,-\alpha)$ o\`u $X^* = \mathrm{Hom}_A (X,A)$.

Si $P$ est un $A$-module projectif de type fini, alors
le rang de $P$  est une fonction localement constante
$\mathrm{rg}_P : \mathrm{Spec} (A) \to \ZZ$ et on d\'efinit
le d\'eterminant de Knudsen-Mumford $\det_A (P)$ en posant :
\[ \det_A(P) = ( \wedge^r P, \mathrm{rg}_P ) \in \mathrm{Ob}(\mathscr{P}(A)). \]

Remarquons que la propri\'et\'e universelle du foncteur $[\cdot]_A$ donne un foncteur
additif $\mathbf{V}(A) \to \mathscr{P}(A)$ qui n'est pas, en g\'en\'eral, une \'equivalence de cat\'egories.

Dans cet article nous n'utilisons que les d\'eterminants sur des produits finis d'anneaux locaux.  Dans ce cas les cat\'egories $\mathbf{V}(A)$ et $\mathscr{P}(A)$ sont \'equivalentes par \cite[section 4.13]{D2} et la construction de Knudsen-Mumford fournit donc un foncteur d\'eterminant universel. Les anneaux typiques auquels nous allons appliquer la th\'eorie pr\'ecedente sont :
\begin{enumerate}
  \item $A=\Qp[G]$ ou bien $A=\Zp[G]$, o\`u $G$ est un groupe ab\'elien fini.
  \item $A = \Zp[[X]]$ ou bien $A = \Qp \otimes_{\Zp} \Zp[[X]]$.
\end{enumerate}

\begin{prop}\label{reguresol}
Si $A$ est un anneau local r\'egulier de dimension $n$, alors : 
\begin{enumerate}
  \item Tout $A$-module de type fini $M$ admet une r\'esolution
projective $P^\bullet$ : $ 0 \to P_m \to \cdots \to P_0 \to M \to 0$ avec $m \leq n$. 
  \item Si $Q(A)$ est l'anneau total des fractions de $A$, et si $M$ est
un $A$-module de torsion, alors le produit tensoriel de $P^{\bullet}$ par
$Q(A)$ donne une suite exacte $ 0 \to Q(A) \otimes_A P_m \to \cdots \to Q(A) \otimes_A P_0 \to 0$. On en d\'eduit une injection canonique $i_A : \det_A(M) \to \det_{Q(A)}(Q(A) \otimes_A P^\bullet) \simeq Q(A)$ et l'image de $\det_A(M)$ dans $Q(A)$ ne d\'epend pas du choix de $P^{\bullet}$ et co\"{\i}ncide avec l'id\'eal
fractionnaire de $M$.
\end{enumerate} 
\end{prop}

\begin{proof}
La premi\`ere assertion est un th\'eor\`eme classique de Serre (voir par exemple \cite[\S 19]{M}). Pour la  deuxi\`eme voir \cite[th\'eor\`eme 3]{KM}.
\end{proof}

\begin{exem}\label{detiw}
En particulier, consid\'erons l'anneau $A=\Zp[[X]]$ qui est local r\'egulier de dimension $2$ et soit $M$ un $A$-module de type fini et de torsion. Il existe alors une suite exacte :
\[ 0 \to (\mathrm{fini}) \to M \to \oplus_{i=1}^{n} A / f_i A \to (\mathrm{fini}) \to 0, \] o\`u les $f_i$ sont des polyn\^omes distingu\'es. On a alors $\det_A(M) =  \mathrm{car}_A(M)^{-1} A$ o\`u $\mathrm{car}_A(M) = \prod_{i=1}^n f_i$ est le polyn\^ome caract\'eristique de $M$.
\end{exem}

\begin{rema} 
L'approche de Deligne a \'et\'e g\'en\'eralis\'ee aux anneaux non-com\-mutatifs par Burns et Flach, voir \cite{BF2}.
\end{rema}

\Subsection{D\'eterminants de la cohomologie galoisienne}\label{detcoh}

On note $\mathbf{M}(G_K)$ la cat\'egorie des $\Zp$-repr\'esentations de $G_K$, c'est-\`a-dire la cat\'egorie des $\Zp$-modules (pas n\'ecessairement de type fini) munis d'une action lin\'eaire et continue de $G_K$. Si $L$ est une extension finie de $K$ on a les foncteurs habituels :
$\Res_{L/K} : \mathbf{M}(G_K) \to \mathbf{M}(G_L)$ et $\Ind_{L/K} : \mathbf{M}(G_L) \to \mathbf{M}(G_K)$, ce dernier foncteur \'etant donn\'e par la formule $\Ind_{L/K} M = \Zp[G_K] \otimes_{\Zp[G_L]} M$. 

Supposons maintenant $L$ ab\'elien sur $K$ et posons $G=\mathrm{Gal}(L/K)$. On note $\iota : \Zp[G] \to \Zp[G]$ l'involution $g \mapsto g^{-1}$. Si $M \in \mathbf{M}(G_K)$, on note (pour simplifier) $\Ind_{L/K} M$ le $\Zp[G_K]$-module $\Ind_{L/K}(\Res_{L/K} M)$. Le module $\Ind_{L/K} M$ a alors une structure naturelle de $\Zp[G]$-module donn\'ee par la formule $\overline{g} (\sigma \otimes m)
= \sigma g^{-1} \otimes g (m)$, et on a des isomorphismes canoniques :
\begin{align*}
M^{G_L} \simeq (\Ind_{L/K} M)^{G_K}, &  \quad m \mapsto \sum_{\overline{g} \in G} g \otimes m; \\
\Ind_{L/K}(M) \simeq (\Zp[G] \otimes_{\Zp} M)^{\iota}, & \quad \sigma \otimes m \mapsto \overline{\sigma} \otimes \sigma (m). 
\end{align*}

Soit $\mathbf{M}(G_K)^\mathrm{ind}$ la sous-cat\'egorie de $\mathbf{M}(G_K)$
dont les objets sont les limites inductives de $\Zp[G_K]$-modules de type fini sur $\Zp$. Pour tout $M \in \mathbf{M}(G_K)^{\mathrm{ind}}$, on note $C^\bullet(G_K, \Ind_{L/K} M)$ le complexe des cocha\^{\i}nes continues de $G_K$  \`a valeurs dans $\Ind_{L/K}(M)$. On obtient ainsi un foncteur de $\mathbf{M}(G_K)^{\mathrm{ind}}$ dans $\mathbf{D}(\Zp[G])$ qui \`a $M$ associe $C^\bullet(G_K, \Ind_{L/K} M)$ et qui induit un foncteur exact :
\[ \RR\Gamma(L,\cdot) : \mathbf{D}(\mathbf{M}(G_K)^{\mathrm{ind}}) \to  \mathbf{D}(\Zp[G]). \]
Le lemme de Shapiro donne un isomorphisme $\RR^i \Gamma(L,M) \simeq H^i(L,M)$.

\begin{prop}\label{pschap}
Si $L/K$ est une extension ab\'elienne finie et si $M$ est un $\Zp[G_K]$-module qui est de type fini sur $\Zp$, alors :
\begin{enumerate}
  \item $\RR\Gamma(L,M) \in \mathbf{D}^p(\Zp[G])$;
  \item si de plus $M$ est de $\Zp$-torsion, alors $\det_{\Zp[G]} \RR\Gamma(L,M) = \det_{\Zp[G]}^{-1} (\Ind_{L/\Qp} M)$ dans $Q(\Zp[G])$.
\end{enumerate}
\end{prop}

\begin{proof}
Voir \cite{K3} et \cite{BF1}.
\end{proof}

Pour terminer ce paragraphe, faisons le lien entre les constructions ci-dessus et la th\'eorie d'Iwasawa des repr\'esentations $p$-adiques. Rappelons que l'on a fix\'e un syst\`eme 
compatible $(\zn)_{n \geq 0}$ de racines primitives $p^n$-i\`emes de l'unit\'e et pos\'e $K_n=K(\zn)$ et $K_\infty = \cup_{n \geq 1} K_n$. Soient $G_n = \mathrm{Gal}(K_n/K)$, $H_K= \mathrm{Gal}(\Qpbar/K_\infty)$, $\Gamma = \mathrm{Gal}(K_\infty/K)$, $\Gamma_n = \mathrm{Gal}(K_\infty/K_n)$ et enfin $\Lambda = \Zp[[\Gamma]]$ l'alg\`ebre d'Iwasawa de $\Gamma$. 

Si $T$ est une $\Zp$-repr\'esentation de $G_K$, alors le module induit $\Ind_{K_\infty/K}(T)$ est isomorphe \`a $(\Lambda \otimes_{\Zp} T)^\iota$ et on pose $\RR\Gamma_{\mathrm{Iw}}(K,T) = \RR\Gamma(K,\Ind_{K_\infty/K}(T))$. La proposition suivante est un cas particulier d'un r\'esultat de Nekov\'a\v{r} (voir \cite[prop 8.4.22]{N}).

\begin{prop}\label{old125}
\begin{enumerate}
  \item On a  des isomorphismes canoniques $\RR^i\Gamma_{\mathrm{Iw}}(K,T) \simeq 
H^i_{\mathrm{Iw}}(K,T)$.
  \item Dans la cat\'egorie $\mathbf{D}(\Zp[G_n])$, on a un isomorphisme canonique : 
  \[ \Zp[G_n] \otimes_\Lambda^{\mathbf{L}}  \RR\Gamma_{\mathrm{Iw}}(K,T) \simeq  \RR\Gamma(K_n,T). \]
  \item On a une suite spectrale d\'eg\'en\'er\'ee : 
  \[  E_2^{ij} = H^i(K_{\infty}/K_n, H_{\mathrm{Iw}}^j(K,T)) \Rightarrow H^{i+j-1} (K_n,T), \] qui donne lieu \`a des suites exactes : 
  \[ 0 \to  H_{\mathrm{Iw}}^j(K,T)_{\Gamma_n} \to H^j(K_n,T) \to H_{\mathrm{Iw}}^{j+1}(K,T)^{\Gamma_n} \to 0. \]
\end{enumerate}
\end{prop}

Enfin, remarquons que la dualit\'e locale fournit un isomorphisme $H_{\mathrm{Iw}}^2(K,T) \simeq H^0(K_\infty,V^*(1)/T^*(1))^{\wedge}$ o\`u ${}^\wedge$ signifie le dual de Pontryagin. En particulier, si $(V^*)^{H_K} = 0$, alors $H_{\mathrm{Iw}}^2(K,T)$ est fini et on a $\sharp H_{\mathrm{Iw}}^2(K,T)^{\Gamma_n} = \sharp H^0 (K_n ,V^*(1)/T^*(1))$.

\Subsection{Constantes locales des repr\'esentations de Weil-Deligne}\label{langlands}

L'objet de ce paragraphe est de fournir des rappels sur la th\'eorie des constantes locales, telle qu'elle est d\'evelopp\'ee dans \cite{D1}, auquel nous renvoyons pour plus de d\'etails. Le corps $K$ est toujours une extension finie de $\Qp$. On fixe une uniformisante $\pi_K$ de $K$ et on note $|\cdot|_K$ la norme de $K$ normalis\'ee par $|\pi_K|_K = q_K^{-1}$ o\`u $q_K$ est le cardinal du corps r\'esiduel $k_K$ de $K$.

On note $\Knr$ l'extension maximale non-ramifi\'ee de $K$ et $\Fr_K$ le Frobenius g\'eom\'etrique de $\Knr$. Le groupe de Weil $W_K$ de $K$ est par d\'efinition le sous-groupe de $G_K$ form\'e des $g \in G_K$ tels que la restriction de $g$ \`a $\Knr$ soit une puissance enti\`ere de $\Fr_K$. On a donc une suite exacte : 
\[ 0 \to I_K \to W_K \overset{\nu}{\to} \ZZ \to 0, \]
o\`u l'application $\nu$ est d\'efinie par la formule $w_{\mid\Knr} = \Fr_K^{\nu (w)}$.

Soit $E$ un corps de caract\'eristique $0$ et contenant toutes les racines de l'unit\'e d'ordre une puissance de $p$ et d'ordre $p-1$. On fixe une mesure de Haar $\mu_K$ sur $K$ et un caract\`ere additif continu $\psi : K \to E^\times$ (le corps $E$ \'etant muni de la topologie discr\`ete). Comme $\psi$ est continu, il est trivial sur un sous-groupe ouvert de $K$ et l'on d\'efinit son conducteur $n(\psi)$ comme \'etant le plus grand entier $n$ tel que $\psi$ est trivial sur $\pi_K^{-n} \OO_K$.

La th\'eorie de Langlands et Deligne (voir \cite{D1}) associe \`a toute repr\'esentation $E$-lin\'eaire $V$ de $W_K$ une constante $\eps(V,\psi,\mu_K)$ v\'erifiant les propri\'et\'es suivantes :
\begin{enumerate}
  \item Si $V$ est de dimension $1$, alors $\eps(V,\psi,\mu_K)$ co\"{\i}ncide avec la constante locale {\og ab\'elienne \fg} d\'efinie par la th\'eorie de Tate (dans \cite{T}). Plus pr\'ecisement, l'isomorphisme de r\'eciprocit\'e $K^\times \to W_K^{\mathrm{ab}}$ permet de voir $V$ comme un quasi-caract\`ere $\eta : K^\times \to E^\times$. On note $a(\eta)$ le conducteur de $\eta$ et on fixe $c \in \OO_K$ v\'erifiant $v_K (c) = a(\eta) + n(\psi)$.
Si $\eta$ est non-ramifi\'e, alors on a : 
\[ \eps(\eta,\psi,\mu_K) = \frac{\eta (c)}{|c|_K} \int_{\OO_K} d\mu_K, \]
et si $\eta$ est ramifi\'e, alors on a :
\[ \eps(\eta,\psi,\mu_K) = \sum_{n \in \ZZ} \int_{\{v_K(x)=n\}} \eta^{-1}(x) \psi(x) d\mu_K = \int_{c^{-1}\OO_K} \eta^{-1}(x) \psi(x) d\mu_K. \]
  \item Pour toute suite exacte de repr\'esentations $0 \to V' \to V \to V'' \to 0$, on a $\eps(V,\psi,\mu_K) = \eps(V',\psi,\mu_K)\eps(V'',\psi,\mu_K)$.
  \item Pour tout $a \in K^\times$, on a $\eps(V,\psi, a\mu_K) = a^{\dim V}\eps(V,\psi,\mu_K)$ et si $m_a$ d\'enote la fonction $x \mapsto ax$, alors $\eps(V,\psi \circ m_a,\mu_K) = \det (V)(a) |a|_K^{-\dim V} \eps(V,\psi,\mu_K)$.
 \item Si $L$ est une extension finie de $K$, alors il existe une constante $\lambda (L/K,\psi, \mu_L,\mu_K)\in E$ telle que pour toute repr\'esentation $V$ de $W_L$ on ait : 
\[ \eps(\Ind_{L/K} (V),\psi, \mu_K) =
\lambda (L/K,\psi,\mu_L,\mu_K)^{\dim V} \eps(V,\psi\circ\mathrm{Tr}_{L/K},\mu_L). \]
 \item Soient $\omega_1 : K^\times \to E^\times$ le quasi-caract\`ere donn\'e par la formule $\omega_1 (a) =|a|_K$ et $\mu^*_K$ la mesure duale de $\mu_K$ relativement  \`a $\psi$. On a alors : 
\[ \eps(V,\psi,\mu_K) \eps(V^* \otimes \omega_1,\psi \circ m_{-1},\mu_K^*) = 1. \]
 \item Pour une repr\'esentation non-ramifi\'ee $W$, on a :
\[ \eps(V \otimes W,\psi,\mu_K) = \det(W) (\pi_K^{a(V)+\dim (V)n(\psi)})
\eps(V,\psi,\mu_K)^{\dim W}, \]
o\`u $a(V)$ est le conducteur d'Artin de $V$.
\end{enumerate}

Rappelons que l'on a fix\'e un syst\`eme compatible $(\zn)_{n \geq 0}$ de racines de l'unit\'e. On note $\psi_0$ l'unique caract\`ere additif de $\Qp$ v\'erifiant $\psi_0(1/p^n)=\zn$ et on pose $\psi_K = \psi_0 \circ \mathrm{Tr}_{K/\Qp}$. On normalise la mesure $\mu_K$ en imposant $\mu_K(\OO_K)=1$. Soit enfin $(\cdot,\cdot)_K :K^\times \times K^\times \to \{\pm1\}$ le symbole de Hilbert. Le lemme suivant est bien connu des experts.

\begin{lemm}\label{lambda}
Si $L/K$ est une extension finie, alors :
\[ \lambda(L/K,\psi_K,\mu_L,\mu_K) = \pm (-1,d_{L/K})_K^{1/2} |d_{L/K}|_p^{[K:\Qp] / 2}, \]
o\`u $d_{L/K}$ est le discriminant de $L/K$. 
\end{lemm}

\begin{proof} 
Les formules (3) et (5), appliqu\'ees \`a la repr\'esentation 
r\'eguli\`ere $\Ind_{L/K}[1]$ donnent : 
\[ \eps(\Ind_{L/K}[1],\psi_K,\mu_K) \eps(\Ind_{L/K}[1] \otimes \omega_1, \psi_K \circ m_{-1}, \mu_K)
=  | d_K|_p^{-[L:K]}. \]
Comme par ailleurs $a(\Ind_{L/K}[1]) = v_K(d_{L/K})$ (voir par exemple \cite[chap IV, prop. 4]{S1}) et $n(\psi_K)=v_K(\mathscr{D}_{K/\Qp})$, on a :
\[ \eps(\Ind_{L/K}[1] \otimes \omega_1, \psi_K \circ m_{-1}, \mu_K) = |d_L|_p \det (\Ind_{L/K}[1])(-1) 
\eps(\Ind_{L/K}[1],\psi_K,\mu_K). \]  
On a $|d_L|_p =  | d_{L/K}|_p^{[K:\Qp]} | d_{K}|_p^{[L:K]}$ et $\det (\Ind_{L/K}[1])(-1) = (-1,d_{L/K})_K$, d'o\`u :
\[ \eps(\Ind_{L/K}[1],\psi_K,\mu_K) = \pm(-1,d_{L/K})_K^{1/2} | d_{L/K}|_p^{[K:\Qp]/2}
| d_L|_p^{-1}. \]
Comme $\eps([1],\psi_L,\mu_L) = | d_L|_p^{-1}$, on en d\'eduit le lemme.
\end{proof}

\begin{rema} 
Il est facile de voir que si $L/K$ est une extension non-ramifi\'ee de degr\'e $f$, alors $\lambda (L/K,\psi,\mu_L,\mu_K) = (-1)^{(f-1)n(\psi)}$. 
\end{rema}

Supposons maintenant que $K$ est une extension non-ramifi\'ee de $\Qp$ de degr\'e $f$, et notons $X(G_n)$ le groupe des caract\`eres de $G_n=\mathrm{Gal}(K_n/K)$ \`a valeurs dans $E$. Pour tout $\eta \in X(G_n)$, on note $e_\eta$ l'idempotent habituel
\[ e_{\eta} =  \frac{1}{\sharp G_n} \sum_{g\in G_n} \eta^{-1}(g)g. \]
On d\'efinit la somme de Gauss $\tau(\eta)$ en posant $\tau(\eta) = 
\sum_{g \in G_k} \eta^{-1}(g) g(\zeta_{p^k}) = \sharp G_k e_{\eta} (\zeta_{p^k})$, o\`u $k=a(\eta)$ est le conducteur de $\eta$.

\begin{lemm}\label{old136}
Pour tout caract\`ere $\eta \in X(G_n)$, on a :
\[ \eps(\eta,\psi_K,\mu_K) = (-1)^{(f-1)a(\eta)} \tau(\eta)^f. \]
\end{lemm}

\begin{proof} 
Comme $K/\Qp$ est non-ramifi\'ee, les groupes de Galois des extensions $\Qp(\zn) /\Qp$ et $K_n/K$ sont isomorphes et $\eta$ peut \^etre vu comme la restriction $\Res_{K/\Qp} \tilde{\eta}$ d'un caract\`ere $\tilde{\eta} : \mathrm{Gal}(\Qp(\zn) / \Qp) \to E^\times$. Comme $\lambda (K/\Qp,\psi_0,\mu_K,\mu_0)=1$, on a :
\[ \eps(\eta,\psi_K,\mu_K) = \eps(\Ind_{K/\Qp} \eta,\psi_0,\mu_0)
\eps(\Ind_{K/\Qp} [1] \otimes  \tilde{\eta},\psi_0,\mu_0) = (-1)^{(f-1) a(\eta)}\eps(\tilde{\eta},\psi_0,\mu_0)^f. \]
Si on suppose que $n=a(\eta)$, alors l'application compos\'ee
$\Qp^\times \to \mathrm{Gal}(\Qp^\mathrm{ab}/\Qp) \to G_n \simeq (\ZZ/p^n\ZZ)^\times$ envoie $u \in \Zp^\times$ sur $u \mod{p^n}$ et $p$ sur $1$, ce qui fait que :
\[ \eps(\tilde{\eta},\psi_0,\mu_0) = p^n \sum_{u \in \Zp^\times / 1+p^n \Zp} \mu_0(1+p^k \Zp)  \tilde{\eta} (u)^{-1} \zn^u = \sum_{u \in \Zp^\times / 1+p^n \Zp} \tilde{\eta}(u)^{-1} \zn^u = \tau(\eta). \]
Le cas g\'en\'eral s'en d\'eduit.
\end{proof}

On appelle repr\'esentation du groupe de Weil-Deligne un couple $(\rho, N)$ form\'e d'une repr\'esentation $\rho : W_K \to \mathrm{Aut}_E(V)$ du groupe de Weil $W_K$
et d'un endomorphisme nilpotent $N : V \to V$ v\'erifiant $\rho(w)^{-1} N \rho(w) = q_K^{\nu (w)}N$ (voir \cite[\S 8]{D1}). On pose alors : 
\[ \eps(V,\psi_K,\mu_K) = \eps(\rho, \psi_K, \mu_K) \det (-\Fr_K \mid V^{I_K}/(V^{I_K})^{N=0}). \]

\Subsection{Constantes locales des repr\'esentations potentiellement semi-stables}\label{fontaine}

Pour plus de d\'etails, voir \cite[chapitre I, \S 1.3]{FP}. On garde les notations et les conventions des paragraphes pr\'ec\'edents. En particulier, $K$ est toujours une extension finie de $\Qp$ et $K_0$ est son sous-corps maximal non-ramifi\'e, dont le degr\'e sur $\Qp$ est $f=[K_0:\Qp]$. Rappelons que l'on a d\'efini ci-dessus un caract\`ere additif $\psi_K$ \`a valeurs dans $\Qp(\zeta_{p^\infty}) = \cup_{n \geq 0} \Qp(\zn)$ en posant $\psi_K(a/p^n)=\zn^{\mathrm{Tr}_{K/\Qp}(a)}$. On fixe une extension ab\'elienne finie $L/K$ et on pose toujours $G=\mathrm{Gal}(L/K)$.

Le lemme suivant est laiss\'e en exercice au lecteur.

\begin{lemm}\label{ddrind}
Si $L/K$ est une extension finie et si $V$ est une repr\'esentation $p$-adique de $G_L$, alors $\ddR^K(\Ind_{L/K}V) \simeq \ddR^L(V)$ et $\dpst(\Ind_{L/K}V) \simeq \Ind_{L/K} \dpst(V)$.
\end{lemm}

Si $V$ est une repr\'esentation potentiellement semi-stable de $G_K$, alors la repr\'esentation $\Ind_{L/K}V \simeq (\Qp [G] \otimes_{\Qp} V)^{\iota}$ est bien-s\^ur elle aussi potentiellement semi-stable et $D = \dpst(\Ind_{L/K}V)$ est un $K_0^{\mathrm{nr}} [G]$-module muni d'une action naturelle $\Fr_K$-semi-lin\'eaire de $W_K$. On munit $D$ d'une action \emph{lin\'eaire} de $W_K$, 
$\rho : W_K \to \mathrm{Aut}_{K_0^{\mathrm{nr}}[G]}(D)$ en posant
$(\rho (w))(d) = w\varphi^{f\nu (w)}(d)$ o\`u l'application $\nu$ est celle d\'efinie au paragraphe \ref{deligne}. Le module $D$ est muni d'un op\'erateur de monodromie $N$ v\'erifiant $N \circ \varphi =p \varphi \circ N$ ce qui fait que $\rho (w)^{-1} N \rho(w) = q_K^{\nu (w)}$ et que $(\rho, N)$ est une repr\'esentation du groupe de Weil-Deligne. On pose alors :
\[ \eps(L/K, V) = \eps(D,\psi_K,\mu_K) = \eps(\rho, \psi_K,\mu_K)
\det(-\Fr_K \mid D^{I_K}/(D^{I_K})^{N=0}). \]
Il est facile de voir (cf \cite[remarque 1.3.3]{FP}) que la repr\'esentation $\rho$ est $\Qp$-rationnelle, d'o\`u l'on tire que $\eps(L/K,V) \in \Qp(\zeta_{p^\infty})[G]$.

Si $E$ est un corps contenant $K^{\mathrm{nr}}$ ainsi que les valeurs des caract\`eres de $G$, alors on a : 
\[ E[G] = \oplus_{\eta \in X (G)} E_{\eta},\quad\text{o\`u $E_{\eta}=e_{\eta}E$,} \]
et le module $D$ se d\'ecompose sur $E$ en produit de ses $\eta$-composantes : \[ E \otimes_{K_0^{\mathrm{nr}}} D  = \oplus_{\eta \in X (G)} D_{\eta},
\quad\text{o\`u $D_{\eta}=e_{\eta}(E \otimes_{K_0^{\mathrm{nr}}} D)$.} \]
On appelle $\eta_0$ le caract\`ere trivial. On d\'eduit de la d\'ecomposition ci-dessus que $\eps(L/K,V) = \sum_{\eta \in X (G)}  \eps(D_{\eta},\psi_{K,\eta},\mu_K)$, avec $\eps(D_{\eta},\psi_{K,\eta},\mu_K)
\in E_\eta$.

Si $V$ est une repr\'esentation potentiellement semi-stable de $G_K$, alors par le lemme \ref{ddrind} ci-dessus, $\ddR^{\Qp}(\Ind_{L/\Qp}(V)) \simeq \ddR^L(V)$ et on a donc un isomorphisme canonique :
\[ \mathrm{comp}_{L/\Qp} : \bdR \otimes_{\Qp} \Ind_{L/\Qp} V  \simeq \bdR \otimes_{\Qp} \ddR^L(V). \]
On en d\'eduit un homomorphisme :
\[ \tilde{\alpha}_{V,L/K}: \det^{-1}_{\Qp[G]}(\ddR^L(V)) \otimes \det_{\Qp[G]} (\Ind_{L/\Qp}V) \to \Qp[G] \otimes_{\Qp} \bdR. \]

On voit que $\det_{\Qp[G]} (\Ind_{L/\Qp}V)$ est une $\Qp[G]$-repr\'esentation de de Rham de rang $1$ et de poids $r=-[K:\Qp]t_H(V)$ et il existe donc une extension ab\'elienne finie $K'/\Qp^{\mathrm{nr}}$ telle que la restriction de  $\det_{\Qp[G]} (\Ind_{L/\Qp}V)$ \`a $G_{K'}$ soit isomorphe \`a $\Qp[G](r)$. On en d\'eduit donc une application $\alpha_{V,L/K} : \det^{-1}_{\Qp[G]} (\ddR^L(V)) \otimes \det_{\Qp[G]} (\Ind_{L/\Qp}V) \to K'[G]$, donn\'ee par la formule $\alpha_{V,L/K} = t^{-r} \tilde{\alpha}_{V,L/K}$, o\`u $t = \log[\eps] \in \bdR$ est l'uniformisante de $\bdR^+$ associ\'ee \`a $\eps=(\zn)_{n \geq 0}$.

Soit $\hat{\sigma}$ l'\'el\'ement de $\mathrm{Gal}(\Qp^{\mathrm{ab}}/\Qp)$ qui op\`ere trivialement sur les racines $p^n$-i\`emes de l'unit\'e et dont la restriction \`a $\Qp^{\mathrm{nr}}$ est \'egale \`a $\sigma$, et soit $a_{V,L/K} = \det_{\Qp[G]} (\Ind_{L/\Qp}V)(\hat{\sigma}) \in \Zp[G]^\times$. 

\begin{defi}\label{defzpvlk}
On pose : $\Zp[G]_{V,L/K} = \{ x \in \hat{\Zp^{\mathrm{nr}}}[G] \mid \sigma(x) = a_{V,L/K} x \}$ et $\Qp[G]_{V,L/K} = \Qp \otimes_{\Zp} \Zp[G]_{V,L/K}$.
\end{defi}

Le module $\Zp[G]_{V,L/K}$ est alors libre de rang $1$ sur $\Zp[G]$ (voir \cite{K2}). Posons
\[ \Gamma^*(i) =
\begin{cases} (i-1)!, & \text{si $i>0$} \\
\frac{(-1)^i}{(-i)!}, & \text{si $i\leqslant 0$,} 
\end{cases} \]
et $\Gamma^*(V) = \prod_{i \in \ZZ} \Gamma^*(-i)^{h_i(V)[K:\Qp]}$. Soit aussi : 
\[ \beta_{V,L/K} = \lambda (K/\Qp)^{-\dim V} \Gamma^*(V)
\eps(L/K,V)^{-1} \alpha_{V,L/K}, \]
o\`u $\lambda (K/\Qp) = \lambda (K/\Qp,\psi_0,\mu_K,\mu_0)$ est la constante d\'efinie dans le paragraphe \ref{deligne}. 

\begin{lemm}
L'application $\beta_{V,L/K}$ induit un isomorphisme :
\[ \beta_{V,L/K} : \det^{-1}_{\Qp[G]} (\ddR^L(V)) \otimes \det_{\Qp[G]} (\Ind_{L/\Qp}V) \to {\Qp[G]}_{L/K,V}. \]
\end{lemm}

\begin{proof}
On note $\chi : G_K \to \Zp^\times$ le caract\`ere cyclotomique. Si on pose $D = \dpst(\Ind_{L/K}(V)) = (\Qp[G] \otimes_{\Qp} \dpst(V))^{\iota}$, alors on a (voir le paragraphe \ref{deligne}) : 
\[ \frac{ \lambda(K/\Qp)^{\dim V} \eps(D,\psi_K,\mu_K)}
{\eps(\Ind_{K/\Qp}D,\psi_0,\mu_0)} \in \Qp[G]. \]
Pour tout $g \in G_{\Qp}$, on a :
\begin{align*}
g(\eps(\Ind_{K/\Qp}D,\psi_0,\mu_0))
&= \eps(\Ind_{K/\Qp}D,\chi (g)\psi_0,\mu_0)\\
&=\det_{\Qp[G]}(\Ind_{L/\Qp}\dpst(V)) (\chi (g))
\eps(\Ind_{K/\Qp}D,\psi_0,\mu_0),
\end{align*}
D'autre part, si $x \in \det^{-1}_{\Qp[G]}(\ddR^L(V)) \otimes 
\det_{\Qp[G]} (\Ind_{L/\Qp}V)$, alors :
\begin{align*}
g(\alpha_{V,L/K}(x)) & = \chi^{-r}(g) \det_{\Qp[G]} (\Ind_{L/\Qp}V) (g) 
\alpha_{V,L/K}(x) \\
& = \det_{\Qp[G]}(\Ind_{L/\Qp}\dpst(V))(g) \alpha_{V,L/K}(x).
\end{align*}
On en d\'eduit le lemme.
\end{proof}

On donne maintenant une formule explicite pour l'application $\beta_{V,K_n/K}$ pour les repr\'esentations absolument cristallines, formule qui est utilis\'ee dans la suite. On suppose donc que $K$ est non-ramifi\'ee, et on \'ecrit comme ci-dessus $f=[K:\Qp]$, $q_K=p^f$ et $d= \dim V$. 

\begin{lemm}\label{old331}
Si $V$ est une repr\'esentation cristalline de $G_K$, alors :
\[ \eps(D_{\eta},\psi_{K,\eta}, \mu_K) = \det(\varphi \mid \dcris (V))^{a(\eta)}
\tau (\eta^{-1})^{fd} \otimes e_{\eta}^{\iota}. \]
\end{lemm}

\begin{proof} 
Comme $D = \dpst (\Ind_{K_n/K}V) = (K^{\mathrm{nr}}[G_n] \otimes_K \dcris (V) )^{\iota}$, on a $D_{\eta} = e_{\eta}^{\iota}\Qp^{\mathrm{ab}} \otimes_K \dcris (V)$, ce qui fait que $G_K$ op\`ere sur $D_{\eta}$ par le caract\`ere
$\eta^{-1}$ et en utilisant le lemme \ref{old136}, on obtient :
\begin{align*}
\eps(D_{\eta},\psi_{\eta}, \mu_K) & = \eps(\eta^{-1},\psi_K,\mu_K)^d
\det (\rho (\Fr_K)^{a(\eta)} \mid \dcris(V)) \otimes e_{\eta}^{\iota} \\
& = \tau (\eta^{-1})^{fd} \det(\varphi \mid \dcris (V))^{a(\eta)} \otimes 
e_{\eta}^{\iota},
\end{align*}
et le lemme est d\'emontr\'e.
\end{proof}

Soit $x_n=\zeta_p+\zeta_{p^2}+\cdots +\zn$ et soit $R_n$ le $\OO_K[G_n]$-r\'eseau de $K_n$ engendr\'e par $x_n$.  On fixe un $\OO_K$-r\'eseau $M$ de $\dcris (V)$ et $T$ un r\'eseau de $V$ et on pose $M_n= R_n \otimes_{\OO_K} M$ ainsi que :
\[ \beta_{V,K_n/K}(M,T)=\beta_{V,K_n/K} (\det^{-1}_{\Zp[G_n]}(M_n)
\otimes \det_{\Zp[G_n]}(\Ind_{K_n/K} T)). \]

\begin{prop}\label{old333} 
Si $V$ est une repr\'esentation cristalline de $G_K$ et $T$ un r\'eseau de $V$,  alors :
\[ \beta_{V,K_n/K}(M,T)=c^{fd} \Gamma^* (V)
q_K^{-n d} \left( \sum_{\eta \neq 1} \det (\varphi \mid \dcris (V))^{-a(\eta)} 
e_{\eta}+(-1)^{fd}q_K^{d} e_{1} \right) \alpha_{V,K}(M,T), \]
o\`u $c$ est la conjugaison complexe $c: \zn \mapsto \zeta_{p^n}^{-1}$.
\end{prop}

\begin{proof} 
Soient $\mathrm{comp}_{\mathrm{cris}} : \bcris \otimes_{\Qp} V \iso \bcris \otimes_K \dcris(V)$ et  $\mathrm{comp}_{K_n/K} : \bdR \otimes_{\Qp} \Ind_{K_n/K}V \iso \bdR \otimes_K \ddR^{K_n}(V)$ les isomorphismes de comparaison. Alors $\mathrm{comp}_{K_n/K,V}$
s'\'ecrit comme le compos\'e :
\begin{align*}
\bdR \otimes_{\Qp} (\Qp[G_n] \otimes_{\Qp}V)^{\iota}   
& \overset{\mathrm{comp}_{\mathrm{cris}}}{\longrightarrow}
\bdR \otimes_K \left( (\Qp[G_n] \otimes_{\Qp} K_n)^{G_n} \otimes_K \dcris(V) \right)^\iota \\ & \iso \bdR \otimes_K \ddR^{K_n}(V). 
\end{align*}

L'isomorphisme $(\Qp[G_n] \otimes_{\Qp} K_n)^{G_n} \iso K_n$
envoie $e_{\eta}^{\iota}\otimes \sharp G_n  e_{\eta}(x_n)$ sur
$e_{\eta}(x_n)$, et donc :
\[ \mathrm{comp}_{K_n/K,V} (e_{\eta}^{\iota} \otimes v) =
(\mathrm{comp}_{\mathrm{cris}} (v)e_{\eta} (x_n)) \otimes
\frac{1}{\sharp G_n  e_{\eta }(x_n)}. \]

Si $k=a(\eta) \neq 0$, alors on a $e_{\eta}(x_n)=e_{\eta}(\zeta_{p^k})$, et la formule bien connue $\tau (\eta) \tau(\eta^{-1}) = \eta (c) p^{k}$ nous donne $1/(\sharp G_n  e_{\eta^{-1}} (x_n)) = p^{-n} \tau(\eta^{-1}) \eta (c)$.

Si au contraire $\eta =1$, alors on a $e_1 (x_n) = (1-p)^{-1}$ d'o\`u $1/(\sharp G_n  e_1(x_n)) = - p^{1-n}$.

La proposition r\'esulte maintenant de la d\'efinition de l'application $\beta_{V,K_n/K}$ et du lemme \ref{old331}.
\end{proof}

\Subsection{La conjecture $C_{\mathrm{EP}}(L/K,V)$}\label{perrinriou}

On commence ce paragraphe par la d\'efinition de la droite d'Euler-Poincar\'e. Rappelons que $V$ est une repr\'esentation $p$-adique de $G_K$, que $L$ est une extension ab\'elienne finie de $K$ et que l'on a pos\'e $G=\mathrm{Gal}(L/K)$.

\begin{defi}\label{eulpoindef}
La droite d'Euler-Poincar\'e $\Delta_{\mathrm{EP}}(L/K,V)$ de $V$ est d\'efinie par la formule suivante : 
\begin{align*}
\Delta_{\mathrm{EP}}(L/K,V) &= 
\det_{\Qp[G]} \RR\Gamma(L,V) \otimes \det_{\Qp[G]}  (\Ind_{L/\Qp} V) \\
& \simeq \otimes_{i=0}^2 (\det_{\Qp[G]} H^i(L,V))^{(-1)^i} \otimes \det_{\Qp[G]}  (\Ind_{L/\Qp} V).
\end{align*}
\end{defi}

Si $T$ est un $\Zp$-r\'eseau de $V$, alors $\Ind_{L/\Qp} T = \Ind_{K/\Qp} (\Zp[G] \otimes_{\Zp}T)^\iota$ et $\RR\Gamma(L,T)$ sont parfaits sur $\Zp[G]$, et par la proposition \ref{pschap} le sous-$\Zp[G]$-module de $\Delta_{\mathrm{EP}}(L/K,V)$ :
\[ \Delta_{\mathrm{EP}}(L/K,T) = \det_{\Zp[G]} \RR\Gamma(L,T) \otimes \det_{\Zp[G]}  (\Ind_{L/\Qp} T) \]
ne d\'epend pas du choix de $T$ et d\'efinit donc un $\Zp[G]$-r\'eseau canonique de
$\Delta_{\mathrm{EP}}(L/K,V)$. 

Revenons aux constructions du paragraphe \ref{kato}. Comme $H^2(L,V)$ est le dual de $H^0(L,V^*(1))$, la suite duale de la suite (eq2) s'ecrit :
\[ 0 \to H_f^1(L,V^*(1))^* \to \dcris^L(V^*(1))^* \oplus t^*_{V^*(1)}(L)
\to \dcris^L(V^*(1))^* \to H^2(L,V) \to 0, \]
et en composant cette suite avec la suite (eq1), on obtient 
une suite exacte de $\Qp[G]$-modules :
\begin{multline*}
0 \to H^0(L,V) \to \dcris^L(V) \to \dcris^L(V) \oplus t_V(L) \to H^1(L,V) \\ 
\to \dcris^L(V^*(1))^* \oplus t^*_{V^*(1)}(L) \to \dcris^L(V^*(1))^* \to H^2(L,V) \to 0. 
\end{multline*}
En utilisant la suite exacte $0 \to t^*_{V^*(1)}(L) \to \ddR^L(V) \to t_{V}(L) \to 0$, on en d\'eduit  des isomorphismes canoniques : 
\begin{align*}
& \delta'_{V,L/K} : \det_{\Qp[G]} \ddR^L(V) \otimes \det_{\Qp[G]} \RR\Gamma(L,V) \iso \Qp[G], \tag{eq3} \\
& \Delta_{EP}(L/K,V) \iso \det^{-1}_{\Qp[G]}(\ddR^L(V)) \otimes \det_{\Qp[G]} (\Ind_{L/\Qp}V).
\end{align*}
En composant le dernier isomorphisme avec l'application $\beta_{V,L/K}$, on obtient une trivialisation canonique de la droite d'Euler-Poincar\'e :
\[ \delta_{V,L/K} : \Delta_{EP}(L/K,V) \iso \Qp[G]_{V,L/K}. \]

Nous pouvons maintenant enfin \'enoncer les conjectures $C_{\mathrm{EP}}(L/K,V)$ et $C_{\mathrm{EP}}(K,V)$ (voir \cite{FP,PR2,K2}).

\begin{conj}[$C_{\mathrm{EP}}(L/K,V)$]\label{new244}
Si $V = \Qp \otimes_{\Zp} T$ est une repr\'esentation potentiellement semi-stable de $G_K$ et si $L/K$ est une extension ab\'elienne finie, alors l'application
$\delta_{V,L/K}$ envoie $\Delta_{EP}(L/K,T)$ sur $\Zp[G]_{V,L/K}$.
\end{conj}

Si $L=K$, alors on peut reformuler cette conjecture en termes des nombres de
Tamagawa locaux (voir \cite{FP} et la d\'efinition \ref{deftamag} ci-dessus). Soit $\omega \in \det_{\Qp} \ddR^K(V)$ une base v\'erifiant $\omega \simeq \omega_2^{-1} \otimes \omega_1$ avec
$\omega_1 \in \det_{\Qp}t_{V}(K)$ et  $\omega_2 \in \det_{\Qp}t_{V^*(1)}(K)$. Soit $\omega_T$ une base de $\Ind_{K/\Qp}(T)$ et soit
$\alpha_{V,K}(\omega,T)=\alpha_{V,K}(\omega^{-1}\otimes\omega_T)$.

\begin{conj}[$C_{\mathrm{EP}}(K,V)$]\label{new245}
Si $V = \Qp \otimes_{\Zp} T$ est une repr\'esentation potentiellement semi-stable de $G_K$, alors : 
\[ \frac{\Tam_{K,\omega_1}^0(T)}{\Tam_{K,\omega_2}^0(T^*(1))}=
| d_K|_p^{\dim V /2}  \left| \Gamma^*(V) \frac{\alpha_{V,K}(\omega, T)}{\eps(K,V)} \right| _p. \]
\end{conj}

La proposition suivante rassemble quelques propri\'et\'es fonctorielles de la conjecture $C_{\mathrm{EP}}(L/K,V)$. 

\begin{prop}\label{old242} 
\begin{enumerate}
\item Les conjectures $C_{\mathrm{EP}}(L/K,V)$ et $C_{\mathrm{EP}}(L/K, V^*(1))$ sont \'equivalentes.
\item Si $0 \to V' \to V \to V'' \to 0$ est une suite exacte de repr\'esentations
potentiellement semi-stables et si la conjecture $C_{\mathrm{EP}}$ est vraie pour deux des repr\'esentations $V'$, $V$ et $V''$, alors elle est vraie pour la troisi\`eme.
\item Si $M/K$ est une extension de $K$ contenue dans $L$ et si $C_{\mathrm{EP}}(L/K,V)$ est vraie, alors les conjectures $C_{\mathrm{EP}}(L/M,V)$ et $C_{\mathrm{EP}}(M/K,V)$ le sont aussi.
\item Si la conjecture $C_{\mathrm{EP}}(L/K,V)$ est vraie, alors pour tout caract\`ere $\eta \in X(G)$, la conjecture $C_{\mathrm{EP}}(K, V(\eta))$ est vraie.
\end{enumerate}
\end{prop}

\begin{proof}
La d\'emonstration se fait comme dans \cite[C.2.9]{PR2}, en utilisant en plus les remarques suivantes :
\begin{enumerate}
\item La dualit\'e locale donne un isomorphisme $\det_{\Zp[G]} \RR \mathrm{Hom}_{\Zp} (\RR\Gamma (L,T),\Zp) \simeq \det_{\Zp[G]} \RR\Gamma (L,T^*(1))$.
\item Pour le triangle exact $\RR\Gamma (L, T') \to \RR\Gamma (L, T) \to \RR\Gamma (L, T'') \to \RR\Gamma (L, T')[1]$, on a un isomorphisme fonctoriel $\det_{\Zp[G]} \RR\Gamma (L, T) \iso \det_{\Zp[G]}
\RR\Gamma (L, T') \otimes \det_{\Zp[G]} \RR\Gamma (L, T'')$ (voir \cite[prop 7]{KM}).
\item Si on pose $H=\mathrm{Gal}(L/M)$ et $D_{L/M} = \dpst(\Ind_{L/M} V)$ et $D_{L/K} = \dpst(\Ind_{L/K}V)$, alors pour tout $\eta \in X(H)$ on a :
\begin{align*}
\lambda (M/K)^{\dim V} \eps(D_{L/M,\eta}, \psi_{M,\eta}, \mu_M)
& = \eps(\Ind_{M/K} (D_{L/M,\eta}),\psi_{K,\eta},\mu_K) \\
& = \prod_{\substack{\hat{\eta} \in X(G) \\ \hat{\eta} \mapsto \eta}}
 \eps(D_{L/K,\hat{\eta}}, \psi_{\hat{\eta}},\mu_K).
\end{align*}
On en d\'eduit que la restriction transforme $\beta_{V,L/K}$ en $\beta_{V,L/M}$ et le fait que  $C_{\mathrm{EP}}(L/K,V)$ implique $C_{\mathrm{EP}}(L/M,V)$ r\'esulte maintenant de la proposition \ref{detuniv}. La deuxi\`eme implication est analogue : on voit facilement que la projection de $\Qp[G]$ sur $\Qp[G/H]$ transforme $\beta_{V,L/K}$ en $\beta_{V,M/K}$.
\item Soit $E$ un corps contenant toutes les valeurs des caract\`eres  $\eta \in X(G)$ et soit $V(\eta) = E(\eta) \otimes_{\Qp} V$. Si on note $\mathcal{A}(G)$ l'ordre maximal de $E[G]$, alors on a des isomorphismes canoniques :
\[ \mathcal{A}(G) \otimes_{\Zp[G]}^{\mathbf{L}}  \RR\Gamma (L,T) \iso
\RR\Gamma (K, \mathcal{A}(G) \otimes_{\Zp} T) \iso \oplus_{\eta\in X(G)} \RR\Gamma (K,T(\eta)). \]
Si la conjecture $C_{\mathrm{EP}}(L/K,V)$ est vraie, alors l'application $\delta_{V,L/K}$  envoie $\Delta_{EP} (K,\mathcal{A}(G) \otimes_{\Zp} T)$ sur $\mathcal{A}(G)_{V,L/K} =\mathcal{A}(G) \otimes_{\Zp[G]} \Zp[G]_{V,L/K}$. En d\'ecomposant cet isomorphisme caract\`ere par caract\`ere, on en d\'eduit les conjectures $C_{\mathrm{EP}}(K, V(\eta))$ pour tous les caract\`eres $\eta \in X(G)$.
\end{enumerate}
\end{proof}

\section{L'exponentielle de Perrin-Riou}\label{expbpr}

Dans tout ce chapitre, on suppose que $K$ est une extension non-ramifi\'ee de $\Qp$. On commence par des rappels et des compl\'ements sur l'exponentielle de Perrin-Riou, ce qui nous permet d'\'enoncer la conjecture $C_{\mathrm{Iw}}(K_\infty/K,V)$. Dans le chapitre suivant, 
on montre que $C_{\mathrm{Iw}}(K_\infty/K,V)$ est \'equivalente \`a $C_{\mathrm{EP}}(K_n/K,V)$ pour tout $n \geq 1$ et finalement, on d\'emontre la conjecture $C_{\mathrm{Iw}}(K_\infty/K,V)$. 

\Subsection{Rappels et compl\'ements}\label{benois}

L'objet de ce paragraphe est de rappeler la construction et certaines propri\'et\'es de l'exponentielle de Perrin-Riou, tout d'abord telle qu'elle a \'et\'e d\'efinie par Perrin-Riou elle-m\^eme dans \cite{PR1}, puis ensuite (dans le paragraphe suivant) telle qu'elle a \'et\'e fait par l'un d'entre nous dans \cite{Ben}. 

Rappelons que $K_n=K(\zn)$, que $K_\infty = \cup_{n \geq 1} K_n$, que $\Gamma = \mathrm{Gal}(K_{\infty}/K)$ et que $G_n=\mathrm{Gal}(K_n/K) \simeq \Gamma / \Gamma_n$. On fixe un g\'en\'erateur topologique $\gamma_1$ de $\Gamma_1$ et on pose $\gamma_n=\gamma_1^{p^{n-1}}$ ce qui fait de $\gamma_n$ un g\'en\'erateur topologique de $\Gamma_n$. Si on note $\Delta_K$ le sous-groupe de torsion de $\Gamma$, alors on a $\Gamma \simeq \Delta_K \times \Gamma_1$ et $\Lambda = \Zp[\Delta_K] \otimes_{\Zp} \Zp[[\Gamma_1]]$. On d\'efinit une action de  $\Gamma$ sur $K[[X]]$ par la formule :
\[ g(X)=(1+X)^{\chi (g)}-1, \]
o\`u $\chi : \Gamma \to \Zp^\times$ est le caract\`ere cyclotomique. On munit par ailleurs $K[[X]]$ d'un Frobenius $\varphi$ et d'un op\'erateur diff\'erentiel $\partial$ en posant :
\begin{align*}
\varphi (\sum_{i=0}^{+\infty} a_i X^i) & = \sum_{i=0}^{+\infty} a_i^{\sigma} \varphi (X)^i, \quad\text{o\`u $\varphi (X)=(1+X)^p-1$,} \\
\partial & = (1+X) \frac{d}{dX}. 
\end{align*}
On v\'erifie facilement que  $\partial \circ \varphi = p \varphi \circ \partial$. Soit $\psi : K[[X]] \to K[[X]]$ l'op\'erateur d\'efini par la formule :
\[ \psi (f(X)) = \frac{1}{p} \varphi^{-1} \left( 
\sum_{\zeta^p=1}  f(\zeta(1+X)-1) \right), \]
qui est compatible avec la d\'efinition du paragraphe \ref{wach}. Il est classique que $\OO_K[[X]]^{\psi=0} = \{f \in \OO_K[[X]] \mid \psi (f) =0 \}$ est un 
$\OO_K[[\Gamma]]$-module libre engendr\'e par $1+X$.

On note $\calH$ l'ensemble des s\'eries formelles $f(X) \in \Qp[[X]]$ qui convergent sur le disque unit\'e ouvert, c'est-\`a-dire $\{x \in \Cp, 
|x|_p<1\}$, et l'on pose $\calH(\Gamma_1) = \{ f(\gamma_1 -1),  f \in \calH \}$ et $\calH(\Gamma) = \Qp[\Delta_K] \otimes_{\Qp} \calH(\Gamma_1)$. Pour tout $\Lambda$-module $N$, l'homomorphisme naturel $N \to N_{\Gamma_n}$ se prolonge en une application $\calH(\Gamma) \otimes_{\Lambda} N \to \Qp \otimes_{\Zp} N_{\Gamma_n}$.

Si $V$ est une repr\'esentation cristalline de $G_K$, alors on pose $\calD(V) = \OO_K[[X]]^{\psi=0} \otimes_{\OO_K} \dcris(V)$. Pour tout $k \in \ZZ$, on d\'efinit une application $\Delta_k : \calD(V) \to \dcris(V)/(1-p^k\varphi)\dcris(V)$ par la formule $\Delta_k(f) = (\partial^k f)(0) \mod{(1-p^k\varphi)\dcris(V)}$. Si on \'ecrit $\Delta = \oplus_{k \in \ZZ} \Delta_k$, alors pour tout $f \in \calD(V)^{\Delta=0}$, l'\'equation $(1- \varphi)F(X)=f(X)$ a une solution dans $\calH(V) = \calH \otimes_{\Qp} \dcris (V)$ et on en d\'eduit une application :
\begin{align*}
\Xi^{\eps}_{V,n} : \calD(V)^{\Delta=0} & \to \ddR^{K_n}(V) / \dcris(V)^{\varphi=1}, \\
f & \mapsto p^{-n}(\sigma \otimes \varphi)^{-n}(F)(\zn -1).
\end{align*}
Dans \cite{PR1}, Perrin-Riou a d\'emontr\'e le r\'esultat suivant.

\begin{theo}\label{old313}
Si $h$ est un entier tel que $\Fil^{-h} \ddR^K(V)= \ddR^K(V)$, alors pour tout $i \in \ZZ$ v\'erifiant $i+h \geq 1$, il existe un $\Lambda$-homomorphisme (appel\'e exponentielle \'elargie, ou exponentielle de Perrin-Riou) :
\[ \Exp_{V(i),h+i}^{\eps} : \calD(V(i))^{\Delta=0} \to \calH (\Gamma) \otimes_{\Lambda} (H^1_{\mathrm{Iw}} (K,T(i)) /  T(i)^{H_K} ), \]
v\'erifiant les propri\'et\'es suivantes:
\begin{enumerate}
  \item Le diagramme ci-dessous est commutatif :
\[  \begin{CD}
\calD(V(i))^{\Delta=0} @> \Exp^{\eps}_{V(i),h}>>  \calH(\Gamma) \otimes_{\Lambda} (H^1_{\mathrm{Iw}}(K,T(i)) /T(i)^{H_K}) \\
@V \Xi^{\eps}_{V(i),n}VV         @VV\mathrm{pr}_{T(i),n}V \\
\ddR^{K_n}(V(i)) @>(h+i-1)! \exp_{V(i),K_n}>> H^1(K_n,V(i)) /H^1(\Gamma_n,V(i)^{H_K}).
\end{CD} \]
  \item Soit $e_1 = \eps^{-1} \otimes t$ le g\'en\'erateur de $\dcris (\Qp(-1))$ associ\'e au choix de $\eps$ et soit :
\[ {\mathrm{Tw}}^{\eps}_{V(i),k} :H^1_{\mathrm{Iw}} (K,V(i)) \to H^1_{\mathrm{Iw}}(K,V(i+k)) \]
l'application d\'efinie par
$\mathrm{Tw}^\eps_{V(i),k}(x) = x \otimes {\eps}^{\otimes k}$.
On a alors :
\[ \Exp^{\eps}_{V(i+1),h+1} = - \mathrm{Tw}^{\eps}_{V(i),1} \circ \Exp^{\eps}_{V(i),h} \circ (\partial \otimes e_1). \]
  \item Si : 
\[ \ell_m = m-\frac{\log (\gamma_1)}{\log \chi (\gamma_1)}, \] 
alors $\Exp^{\eps}_{V(i),h+1} = \ell_h \Exp^{\eps}_{V(i),h}$.
\end{enumerate}
\end{theo}

On d\'eduit de ce th\'eor\`eme plusieurs formules qui nous sont utiles. Tout d'abord, en it\'erant (2), on obtient :
\[ \Exp_{V(i),h+i}^\eps = (-1)^i {\mathrm{Tw}}^\eps_{V,k} \circ \Exp_{V,h}^\eps \circ (\partial^i \otimes e_i). \]

Soit $\calK(\Gamma)$ l'anneau total des fractions de $\calH(\Gamma)$ (il suffit en fait d'inverser les $\ell_i$). Le (3) permet de d\'efinir pour tout $h \in \ZZ$ :
\[ \Exp_{V,h}^{\eps} : \calD(V)^{\Delta=0} \to \calK (\Gamma) \otimes_{\Lambda}
\left (H^1_{\mathrm{Iw}}(K,T)/T^{H_K}\right). \]
En particulier, si $\Fil^0 \dcris (V) = \dcris (V)$, alors on dispose de l'application : 
\[ \Exp_{V,0}^{\eps} : \calD(V)^{\Delta=0} \to \calH (\Gamma) \otimes_{\Lambda}
\left(H^1_{\mathrm{Iw}}(K,T)/T^{H_K}\right), \]
qui est telle que pour tout $i \geq 1$, 
si l'on pose $\Xi^{(k),\eps}_{V,n} = \Xi^{\eps}_{V(i),n} \circ (\partial^{-k} \otimes e_{-k})$, alors le diagramme ci-dessous commute :
\[ \begin{CD}
\calD(V)^{\Delta=0}(i) @>{\mathrm{Tw}}^\eps_{V,i} \circ \Exp^{\eps}_{V,0}>> 
\calH (\Gamma) \otimes_{\Lambda} \left(H^1_{\mathrm{Iw}}(K,T(i))/T(i)^{H_K} \right) \\
@V\Xi^{(k),\eps}_{V,n}VV         @VV\mathrm{pr}_{T(i),n}V \\
\ddR^{K_n}(V(i)) @>(-1)^i(i-1)!\exp_{V(i),K_n}>> 
H^1(K_n, V(i)) / H^1(\Gamma_n,V(i)^{H_K}).
\end{CD} \]

Rappelons \`a pr\'esent quelques r\'esultats techniques concernant l'application $\Xi^{\eps}_{V,n}$ et qui sont d\'emontr\'es dans \cite[\S 3.4]{PR1}. L'homomorphisme $\Delta$ donne lieu \`a une suite exacte courte :
\[ 0 \to \calD(V)^{\Delta =0} \to \calD(V) \overset{\Delta}{\to}
\bigoplus_{j \in \ZZ} \left( \frac{\dcris (V)}{(1-p^j\varphi)\dcris (V)} \right) (j) \to 0, \]
qui induit une suite exacte :
\[ 0 \to \frac{\dcris (V)}{(1-\varphi)\dcris (V)} \to(\calD(V)^{\Delta=0})_{\Gamma_n}
\to\calD(V)_{\Gamma_n}\to\frac{\dcris (V)}{(1-\varphi)\dcris (V)}\to 0. \tag{eq4} \]
La deuxi\`eme fl\`eche de cette suite est donn\'ee par la formule
$\tilde{d} \mapsto d \otimes (\gamma_n-1)(1+X)$ si $\tilde{d} = d \mod{(1-\varphi)\dcris (V)}$.

L'application $\Xi^{\eps}_{V,n}$ se factorise par $(\gamma_n-1) \calD(V)^{\Delta=0}$ et on note :
\[ \tilde{\Xi}^{\eps}_{V,n} : (\calD(V)^{\Delta=0})_{\Gamma_n} \to
\ddR^{K_n} (V)/\dcris(V)^{\varphi=1} \]
la fl\`eche qui s'en d\'eduit. Soit :
\[ \Exp^{\eps}_{V,h,n}:(\calD(V)^{\Delta=0})_{\Gamma_n} \to
(\Qp\otimes_{\Zp} H^1_{\mathrm{Iw}}(K,T)/T^{H_K})_{\Gamma_n} \]
l'application d\'eduite de $\Exp^{\eps}_{V,h}$.

\begin{prop}\label{old316} 
\begin{enumerate}
\item La suite :
\[ 0 \to \dcris(V)/(1-\varphi)\dcris (V) \to \ker \tilde{\Xi}^{\eps}_{V,n} \overset{f}{\to}
\dcris (V)^{\varphi=p^{-1}} \to 0, \]
o\`u $f(\alpha(X))=\alpha (0)$, est exacte.
\item L'application $\alpha \mapsto (1-\varphi)\mathrm{Tr}_{K_n/K} (\alpha)$
induit un isomorphisme :
\[ \mathrm{coker}(\tilde{\Xi}^{\eps}_{V,n}) \iso \dcris (V) / (1-p^{-1}\varphi^{-1}) \dcris (V). \]
\item On a une suite exacte :
\begin{multline*}
0 \to \ker (\tilde{\Xi}^{\eps}_{V,n})\to \ker (\Exp^{\eps}_{V,h,n})\to\Fil^0 \ddR^{K_n} (V)/V^{G_K} \\
\to \dcris (V)/(1-p^{-1}\varphi^{-1})\dcris (V)
\to (V^*(1)^{G_K})^* \to 0. \tag{eq5}
\end{multline*}
\end{enumerate}
\end{prop}

\begin{proof}
Voir \cite[3.4.4-3.4.5]{PR1}.
\end{proof}

\Subsection{L'application exponentielle et les $(\varphi,\Gamma)$-modules}\label{altexpo}

Nous rappelons maintenant la construction de l'exponentielle en termes de $(\varphi,\Gamma)$-modules qu'a donn\'ee l'un d'entre nous (dans \cite{Ben}). On suppose d\'esormais que $V$ est une repr\'esentation cristalline de $G_K$ qui est positive, c'est-\`a-dire que les oppos\'es des poids de Hodge-Tate de $V$ sont $0 = r_1 \leq  r_2 \leq \cdots \leq r_d = h$. On fixe un $\Zp$-r\'eseau $T$ de $V$ stable par $G_K$, et on d\'efinit un $\OO_K$-r\'eseau $M$ de $\dcris (V)$ par :
\[  M = \{ f(X) \in (\bhol \otimes_{\aplus_K} \nwach(T))^\Gamma \mid f(0) \in \nwach(T)/ X \nwach(T) \}. \]
La proposition V.1 de \cite{Ber2} nous dit que le d\'eterminant de l'isomorphisme de comparaison $\bdR  \otimes_{\Qp} \Ind_{K/\Qp} V \simeq \bdR  \otimes_{\Qp} \dcris (V)$,  
calcul\'e dans des bases de $T$ et de $M$, appartient \`a $\hat{\OO}_{K^{\mathrm{nr}}}^\times t^{r_1+\hdots +r_d}$, c'est-\`a-dire que dans les notations de la section \ref{fontaine}, on a $\alpha_{K,V}(M,T) \in \hat{\OO}_{K^{\mathrm{nr}}}^\times$.

L'anneau $\OO_K[[X]]$ est muni comme ci-dessus des op\'erateurs $\psi$ et $\partial = (1+X) d/dX$, et on pose $\calD(T) = \OO_K[[X]]^{\psi=0} \otimes_{\OO_K} M$, o\`u $M$ est le r\'eseau de $\dcris(V)$ que l'on vient de d\'efinir. 

Pour des raisons techniques, on remplace le complexe
$C_{\varphi,\gamma_n}(K_n,T)$ par le complexe
$\varphi^{-n}(C_{\varphi,\gamma_n}(K_n,T))$ de $(\varphi,\Gamma)$-modules
sur $\afont_{K_n}=\varphi^{-n}(\afont_{K})$, complexe qui est isomorphe \`a $C_{\varphi,\gamma_n}(K_n,T)$. On pose $X_n = [\eps^{1/p^n}]-1$.

Supposons d'abord que $V^{H_K}=0$. 
Soit $n_1,\hdots ,n_d$ une base de $\nwach(T)$ et soit
$m=\sum_{i=1}^d a_i(X) \otimes n_i$, avec $a_i(X) \in \bhol$, un 
\'el\'ement de $M$. Si $\gamma \in \Gamma$, alors un petit calcul qui utilise les congruences $\gamma(n_i) \equiv n_i \mod{X \nwach(T)}$ montre que si l'on pose $c_k = \prod_{j=1}^k (\chi^j (\gamma)-1)$ pour un g\'en\'erateur topologique $\gamma$ de $\Gamma$ et pour tout $k \geq 1$, alors $a_i(X)$ appartiennent \`a l'anneau $\afont'_K = \aplus_K [[ X^k / c_k, k\geq 0 ]]$. Si $\alpha = f(x) \otimes m \in \calD(T)$, alors on pose $E_{k,n}(f) \otimes m = \sum_{i=1}^d (a_i(X) E_{k,n}(f)) \otimes n_i$, o\`u :
\[ E_{k,n}(f) = \sum_{j=1}^\infty \frac{(1-k)(2-k)\cdots
 \cdot (j-k-1)}{t^j}
p^{n(j-1)} \partial^{-j}(f(X_n)). \]
On tronque les s\'eries $a_i(X)/t^j$ modulo $X$ et on note $\overline{a_iE_{k,n}(f)}$ les s\'eries que l'on obtient ainsi. Soit :
\[ \calE_{T,k,n} (\alpha) = \sum_{i=1}^d
n_i\otimes \overline{a_iE_{k,n}(f)}\otimes \eps^{\otimes k}. \]
On v\'erifie que $\calE_{T,k,n} (\alpha) \in \varphi^{-n}(\dfont(T(k)))$ et on d\'efinit $\calF_{T,k,n} (\alpha) \in \varphi^{-n}(\dfont(T(k)))$ par :
\[ (1-\varphi)\calF_{T,k,n}(\alpha) = (1-\gamma_n) \calE_{T,k,n}(\alpha), \]
ce qui fait que $(\calE_{T,k,n}(\alpha),\calF_{T,k,n}(\alpha))$ d\'efinit une classe de cohomologie dans l'espace $H^1(\varphi^{-n}(C_{\varphi,\gamma_n}(K_n,T(k))))$. En composant avec l'isomorphisme :
\[ H^1(\varphi^{-n}(C_{\varphi,\gamma_n}(K_n,T(k))))\iso H^1(K_n,T(k)), \]
on obtient un homomorphisme :
\[ \Omega_{T,k,n}^\eps : \calD(T) \to H^1(K_n,T(k)). \]

Revenons maintenant au cas g\'en\'eral (on ne suppose plus que $V^{H_K}=0$) et posons : 
\[ \calH(T) = \{ \alpha \in \OO_K[[X]] \otimes_{\OO_K} M \mid \psi(\alpha)=\alpha \}. \]
Un petit calcul montre qu'on a une suite exacte : 
\[ 0 \to M^{\varphi=1} \to \calH(T) \overset{1-\varphi}{\longrightarrow} \calD(T)^{\Delta_0=0} \to 0, \] 
et la m\^eme construction qu'avant fournit un homomorphisme :
\[ \Sigma_{T,k,n}^{\eps} : \calH(T) \to H^1(K_n,T(k)) \]
qui s'inscrit dans un diagramme :
\[ \begin{CD}
\calH(T) @>{\Sigma_{T,k,n}^{\eps}}>> H^1(K_n,T(k)) \\
@V{1-\varphi}VV @VVV \\
\calD(T)^{\Delta_0=0} @>>> H^1(K_n,T(k)) / H^1(\Gamma_n,T(k)^{H_K}).
\end{CD} \]
En particulier, si $V^{H_K}=0$, alors on a :
\[ \Sigma_{T,k,n}^{\eps}(\alpha) = \Omega_{T,k,n}^{\eps}((1-\varphi)\alpha). \]

Les r\'esultats suivants sont d\'emontr\'es dans \cite[th\'eor\`emes 4.3 et 5.1.2]{Ben}.

\begin{prop}\label{old423}
Si $V$ est une repr\'esentation positive v\'erifiant $V^{H_K}=0$ et si 
$\alpha \in \calD(T)$, alors :
\begin{enumerate}
  \item Pour tous $k \geq 1$ et $n\geq 1$, on a :
\[ \Omega_{T,k,n}^\eps(\alpha) = (-1)^k (k-1)! \exp_{V(k),K_n} (F_k(\zn-1)), \] 
o\`u $F_k(X)$ est une solution de l'\'equation
$(1-\varphi)F_k =(\partial^{-k}\otimes e_{-k})(\alpha)$ et $e_{-k}$ est le g\'en\'erateur de $\dcris(\Qp(k))$ associ\'e \`a $\eps$.
  \item Plus g\'en\'eralement, pour tous $k \in \ZZ$ et $n \geq 1$, on a : 
\[ \Omega_{T,k,n}^\eps((\varphi \otimes \sigma)^{-n}(\alpha)) =
\mathrm{pr}_{V(k),n} \circ \mathrm{Tw}_{V,k}^{\eps} \circ \Exp_{V,0}^{\eps} (\alpha). \]
  \item Soit $(\cdot,\cdot)_{T(k),n} : H^1(K_n,T(k)) \times H^1(K_n,T^*(1-k)) \overset{\cup}{\to} \Qp$ l'accouplement fourni par la dualit\'e locale. On a alors 
\begin{multline*}
(\Omega_{T,k,n}^\eps (\alpha), \Omega_{T^*(-h),h-k+1,n}^\eps (\beta))_{T(k),n}=
\\ (-1)^k p^{nh}  \prod_{m=1}^h (k-m) 
\mathrm{Tr}_{K/\Qp} \mathrm{res}
\left(\frac{1}{X} \left[\partial^{-k}\alpha(X_n),(\partial^{k-h-1}\otimes e_{-h})
\beta (X_n)\right] \frac{dX_n}{1+X_n} \right). 
\end{multline*}
\end{enumerate}
\end{prop}

\begin{rema}\label{old424} 
Cette proposition entra\^{\i}ne la loi de r\'eciprocit\'e de Perrin-Riou (le th\'eor\`eme \ref{old322} ci-dessous).
\end{rema}

Rappelons qu'une repr\'esentation cristalline $W$ telle que $W=W^{H_K}$ est n\'ecessairement de la forme $W= \oplus_{i \in \ZZ} \Qp(i)^{d_i}$.

\begin{prop}\label{old425} 
Si $V$ n'a pas de sous-quotient isomorphe \`a $\Qp(m)$, avec $m \in \ZZ$, alors
pour tout $k\notin [1,h]$, l'application $\Omega_{T,k,1}^\eps$ induit un isomorphisme
\[ \calD(T)(k)_{\Gamma_1} \iso (\varphi^* \nwach(T)(k))^{\psi=1}_
{\Gamma_1} \hookrightarrow H_{\mathrm{Iw}}^1(K,T(k))_{\Gamma_1}. \]
\end{prop}

La d\'emonstration de cette proposition fait l'objet du reste de ce paragraphe. Pour simplifier la notation on pose $\Omega_{T,k}^\eps = \Omega_{T,k,1}^\eps$.

\begin{lemm}\label{old426} 
Le $\Lambda$-module $(\varphi^* \nwach(T))^{\psi=1}$ est libre de rang $d=\dim(V)$ et 
pour tout $k \notin [1,h]$, l'application $\Omega_{T,k}^\eps$ induit une injection :
\[ \calD(T)(k)_{\Gamma_1} \overset{\Omega_{T,k}^\eps}{\hookrightarrow} 
(\varphi^* \nwach(T)(k))^{\psi=1}_
{\Gamma_1} \hookrightarrow H_{\mathrm{Iw}}^1 (K,T(k))_{\Gamma_1}. \]
\end{lemm}

\begin{proof} 
Rappelons que $\afont \subset W(\et)$ et soit $\afont^{>0}$ l'ensemble des $x = \sum_{k=0}^\infty p^k [x_k] \in \afont$ tels que $x_k \in \mathfrak{m}_{\et}$ pour tout $k \geq 0$. On a une suite exacte scind\'ee $0 \to \afont^{>0} \to \afont^+ 
\leftrightarrows W(\overline{k}_K) \to 0$. Si $\dfont^{>0}(T) = (\afont^{>0} \otimes_{\Zp} T)^{H_K}$, alors on a une suite exacte :
\[ 0 \to \dfont^{>0}(T)  \to \dfont^+(T) \leftrightarrows (W(\overline{k}_K) \otimes_{\Zp} T)^{H_K} \to 0. \]

Comme la restriction de $\psi$ \`a $W(\overline{k}_K)$ co\"{\i}ncide avec $\varphi^{-1}$, on a $((W(\overline{k}_K) \otimes_{\Zp} T)^{\psi=1})^{H_K} = T^{H_K}=0$ par hypoth\`ese , d'o\`u $\dfont^+(T)^{\psi=1}=\dfont^{>0}(T)^{\psi=1}$.
Comme $T$ est positive, on a $\varphi^*(\nwach(T)) \subset \nwach(T)\subset \dfont^+(T)$,
d'o\`u $(\varphi^*\nwach(T))^{\psi=1}\subset \dfont^{>0}(T)^{\psi=1}$.
L'op\'erateur $1-\varphi$ est inversible sur $\dfont^{>0}(T)$, d'inverse $\sum_{j \geq 0} \varphi^j$,  et un petit calcul montre qu'il donne lieu \`a un isomorphisme :
\[ (\varphi^*\nwach(T))^{\psi=1} \iso (\varphi^*\nwach(T) )^{\psi=0}. \]

Pour montrer que le $\Lambda$-module $(\varphi^* \nwach(T))^{\psi=1}$ est libre de rang $d=\dim V$, il reste enfin \`a remarquer que $(\varphi^*\nwach(T))^{\psi=0}$ est un
$\Lambda$-module libre engendr\'e par les \'el\'ements $\varphi (n_i) \otimes (1+X)$, o\`u $n_1,\hdots ,n_d$ est une base de $\nwach(T)$.

Posons maintenant $\alpha = f \otimes m \in \calD(T)$, o\`u $m=\sum_{i=1}^d 
a_i(X) \otimes n_i \in (\bhol \otimes_{\aplus_K} \nwach(T))^{\Gamma}$ et $f \in \OO_K[[X]]^{\psi=0}$. La proposition 3.1.3 de \cite{Ben} montre que l'on a alors :
\begin{align*}
(1-\gamma_1)\calE_{T,k}(\alpha) & \equiv (1-\gamma_1)(E_{k,1}(f)\otimes m \otimes 
\eps^{\otimes k})\\
& \equiv \frac{1-\chi (\gamma_1)^k}{p^k}f(X_1) \otimes m \otimes \eps^{\otimes k}
\mod {X_1 \Qp [[X_1]] \otimes_{\aplus_K} \nwach(T)(k)}\\
& \equiv  \sum_{i=1}^d  a_i(0)
\frac{1-\chi (\gamma_1)^k}{p^k}f(X_1) \otimes n_i \otimes \eps^{\otimes k} \mod {X_1 \nwach(T)(k)}
\end{align*}
ce qui fait que $(1-\gamma_1) \calE_{T,k}(\alpha) \in \afont_{K_1}\otimes_{\afont_{K_1}^+}  \nwach(T)(k)$. D'autre part, soit $\psi_1$ l'op\'erateur $\psi$ agissant sur $\afont_{K_1}$. Comme $\psi_1(\partial^{-j} f(X_1))=0$
et $\psi_1(X^m x)=X_1^m \psi_1(x)$, on a $\psi_1(\calE_{T,k}(\alpha))=0$ et donc $(1-\gamma_1)\calE_{T,k}(\alpha)\in(\afont_{K_1} \otimes_{\afont_{K_1}^+}\nwach(T)(k))^{\psi_1=0}$, d'o\`u :
\[ \calF_{T,k}(\alpha) = (1-\varphi)^{-1}(1-\gamma_1)\calE_{T,k}(\alpha)\in
(\afont_{K_1}^+ \otimes_{\afont_{K}^+} \nwach(T)(k))^{\psi_1=1}, \]
et la formule $\alpha \mapsto \varphi (\calF_{T,k}(\alpha))$ d\'efinit 
un homomorphisme $\calD(T)(k)\to (\varphi^*\nwach(T)(k))^{\psi=1}$ qui induit un diagramme commutatif :
\[ \xymatrix{
\calD(T)(k)_{\Gamma_1} \ar[r] \ar[ddr]_{\Omega_{T,k}^{\eps}} & (\varphi^*\nwach(T)(k)
)^{\psi=1}_{\Gamma_1} \ar[d]\\
&H_{\mathrm{Iw}}^1(K,T(k))_{\Gamma_1} \ar@{_{(}->}[d]\\
& H^1(K_1,T(k)). }  \]

Le (3) de la proposition \ref{old424}  entra\^{\i}ne l'injectivit\'e de $\Omega_{T,k}^\eps$
pour $k\notin [1,h]$ et l'application $\calD(T)(k)_{\Gamma_1}\to (\varphi^*\nwach(T)(k)
 )^{\psi=1}_{\Gamma_1}$ est donc injective. D'autre part, $\calD(T)(k)$
et  $(\varphi^*\nwach(T)(k))^{\psi=1}$ sont des $\Lambda$-modules libres de 
m\^eme rang, donc $\calD(T)(k)_{\Gamma_1}$ et $(\varphi^*\nwach(T)(k))^{\psi=1}_{\Gamma_1}$ sont des $\Zp$-modules libres 
de m\^eme rang et $\calD(T)(k)_{\Gamma_1}$ est un r\'eseau de $(\varphi^*\nwach(T)(k))^{\psi=1}_{\Gamma_1}$.  On en d\'eduit que l'application 
$(\varphi^*\nwach(T)(k))^{\psi=1}_{\Gamma_1} \to H_{\mathrm{Iw}}^1(K,T(k))_{\Gamma_1}$
est injective et le lemme est d\'emontr\'e.
\end{proof} 

\begin{rema}\label{old427} 
On donnera plus bas une autre preuve de l'injectivit\'e de l'application $(\varphi^*\nwach(T)(k))^{\psi=1}_{\Gamma_1} \to H_{\mathrm{Iw}}^1(K,T(k))_{\Gamma_1}$ (voir Proposition \ref{old432}).
\end{rema}

\begin{lemm}\label{old428} 
\begin{enumerate}
\item Si $f(X_1)$, $g(X_1) \in \afont_{K_1}^+$, alors pour tout $k \in \ZZ$ et $n \geq 0$, on a :
\[ \mathrm{res} \left(E_{k,1}(f) t^n g(X_1)\frac{dX_1}{1+X_1}\right)
\equiv 0 \mod{\left( p^n\frac{\Gamma^* (k)}{\Gamma^* (k-n)}\right)}. \]
\item Si $m \otimes f \in \calD(T)$ et $\beta \in (\varphi^*\nwach(T(-h)))^{\psi=1}(h-k-1)$, alors :
\[ \mathrm{res} \left(E_{k,1}(f)[m,\varphi^{-1}(\beta)]_{V(k)}\frac{dX_1}{1+X_1}\right)
\equiv 0 \mod{\left( p^h \frac{\Gamma^* (k)}{\Gamma^* (k-h)}\right)}, \]
o\`u $h$ est le dernier saut de la filtration de Hodge de $\ddR(V)$
et o\`u on a $\Gamma ^*(k) / \Gamma^*(k-n)=(k-1)\times \cdots \times (k-n)$, m\^eme si $k\leq 0$.
\end{enumerate}
\end{lemm}

\begin{proof} 
On commence par montrer le (1). Si $\partial_1 = (1+X_1) d/dX_1$, alors $\partial_1 = p \partial$ et on a :
\[ \partial_1 \left (\frac{h(X_1)}{t^n}\right ) =
\frac{\partial_1 h(X_1)}{t^n}-np\frac{h(X_1)}{t^{n+1}}. \]
On en d\'eduit que :
\[ \mathrm{res} \left (\frac{h(X_1)}{t^{n+1}}\frac{dX_1}{1+X_1}\right )
=
\frac{1}{np}\mathrm{res}\left (\frac{\partial_1 h(X_1)}{t^n}\frac{dX_1}{1+X_1}\right ), \]
et par r\'ecurrence on obtient la formule :
\[ \mathrm{res}\left (\frac{h(X_1)}{t^n}\frac{dX_1}{1+X_1}\right )
= \frac{1}{(n-1)!p^{n-1}}\mathrm{res}\left (\frac{\partial_1^{n-1}h(X_1)}{X}\frac{dX_1}{1+X_1}\right ). \]
On applique cette formule au calcul du r\'esidu :
\[ \mathrm{res} \left (E_{k,1}(f)[m,\varphi^{-1}(\beta)]_{V(k)}
\frac{dX_1}{1+X_1}\right ). \]
La s\'erie $E_{k,1}$ est d\'efinie par :
\[ E_{k,1}(f)=
\sum_{j=1}^\infty \frac{(1-k)(2-k)\cdots (j-k-1)}{t^j}
p^{j-1} \partial^{-j}(f(X_1)). \]
Si $n\leq j-1$, on a $\frac{(n+1-k)\times \cdots \times (j-1-k)}{(j-n-1)!}\in \ZZ$, d'o\`u :
\begin{multline*}
\mathrm{res}\left (\frac{(1-k)(2-k)\cdots (j-k-1)}{t^j}
p^{j-1} \partial^{-j}f(X_1) 
t^n g(X_1)\frac{dX_1}{1+X_1}\right )=\\
\frac{(1-k)(2-k)\cdots (j-k-1)p^{j-1}}{(j-n-1)!p^{j-n-1}}
\mathrm{res}\left (\frac{\partial^{j-n-1}(\partial^{-j}f(X_1)g(X_1))}{X}
\frac{dX_1}{1+X_1}\right),
\end{multline*}
et on en d\'eduit la congruence (1) du lemme.

Montrons maintenant le (2); rappelons qu'on a pos\'e $c_j=(\chi (\gamma)-1)\cdots (\chi (\gamma)^j-1)$, o\`u $\gamma$ est un g\'en\'erateur topologique de $\Gamma$. On a $X=\exp(t)-1 = t+t^2/2!+\cdots$, d'o\`u :
\[ X^{j+h} = t^{j+h} \sum_{s_1,\cdots ,s_{j+h}\geq 0} \frac{t^{s_1+\cdots +s_{j+h}}}{(s_1+1)! \cdots (s_{j+h}+1)!}. \tag{eq10} \]
Comme $m \in M$ et $\varphi^{-1}(\beta)\in \afont_{K_1}^+ \otimes_{\aplus_K}
\nwach(T^*(-h))$, on a :
\[ [m,\varphi^{-1}(\beta)]_V =  X^h \sum_{j=0}^\infty  \frac{b_j(X_1)X^j}{c_j}, \tag{eq11} \]
pour des \'el\'ements $b_j(X_1)\in \afont_{K_1}^+$. En utilisant la congruence (1)
du lemme et les congruences \'evidentes $p^s/(s+1)! \equiv 0 \mod{p}$ et $p^s/c_s  \equiv 0 \mod{p}$ pour tout $s \geq 1$, on montre que :
\[ \mathrm{res} \left (E_{k,1}(f)b_j(X_1)
\frac{t^{s_1+\ldots +s_{j+h}}}{(s_1+1)!\cdots (s_{j+h}+1)!c_j}
\frac{dX_1}{1+X_1}\right )
\equiv 0  \mod{\left( p^h \frac{\Gamma^* (k)}{\Gamma^* (k-h)}\right)}, \]
et la congruence (2) d\'ecoule maintenant des formules (eq10) et (eq11).
\end{proof}

\begin{prop}\label{old429}
Soient $(\cdot,\cdot)_{T(k)} : H^1(K_1,T(k)) \times H^1(K_1,T^*(1-k))
\to \Zp$ l'accouplement fourni par la dualit\'e locale, $\alpha \in \calD(T)(k)$,
$\beta \in (\varphi^*\nwach(T^*(-h)))^{\psi=1}(h-k+1)$ et $\mathrm{cl} (\beta)\in H^1(K_1,T^*(1-k))$ la classe de cohomologie associ\'ee \`a $\beta$ 
via l'injection $(\varphi^*\nwach(T^*(-h)))^{\psi=1}(h-k+1) \hookrightarrow H^1_{\mathrm{Iw}}(K,T^*(1-k))$ suivie de la projection. On a alors : 
\[  (\Omega_{T,k}^{\eps}(\alpha), \mathrm{cl} (\beta))_{T(k)}
\equiv 0 \mod{\left( p^h \frac{\Gamma^* (k)}{\Gamma^* (k-h)}\right)}. \]
\end{prop}

\begin{proof} 
Soit $A \in \varphi^{-1}(\dfont(T^*(1-k))^{\psi=0})$ une solution de l'\'equation $(\gamma_1-1)A = \varphi^{-1}(\varphi-1)\beta$. La formule du cup-produit en termes de $(\varphi,\Gamma)$-modules (voir \cite[proposition 4.4]{H2}) s'\'ecrit :
\[ (\Omega_{T,k}^{\eps}(\alpha),\mathrm{cl} (\beta) )_{T(k)} =
-\mathrm{cl} \left ( [\gamma_1 \calE_{T,k}(\alpha),\varphi^{-1}(\beta)]_{V(k)}-
[\varphi \calF_{T,k}(\alpha),A]_{V(k)}\right ). \]
Pour calculer cette classe, on reprend les arguments de la preuve du th\'eor\`eme 
5.1.2 de \cite{Ben}. En termes de $(\varphi,\Gamma)$-modules, l'isomorphisme canonique $H^2(K_n,\Zp(1)) \simeq \Zp$ est donn\'e par la formule (voir \cite[th\'eor\`eme 2.2.6]{Ben}) :
\begin{align*}
\mathrm{TR}_n : H^2(\varphi^{-n}(C_{\varphi,\gamma_n}(K_n,\Zp(1)))) & \iso \Zp \\
(\mathrm{cl} (h(X_n) \otimes \eps)) & \mapsto
-\frac{p^n}{\log \chi (\gamma_n)}\mathrm{Tr}_{K/\Qp}
\left (\mathrm{res} \frac{h(X_n)dX_n}{1+X_n}\right).
\end{align*}
Si $\alpha = f \otimes m \in  \OO_K[[X]]^{\psi=0} \otimes_{\OO_K} M$, alors il existe $y$ tel que $\calE_{T,k}(\alpha) =  E_{k,1}(f) \otimes m +(\varphi-1)y$, $\calF_{T,k}(\alpha)=F_{T,k}(\alpha)+(1-\gamma_1)y$ et $(\varphi-1)F_{T,k}(\alpha)=(\gamma_1-1)(E_{k,1}(f) \otimes m)$.

On en d\'eduit que :
\begin{multline*}
(\Omega_{T,k}^{\eps}(\alpha),\mathrm{cl} (\beta))_{T(k)}  = \\
\frac{p}{\log \chi (\gamma_1)}\mathrm{Tr}_{K/\Qp}
\left (\mathrm{res} 
 ( [\gamma_1  E_{T,k}(f)\otimes m,\varphi^{-1}(\beta)]_{V(k)}-
 [\varphi  F_{T,k}(\alpha),A]_{V(k)} )
\frac{dX_1}{1+X_1}\right).
\end{multline*}
Comme $\psi_1(A)=0$, on a $\psi_1 ( [\varphi F_{T,k}(\alpha),A]_{V(k)} )=0$, d'o\`u (cf. \cite[lemme 2.2.2.1]{Ben}) :
\[ \mathrm{res} 
\left ([\varphi F_{T,k}(\alpha),A]_{V(k)}\frac{dX_1}{1+X_1}\right )=0. \]
D'autre part, comme $(\gamma_1-1)E_{k,1}(f)
\in \Qp[[X_1]]$ le (2) du lemme \ref{old428} nous donne :
\begin{multline*}
\mathrm{res} \left (
[\gamma_1  E_{T,k}(f)\otimes m,\varphi^{-1}(\beta)]_{V(k)}
\frac{dX_1}{1+X_1}\right )= \\
\mathrm{res} \left (
[  E_{T,k}(f)\otimes m,\varphi^{-1}(\beta)]_{V(k)}
\frac{dX_1}{1+X_1}\right )\equiv 0 
\mod{\left( p^h \frac{\Gamma^* (k)}{\Gamma^* (k-h)}\right)}.
\end{multline*}
\end{proof}

\begin{proof}[D\'emonstration de la proposition \ref{old425}]
La (3) de la proposition \ref{old424} nous donne :
\begin{multline*} 
\det_{\Zp} \left (
\Omega_{T,k}^\eps (\calD(T)(k)_{\Gamma_1}),
\Omega_{T^*(-h),1+h-k}^\eps (\calD(T^*(-h))(1+h-k)_{\Gamma_1})\right )_{T(k)}
\\ = \left (
{p^h\frac{\Gamma^*(k)}{\Gamma^*(k-h)}} \right )^{[K_1:\Qp]d}
\Zp,
\end{multline*}
o\`u $d= \dim (V)$. Par ailleurs, la proposition \ref{old429} donne l'inclusion :
\begin{multline*} 
\det_{\Zp} \left ( (\varphi^*\nwach(T)(k))^{\psi=1}_{\Gamma_1},
\Omega_{T^*(-h),1+h-k}^\eps (\calD(T^*(-h))(1+h-k)_{\Gamma_1})\right )_{T(k)}
\\Ê\subset \left (
{p^h\frac{\Gamma^*(k)}{\Gamma^*(k-h)}} \right )^{[K_1:\Qp]d}\Zp. 
\end{multline*}
Comme l'inclusion  $\calD(T)(k)_{\Gamma_1} \hookrightarrow (\varphi^* \nwach(T)(k) )^{\psi=1}_
{\Gamma_1}$ est d\'ej\`a \'etablie, on en d\'eduit que $\calD(T)(k)_{\Gamma_1} \iso (\varphi^* \nwach(T)(k) )^{\psi=1}_{\Gamma_1}$ et la proposition \ref{old425} est d\'emontr\'ee.
\end{proof}

\begin{coro}\label{old4211} 
Si $k\notin [1,h]$, alors :
\begin{multline*}
\left [\dfont(T(k))^{\psi=1}_{\Gamma_1}:(\varphi^*\nwach(T)(k))^{\psi=1}_{\Gamma_1} \right ]
\left [\dfont(T^*(1-k))^{\psi=1}_{\Gamma_1}:(\varphi^*\nwach(T^*(-h))(h+1-k))^{\psi=1}_{\Gamma_1} \right] \\ = \left(
{p^h\frac{\Gamma^*(k)}{\Gamma^*(k-h)}} \right )^{[K_1:\Qp]\dim (V)}.
\end{multline*}
\end{coro}

\begin{proof} 
Soit $\tilde{H}^1(K_1,T(k))=H^1(K_1,T(k)) / H^1(K_1,T(k))_{\mathrm{tor}}$. Comme par hypoth\`ese $V(k)^{G_K}=0$, on a $H^1(K_1,T(k))_{\mathrm{tor}} \simeq H^0(K_1, V(k)/T(k))$, d'o\`u :
\begin{multline*}
[H_{\mathrm{Iw}}^1(K,T(k))_{\Gamma_1}:(\varphi^*\nwach(T)(k))^{\psi=1}_{\Gamma_1}]
= \\ \frac{\sharp H^0(K_1,V(k)/T(k))}{\sharp H^0(K_1,V^*(1-k)/T^*(1-k))} [\tilde{H}^1(K_1,T(k)):(\varphi^*\nwach(T)(k))^{\psi=1}_{\Gamma_1}].
\end{multline*}
Le corollaire r\'esulte alors du fait que $\dfont(T(k))^{\psi=1} \simeq H_{\mathrm{Iw}}^1(K,T(k))$.
\end{proof}

\section{La conjecture $C_{\mathrm{Iw}}(K_\infty/K,V)$}\label{nobalt}

Dans ce chapitre, on \'enonce la conjecture $C_{\mathrm{Iw}}(K_\infty/K,V)$ puis on montre qu'elle est \'equivalente $C_{\mathrm{EP}}(K_n/K,V)$ pour tout $n \geq 1$ et finalement, on d\'emontre la conjecture $C_{\mathrm{Iw}}(K_\infty/K,V)$. 

\Subsection{Enonc\'e de la conjecture}\label{nobody}

Dans ce paragraphe, on \'enonce la conjecture $C_{\mathrm{Iw}}(K_\infty/K,V)$ (c'est la conjecture que Perrin-Riou appelle $\delta_{\Zp}(V)$). On commence par des rappels et des compl\'ements sur la loi de r\'eciprocit\'e explicite. 

Soit $(\cdot,\cdot)_{T,n} : H^1(K_n,T) \times H^1(K_n,T^*(1)) \to \Zp$ l'accouplement fourni par la dualit\'e locale. On d\'efinit une application bilineaire :
\[ \pscal{\cdot,\cdot}_T : H^1_{\mathrm{Iw}} (K,T) \times H^1_{\mathrm{Iw}}(K,T^*(1))^{\iota} \to \Lambda, \]
en imposant que pour tout $n \geq 1$, on ait :
\[ \pscal{x,y}_T \equiv \sum_{\tau \in G_n} (\tau^{-1}x_n,y_n)_{T,n}\tau \mod {(\gamma_n -1)}. \]
Par lin\'earit\'e, on obtient un accouplement :
\[ \pscal{\cdot,\cdot}_V : \calK(\Gamma) \otimes_{\Lambda} H^1_{\mathrm{Iw}}(K,T) \times
\calK(\Gamma) \otimes_{\Lambda} H^1_{\mathrm{Iw}}(K,T^*(1))^{\iota} \to \calK(\Gamma). \]
 
D'autre part, en posant $(1+X) \star (1+X)=1+X$, on \'etend la dualit\'e canonique $\ddR^K(V) \times \ddR^K(V^*(1)) \to K$ \`a une forme $\Lambda$-bilin\'eaire :
\[ \star_{\calD(V)} : \calD(V) \times \calD(V^*(1)) \to K \otimes_{\OO_K} \OO_K[[X]]^{\psi=0}. \]
Le th\'eor\`eme suivant est la loi de r\'eciprocit\'e de Perrin-Riou (la conjecture $\mathrm{Rec}(V)$ de \cite{PR1}).

\begin{theo}\label{old322}
Si $V$ est une repr\'esentation cristalline de $G_K$, alors pour tout $h$ on a :
\[ \pscal{\Exp_{V,h}^{\eps}(f),\Exp_{V^*(1),1-h}^{\eps^{-1}}(g^{\iota})}_V (1+X) = (-1)^{h-1}\mathrm{Tr}_{K/\Qp} (f \star_{\calD(V)}g). \]
\end{theo}

On dispose de plusieurs d\'emonstrations de ce r\'esultat : voir  \cite{C1,KKT,Ben,Ber3}. L'approche de \cite{Ben}, qui repose sur la th\'eorie des $(\varphi,\Gamma)$-modules et qui a \'et\'e r\'esum\'ee dans le paragraphe \ref{altexpo}, joue un r\^ole important dans la suite de cet article. 
On note $\Lambda_{\Qp}$ l'anneau $\Qp \otimes_{\Zp} \Lambda$, et on pose : 
\begin{align*}
\Delta_{\mathrm{PR}}(K_{\infty}/K,V) & =
\det_{\Lambda_{\Qp}}\RR\Gamma_{\mathrm{Iw}} (K,V) \otimes \det_{\Lambda_{\Qp}} \calD(V) \\
& \simeq \otimes_{i=1}^2 (\det_{\Lambda_{\Qp}}
H^i_{\mathrm{Iw}}(K,V))^{(-1)^i}\otimes \det_{\Lambda_{\Qp}} \calD(V).
\end{align*}
L'exponentielle de Perrin-Riou induit alors une application :
\[ \det_{\Lambda} (\Exp^{\eps}_{V,h}) : \Delta_{\mathrm{PR}} (K_{\infty}/K,V) \to \calH(\Gamma). \]
Soient $\mathbf{\Gamma}_h(V) = \prod_{j > -h} (\ell_{-j})^{\dim_{\Qp} \Fil^j \dcris (V)}$ et $\delta'_{V,K_{\infty}/K} = \mathbf{\Gamma}_h(V)^{-1} \det_{\Lambda} (\Exp^{\eps}_{V,h})$. Un petit calcul montre que l'application $\delta'_{V,K_{\infty}/K} : \Delta_{\mathrm{PR}} (K_{\infty}/K,V) \to \calK(\Gamma)$ ne d\'epend pas de $h$, et le th\'eor\`eme \ref{old322} entra\^{\i}ne le r\'esultat suivant (c'est l'ancienne conjecture $\delta_{\Qp}(V)$ de \cite{PR1}).

\begin{theo}\label{old324}
On a  $\delta'_{V,K_{\infty}/K}(\Delta_{\mathrm{PR}}(K_{\infty}/K,V)) = \Lambda_{\Qp}$.
\end{theo}

\begin{proof} 
Voir le th\'eor\`eme 3.4.2 et la proposition 3.6.6 de \cite{PR1} .
\end{proof}

Nous allons maintenant donner une version enti\`ere de la conjecture $\delta_{\Qp}(V)$ ci-dessus, c'est la conjecture $\delta_{\Zp}(V)$ de Perrin-Riou. Soient :
\begin{align*}
\Delta_{\mathrm{Iw}}(K_{\infty}/K,T) & =
\det_{\Lambda} \RR\Gamma_{\mathrm{Iw}} (K,T) \otimes_{\Lambda}
\det_{\Lambda} (\Ind_{K_{\infty}/\Qp} T) \\
& \simeq \otimes_{i=1}^2 (\det_{\Lambda}
H^i_{\mathrm{Iw}}(K,T))^{(-1)^i} \otimes_{\Lambda}
\det_{\Lambda} (\Ind_{K_{\infty}/\Qp} T),
\end{align*}
et $\Delta_{\mathrm{Iw}}(K_{\infty}/K,V) =
\Delta_{\mathrm{Iw}}(K_{\infty}/K,T)\otimes_{\Zp}\Qp$.

Soit $a_{V,K} \in \Zp^\times$ l'\'el\'ement d\'efini au paragraphe \ref{fontaine} et soit :
\[ \Lambda_{V,K_\infty/K} = \{ f \in \hat{\OO}_{K^{\mathrm{nr}}} \hat{\otimes}_{\Zp}  \Lambda \mid \sigma(f) = a_{V,K} f \}. \]
On a $\Ind_{K_\infty/\Qp} (V) \simeq (\Lambda \otimes_{\Zp} \Ind_{K/\Qp} (V))^\iota$ et $\calD(V) \simeq \Lambda \otimes_{\Zp} \dcris(V)$. Comme $V$ est cristalline, on a $\dpst(V) =  K^{\mathrm{nr}} \otimes_K \dcris(V)$ et le lemme \ref{old331} donne $\eps(K,V)=1$. L'application $\alpha_{V,K}$ induit donc par lin\'earit\'e un homomorphisme :
\[ \alpha_{V,K_\infty / K} : \det^{-1}_{\Lambda_{\Qp}} \calD(V) \otimes 
\det_{\Lambda_{\Qp}}(\Ind_{K_\infty/\Qp} (V)) \to \Lambda_{V,K_\infty/K} 
\otimes_{\Zp} \Qp. \]
En le composant avec $\delta'_{V,K_\infty/K}$, on obtient une trivialisation canonique : 
\[ \delta_{V,K_{\infty}/K} :  \Delta_{\mathrm{Iw}}(K_{\infty}/K,V) \iso 
\Lambda_{V,K_\infty/K}  \otimes_{\Zp} \Qp. \]

\begin{conj}[$C_{\mathrm{Iw}}(K_\infty/K,V)$]\label{old325}
On a $\delta_{V,K_{\infty}/K}(\Delta_{\mathrm{Iw}}(K_{\infty}/K,T))=
\Lambda_{V,K_\infty/K}$.
\end{conj}

Cette conjecture est d\'emontr\'ee dans le paragraphe \ref{main}, c'est le th\'eor\`eme \ref{new345}.

\Subsection{\'Equivalence de $C_{\mathrm{Iw}}$ et de $C_{\mathrm{EP}}$ : \'etude de $\Xi^{\eps}_{V,n}$}\label{same1}

Ce paragraphe et le suivant sont consacr\'es \`a la d\'emonstration du th\'eor\`eme suivant.

\begin{theo}\label{old326} 
Pour tout $n \geq 1$, la conjecture $C_{\mathrm{Iw}}(K_\infty/K,V)$ est \'equivalente \`a la conjecture $C_{\mathrm{EP}}(K_n/K,V)$.
\end{theo}

Afin de montrer le th\'eor\`eme ci-dessus, nous avons besoin de r\'esultats de descente. La technique g\'en\'erale de descente des complexes a \'et\'e d\'evelopp\'ee par Nekov\`a\v{r} (voir \cite[\S 11.6]{N} ainsi que \cite[lemme 8.1]{BG}). Nous avons besoin d'un cas tr\`es particulier de cette th\'eorie, qui est sans doute bien connu, et qui en tout cas se d\'emontre facilement.

Dans cette section, on pose $\calH=\calH(\Gamma)$ pour all\'eger la notation. Si $M$ et $N$ sont deux $\calH$-modules libres de m\^eme rang et si $f : M \to N$ est un homomorphisme injectif, on note : 
\[ i_{f,\infty} : \det_{\calH} M \otimes \det_{\calH}^{-1} N \to \calH \]
l'homomorphisme qui s'en d\'eduit. Pour tout $n \geq 0$, on a une projection naturelle $\calH \to \Qp[G_n]$. L'alg\`ebre $\Qp[G_n]$ se d\'ecompose en produit de corps $\Qp[G_n] \simeq \oplus_{\lambda} E_\lambda$, et on note $\mathfrak{p} (\lambda)$ le noyau de la projection $\calH \to E_\lambda$. On pose $M_\lambda = E_{\lambda}
\otimes_{\calH} M$, et on note $f_\lambda : M_\lambda \to N_\lambda$ l'homomorphisme de $E_\lambda$-modules qui se d\'eduit de $f$. Le diagramme commutatif :
\[ \begin{CD}
0 @>>> M @>{\gamma_n-1}>> M @>>> M_n @>>> 0 \\
& & @VVfV @VVfV @VV{f_n}V \\
0 @>>> N @>{\gamma_n-1}>> N @>>> N_n @>>> 0
\end{CD} \]
donne lieu \`a des isomorphismes canoniques :
\begin{align*}
\ker (f_n) & \simeq \mathrm{coker} (f)^{\Gamma_n},\\
\mathrm{coker}(f_n) & \simeq \mathrm{coker} (f)_{\Gamma_n}.
\end{align*}

On dit que $f$ est $\lambda$-semi-simple si la $\lambda$-composante de l'application : 
\[ B_{f,n} : \ker (f_n) \hookrightarrow \mathrm{coker}(f) \to
 \mathrm{coker}(f_n) \]
est un isomorphisme. Dans ce cas on a un isomorphisme canonique :
\[ i_{f,\lambda} : \det_{E_\lambda} M_\lambda \otimes 
\det^{-1}_{E_\lambda}N_\lambda \simeq
\det_{E_\lambda} \ker (f_\lambda)
\otimes \det^{-1}_{E_\lambda}\mathrm{coker}(f_\lambda)
\simeq E_\lambda, \]
le deuxieme isomorphisme \'etant induit par $B_{f,n}$.

On dit que $f$ est $\lambda$-admissible si $\mathrm{Im}(i_{f,\infty})$ s'\'ecrit sous la forme $(\gamma_n-1)^r h \calH$, avec $h \in \calH_{\mathfrak{p}(\lambda)}$ et $r = \dim_{E_\lambda} (\ker f_\lambda)$. On dit que $f$ est admissible si elle est $\lambda$-admissible pour tout $\lambda$.

\begin{lemm}\label{old328} 
On conserve les hypoth\`eses concernant $f : M \to N$.  Si $f$ est $\lambda$-admissible, alors $f$ est $\lambda$-semi-simple et le diagramme suivant est commutatif :
\[ \begin{CD}
\det_{{\calH}_{\mathfrak{p}(\lambda)}} M_{\mathfrak{p}(\lambda)} \otimes 
\det_{{\calH}_{\mathfrak{p}(\lambda)}}^{-1} N_{\mathfrak{p}(\lambda)} @>{(\gamma_n-1)^{-r} i_{f,\infty}}>> {\calH}_{\mathfrak{p}(\lambda)} \\
@VVV @VVV \\
\det_{E_\lambda} M_\lambda \otimes
\det_{E_\lambda}^{-1} N_\lambda @>{i_{f,\lambda}}>> E_\lambda.
\end{CD} \]
\end{lemm}

\begin{proof} 
Soit $X_1, \hdots ,X_r \in M_{\mathfrak{p}(\lambda)}$ un rel\`evement d'une base
$x_1,\hdots  x_r$  de $\ker f_\lambda$. Comme $M_\lambda$ est un facteur
direct de $M_n$ on peut choisir $X_i$ de telle fa\c{c}on que
$f(X_i)=(\gamma_n-1)A_i$ o\`u $A_i\in N_{\mathfrak{p}(\lambda)}$.
On fixe un compl\'ement  $X_{r+1},\hdots , X_m$
de $X_1,\hdots ,X_r$  \`a une base de $M_{\mathfrak{p}(\lambda)}$. Soit 
$Y_1,\hdots,Y_r\in N_{\mathfrak{p}(\lambda)}$ un rel\`evement d'une base $y_1,\hdots, y_r$ de $\mathrm{coker}(f_{\lambda})$. Les \'el\'ements $Y_1,\hdots ,Y_r,
Y_{r+1}=f(X_{r+1}),\hdots ,Y_m=f(X_m)$ forment alors une base de $N_{\mathfrak{p}(\lambda)}$ et on a :
\[ i_{f,\infty} (\wedge_{i=1}^m X_i \otimes    
\wedge_{j=1}^m Y_j^* ) = \det (\wedge_{i=1}^r X_i \otimes    
\wedge_{j=1}^r Y_j^* ) =(\gamma_n -1)^r \det_{1\leq i,j \leq r} (Y_j^*(A_i)), \]
ce qui montre le lemme.
\end{proof}

En particulier, soit $N$ un $\Lambda_{\Qp}$-module de torsion et de type fini; il admet une r\'esolution projective $0 \to P_1 \overset{f}{\to} P_0 \to N \to 0$, 
o\`u $ \mathrm{rg}(P_0)=\mathrm{rg}(P_1)$. Pour tout $\lambda$, on note $\Lambda_{\mathfrak{p}(\lambda)}$ la localisation de $\Lambda_{\Qp}$ en $\mathfrak{p}(\lambda)$. On dit que $N$ est $\lambda$-admissible si $f : P_1 \to P_0$ l'est. Dans ce cas, le lemme \ref{old328} fournit un diagramme commutatif :
\[ \begin{CD}
\det_{\Lambda_{\mathfrak{p}(\lambda)}}^{-1} N_{\mathfrak{p}(\lambda)} @>{(\gamma_n -1)^{-r} i_{f,\infty}}>> \Lambda_{\mathfrak{p}(\lambda)} \\
@VVV @VVV \\
\det_{E_ \lambda}^{-1} {(N_{\Gamma_n})}_ \lambda \otimes
\det_{E_ \lambda} {(N^{\Gamma_n})}_ \lambda  @>>{\mathrm{id}}> E_ \lambda, 
\end{CD} \]
o\`u la deuxi\`eme ligne est induite par la projection $N^{\Gamma_n} \to N_{\Gamma_n}$.

On pose maintenant $\Lambda_{\Gamma_1}=\Zp[[\Gamma_1]]$, et on fixe un isomorphisme $\Lambda_{\Gamma_1} \simeq \Zp[[T]]$ en envoyant $\gamma_1$ sur $1+T$. On pose :
\[ \delta_i =\frac{1}{\sharp \Delta_K} \sum_{g \in \Delta_K} \chi^i(g^{-1}) g, \] 
ce qui fait que $\Lambda = \oplus_{i=0}^{p-2} \Lambda_i$,  o\`u $\Lambda_i = \delta_i \Lambda_{\Gamma_1}$.

Il est clair que $\delta_i (\Zp(j))=0$ si $i \neq j \mod {(p-1)}$. Sinon, on a une suite exacte $0 \to \Zp[[T]] \overset{f_j}{\to} \Zp[[T]] \to \Zp(j) \to 0$, o\`u  $f_j$ est la multiplication par $(\chi(\gamma_1)^j-1)-T$, ce qui fait que $\det^{-1}_{\Lambda_i} \Zp(j) = ((\chi(\gamma_1)^j-1)T) \Lambda_i$ 
si $i = j \mod {(p-1)}$.

\begin{prop}\label{old3211}
Si $h \geq 1$ est un entier tel que $\Fil^{-h} \ddR^K(V) = \ddR^K(V)$, alors l'application $\Exp_{V,h}^{\eps}$ est admissible.
\end{prop}

Nous allons d\'eduire cette proposition du th\'eor\`eme \ref{old324}. Pour all\'eger les notations, posons $a_j = \dim_{\Qp} \Fil^j \ddR^K(V)$, $b_j=\dim_{\Qp} \dcris(V)^{\varphi=p^{-j}}$ et $\omega_k(T) = (1+T)^{p^k} - 1$. On commence par un lemme purement technique.

\begin{lemm}\label{old3212}
Pour tout $n \geq 1$, on a $\mathbf{\Gamma}_h(V) \equiv  {\mathbf{\Gamma}}^*_h(V) \omega_{n-1}(T)^{a_0}  \mod{\omega_{n-1}(T)^{a_0+1}}$, o\`u $\mathbf{\Gamma}_h^*(V) = \pm (h-1)!^{\dim_{\Qp} \dcris (V)} (\log \chi (\gamma_n))^{-a_0}\Gamma^*(V)^{-1}$.
\end{lemm}

\begin{proof}
Comme :
\[ \frac{\log (\gamma)}{\log \chi (\gamma)}=\frac{\log (\gamma_n)}{\log \chi (\gamma_n)}
\equiv \log^{-1} \chi (\gamma_n) \omega_{n-1}(T) \mod{\omega_{n-1}(T)^2}, \]
on a :
\[ \mathbf{\Gamma}_h (V) \equiv (-1)^{a_0} (\log \chi (\gamma_n))^{-a_0} 
 \prod_{\substack{j>-h \\ j \neq 0}} j^{a_j}  \omega_{n-1}(T)^{a_0} 
\mod{\omega_{n-1}(T)^{a_0+1}}. \]
Un calcul facile montre alors que :
\begin{align*}
(h-1)!^{\dim_{\Qp} \dcris (V)} & = \pm \prod_{\substack{j>-h \\ j \neq 0}}
 j^{a_j}  \prod_{i \geq 0}
\Gamma^*(h-i)^{[K:\Qp] h_{i-h}(V)} \\
& = \pm  \prod_{\substack{j>-h \\ j \neq 0}} j^{a_j} \Gamma^*(V),
\end{align*}
d'o\`u le lemme.
\end{proof}

\begin{proof}[D\'emonstration de la proposition \ref{old3211}] 
Il r\'esulte de la suite exacte (eq4) que l'image de :
\[ \det_{\Lambda_{\Qp}} \calD(V)^{\Delta=0} \otimes
\det_{\Lambda_{\Qp}}^{-1}
\left ( \frac{H^1_{\mathrm{Iw}}(K,V)}{V^{H_K}} \right) \]
dans $\calH(\Gamma)$ est \'egale \`a :
\[ v= \mathbf{\Gamma}_h (V) \det_{\Lambda_{\Qp}}(V^{H_K})  
\det^{-1}_{\Lambda_{\Qp}}(V^*(1)^{H_K})^*  \prod_{j \in \ZZ} 
\left(\det^{-1}_{\Lambda_{\Qp}}\Qp(j) \right)^{b_j}. \]

Fixons un $\lambda$ et notons $n$ le plus petit entier tel que $E_{\lambda} \subset \Qp[G_n]$.

Supposons d'abord que $\lambda\neq\lambda_0$ ($\lambda_0$ correspond \`a l'inclusion de $\Qp$ dans $\Qp[G_n]$). Alors 
$\det_{\Lambda_{\Qp}}(V^{H_K})$, 
$\det_{\Lambda_{\Qp}}(V^*(1)^{H_K})^*$
et $\det_{\Lambda_{\Qp}}\Qp(j)$
sont des unit\'es de $\calH_{\mathfrak{p}(\lambda)}$ et $v$ est congru \`a :
\[ \mathbf{\Gamma}_h^*(V)
\det_{\Lambda_{\Qp}}(V^{H_K})  
\det^{-1}_{\Lambda_{\Qp}}(V^*(1)^{H_K})^* \prod_{j \in \ZZ} \det_{\Lambda_{\Qp}}\Qp(j)^{-b_j}
\omega_{n-1}(T)^{a_0} \mod{\omega_{n-1}^{a_0+1}(T)}. \]
D'autre part, la suite exacte (eq5) nous donne $\ker (\Exp_{V,h,n}^{\eps})_{\lambda} \simeq \Fil^0 \ddR^{K_n}(V)_{\lambda}$, ce qui fait que $\dim_{E_{\lambda}}(\ker (\Exp_{V,h,n}^{\eps})_{\lambda})=a_0$
et $\Exp_{V,h}^{\eps}$ est bien $\lambda $-semi-simple.

Supposons maintenant que $\lambda = \lambda_0$. Le m\^eme calcul montre alors que $v$ est congru \`a : 
\[ \mathbf{\Gamma}_h^*(V) \det_{\Lambda_{\Qp}}(V^{H_K}/V^{G_K})  
\det^{-1}_{\Lambda_{\Qp}}(V^*(1)^{H_K}/V^*(1)^{G_K})^* 
\prod_{j\in \ZZ} \det_{\Lambda_{\Qp}}\Qp(j)^{-b_j}
T^r \mod{T^{r+1}}, \]
o\`u :
\begin{align*}
r & = \dim_{\Qp}(\ker (\Exp_{V,h,n}^{\eps})_{\lambda_0}) \\
& = \dim_{\Qp} (\dcris(V)^{\varphi=1})+
\dim_{\Qp} (\Fil^0 \dcris (V))-
\dim_{\Qp} (V^{G_K})+
\dim_{\Qp} (V^*(1)^{G_K}),
\end{align*}
et $\Exp_{V,h}^{\eps}$ est bien $\lambda$-semi-simple dans ce cas aussi.
\end{proof}

Nous allons maintenant calculer le d\'eterminant de l'application $\Xi^{\eps}_{V,n}$; ces calculs g\'en\'eralisent (et corrigent\dots) ceux de \cite[lemme 3.5.7]{PR1}. Si : 
\[ \Delta_0 : \calD(V)\to\dcris (V)/(1-\varphi)\dcris (V) \] 
est l'application d\'efinie ci-dessus par $\Delta_0(\alpha (X))=\alpha (0) \mod {(1-\varphi)\dcris (V)}$, alors on voit facilement que 
$\calD(V)^{\Delta_0=0}_{\Gamma_n}=\calD(V)^{\Delta=0}_{\Gamma_n}$
et que $\Delta_0$  induit une suite exacte :
\[ 0\to \frac{\dcris (V)}{(1-\varphi)\dcris (V)} \to(\calD(V)^{\Delta=0})_{\Gamma_n}
\to\calD(V)_{\Gamma_n}\to\frac{\dcris (V)}{(1-\varphi)\dcris (V)}\to 0. \]

\begin{lemm}\label{old342}
On a :
\[ \tilde{\Xi}^{\eps}_{V,n} (\alpha (X))=
p^{-n} \left(\sum_{k=1}^n (\sigma \otimes \varphi)^{-k}\alpha (\zeta_{p^k}-1)+
(1-\varphi)^{-1}\alpha(0) \right )
\mod {\dcris (V)^{\varphi=1}}. \]
\end{lemm}

\begin{proof} 
Si $\alpha (X)\in \calD(V)^{\Delta_0=0}$,  alors $\alpha (0)\in (1-\varphi)\dcris (V)$ et l'\'el\'ement
$(1-\varphi)^{-1} \alpha (0) \in \dcris (V)/\dcris (V)^{\varphi=1}$
est bien d\'efini. Soit $F(X) \in \calH \otimes_K \dcris (V)$ telle que
$(1-\varphi )F(X)=\alpha (X)$. Pour tout $0 \leq k \leq n-1$, on a :
\[ (\varphi\otimes \sigma)^k F(\zeta_{p^{n-k}}-1)-
(\varphi\otimes \sigma)^{k+1} F(\zeta_{p^{n-k-1}}-1)=
(\varphi\otimes \sigma)^k \alpha (\zeta_{p^{n-k}}-1) \]
ainsi que $(1-\varphi)F(0)=\alpha (0)$, ce qui fait que : 
\begin{align*}
\tilde{\Xi}^{\eps}_{V,n} (\alpha (X))&=p^{-n}(\varphi\otimes \sigma)^{-n}
F(\zeta_{p^n}-1)\\
&= p^{-n} \left (\sum_{k=1}^n (\varphi \otimes \sigma)^{-k}\alpha (\zeta_{p^k}-1)+ (1-\varphi)^{-1}\alpha (0) \right ) 
\mod {\dcris (V)^{\varphi=1}}.
\end{align*}
\end{proof}

Le (2) de la proposition \ref{old316} donne une suite exacte :
\[ 0\to\ker \tilde{\Xi}_{V,n}^{\eps}
\to \calD(V)_{\Gamma_n}^{\Delta_0=0}\overset{\tilde{\Xi}^{\eps}_{V,n}}{\longrightarrow} \frac{\ddR^{K_n}(V)}{\dcris(V)^{\varphi=1}}
\to \frac{\dcris (V)}{(1-p^{-1}\varphi^{-1})\dcris (V)}\to 0. \]
En composant cette suite avec les suites tautologiques :
\begin{align*}
&0\to\dcris (V)^{\varphi=1}\to\dcris (V)\overset{1-\varphi}{\longrightarrow} \dcris (V)\to
\frac{\dcris (V)}{(1-\varphi)\dcris (V)}\to 0,\\ \tag{eq6} 
&0\to\dcris (V)^{\varphi=p^{-1}}\to\dcris (V) \overset{1-p^{-1}\varphi^{-1}}{\longrightarrow} \dcris (V)\to
\frac{\dcris (V)}{(1-p^{-1}\varphi^{-1})\dcris (V)}\to 0,
\end{align*}
et en utilisant le (1) de la proposition \ref{old316}, on obtient un isomorphisme canonique :
\begin{multline*}
\kappa_{V,n}:\det_{\Qp[G_n]}\calD(V)^{\Delta=0}_{\Gamma_n}
\otimes \det^{-1}_{\Qp[G_n]} \ddR^{K_n}(V)\simeq \\
 \left (\det_{\Qp[G_n]}^{-1} \dcris (V)^{\varphi=1} \otimes \det_{\Qp[G_n]} \frac{\dcris (V)}{(1-\varphi)\dcris (V)}\right )
\\Ê\otimes \left (
\det_{\Qp[G_n]} \dcris (V)^{\varphi=p^{-1}} \otimes \det_{\Qp[G_n]}^{-1} \frac{\dcris (V)}{(1-p^{-1}\varphi^{-1})\dcris (V)}
\right ) 
\simeq \Qp[G_n]. 
\end{multline*}

Rappelons que l'on note $R_n$ le $\OO_K[G_n]$-module libre engendr\'e
par $x_n=\zeta_p+\cdots+\zn$ et que l'on appelle $\eta_0$ le caract\`ere trivial. On fixe un r\'eseau $M$ de $\dcris (V)$ et l'on pose $\calD_M(V)=\OO_K[[X]]^{\psi=0} \otimes M$ et $M_n=R_n \otimes_{\Zp} M$. 

\begin{prop}\label{old344} 
L'isomorphisme $\kappa_{V,n}$ envoie $\det_{\Zp[G_n]} \calD_M(V)^{\Delta_0=0}_{\Gamma_n}
\otimes\det^{-1}_{\Zp[G_n]} M_n$ sur le r\'eseau engendr\'e par : 
\begin{multline*} p^{-n d} \sum_{\eta\neq\eta_0} (\eta (\gamma_1)-1)^{\delta (\eta)\dim_{\Qp}\dcris (V)^{\varphi=1}} \det (\varphi  \mid \dcris (V))^{-a(\eta)} 
e_{\eta} 
\\ + (-1)^{d} p^{(1-n)d+\dim_{\Qp}\dcris (V)^{\varphi=1}} e_{\eta_0}, 
\end{multline*}
o\`u $\delta (\eta)= 1$ si $\eta \in X(\Gamma_1)$ et $\delta (\eta)=0$ sinon.
\end{prop}

\begin{proof}
Nous allons commencer par montrer qu'il suffit de d\'emontrer la proposition caract\`ere par caract\`ere. Posons $\Lambda'  = T \Lambda_0 \oplus (\oplus_{i=1}^{p-2} \Lambda_i)$ ce qui fait que $(X\Zp[[X]])^{\psi=0} = \Lambda'(1+X)$ et donc que $(X\Zp[[X]])^{\psi=0}$ est un $\Lambda$-module libre engendr\'e par
$\lambda  = T \delta_0 (1+X)+\sum_{i=1}^{p-2} \delta_i (1+X)$.

Si $N=(1-\varphi)\dcris (V) \cap M$, alors $M/N$ est sans torsion et il existe $N' \subset M$ tel que $M= N \oplus N'$. On a :
\[ \calD_M(V)^{\Delta_0=0}\simeq  (\Zp[[X]]^{\psi=0} \otimes_{\Zp} N)\oplus
((X\Zp[[X]])^{\psi=0} \otimes_{\Zp} N'). \]

Soient $A$ et $B$ deux $\Zp[G_n]$-modules libres de m\^eme rang et supposons que l'on se donne un isomorphisme :
\[ \kappa : \det_{\Qp[G_n]} (A_{\Qp}) \otimes \det_{\Qp[G_n]}^{-1} (B_{\Qp}) \to\Qp[G_n]. \]
Soit $E$ une extension de $\Qp$ contenant toutes les valeurs des caract\`eres de $G_n$. Pour tout $\eta \in X(G_n)$, posons $A_\eta = e_\eta (\OO_E \otimes_{\Zp} A)$ et $\OO_\eta = e_\eta \OO_E$. Pour tous $a \in \det_{\Zp[G_n]}A$ et $b\in \det_{\Zp[G_n]}B$, on a alors $e_\eta (\kappa (a\otimes b^{-1}))=\kappa_\eta (e_\eta (a)\otimes e_\eta (b)^{-1})$. 
On en conclut qu'il suffit de d\'emontrer la proposition \ref{old344} caract\`ere par caract\`ere.

On fixe une base $(n_i)$ (resp. $(n'_j)$) de $N$ (resp. $N'$) et l'on pose :
\begin{alignat*}{2}
&\alpha_i=n_i\otimes (1+X) & \qquad &\beta_i=n_i\otimes x_n\\
&\alpha'_j=n'_j\otimes \lambda  & &\beta'_j=n'_j \otimes x_n\\
&\tilde{\alpha}=\wedge_i \alpha_i & &\tilde{\beta} = \wedge_j \beta_j\\
&\tilde{\alpha}'=\wedge_i \alpha'_i & &\tilde{\beta}' = \wedge_j \beta'_j\\
&\tilde{a} = \tilde{\alpha}  \wedge \tilde{\alpha}' &&
\tilde{b} = \tilde{\beta}  \wedge \tilde{\beta}'. 
\end{alignat*}

Si $\eta \neq \eta_0$ est un caract\`ere de conducteur $p^k$, alors il s'\'ecrit sous la forme $\eta = \omega^i \eta '$, o\`u $\omega$ est le caract\`ere de Teichm\"uller et
$\eta' \in X(\Gamma_1/\Gamma_k)$. Si $\eta \not \in X(\Gamma_1)$, alors $i \neq 0$, $e_\eta (\lambda) = e_\eta (1+X)$, et la proposition \ref{old342} nous donne :
\begin{align*}
e_\eta \left ( \tilde{\Xi}^{\eps}_{V,n} (\alpha_i) \right)&=p^{-n}\varphi^{-k} (n_i) 
e_\eta (\zeta_{p^k}),\\
e_\eta \left ( \tilde{\Xi}^{\eps}_{V,n} (\alpha'_j )\right )&=p^{-n}\varphi^{-k} (n'_j) 
e_\eta (\zeta_{p^k}).
\end{align*}

Comme $e_\eta (x_n) = e_\eta(\zeta_{p^k})$, on a pour $\eta \not \in X(\Gamma_1)$ :
\[ \kappa_{V,\eta} (e_{\eta} (\tilde{a} \otimes \tilde{b}^{-1})) =
p^{-n \dim_{\Qp}\dcris (V)} \det_{\Qp} (\varphi^{-k}  \mid
\dcris (V)), \]
tandis que si $\eta \in X(\Gamma_1)$, alors on a :
\begin{align*}
e_\eta (\lambda) & =
(\eta (\gamma_1)-1)e_\eta (1+X), \\
 e_\eta \left ( \tilde{\Xi}^{\eps}_{V,n} (\alpha'_j )\right ) & = p^{-n}
(\eta (\gamma_1)-1)  \varphi^{-k} (n'_j) 
e_\eta (\zeta_{p^k}), \tag{eq7} \\
\kappa_{V,\eta} \left (e_{\eta} (\tilde{a} \otimes \tilde{b}^{-1}) \right ) & =
p^{-n \dim_{\Qp}\dcris (V)}
(\eta (\gamma_1)-1)^{\dim_{\Qp}\dcris (V)^{\varphi=1}}
 \det_{\Qp} (\varphi^{-k}  \mid \dcris (V)). 
\end{align*}

Supposons maintenant que $\eta =\eta_0$. Dans ce cas, $e_{\eta_0}(\zeta_p)=(1-p)^{-1}$ et la proposition \ref{old342} nous donne :
\begin{align*}
e_{\eta_0} \left ( \tilde{\Xi}^{\eps}_{V,n} (\alpha_i) \right)&=
\frac{1}{[K_n:K]}
(1-p^{-1}\varphi^{-1})(1-\varphi)^{-1}n_i,
\\
e_{\eta_0} \left ( \tilde{\Xi}^{\eps}_{V,n} (\alpha'_j )\right )&=0.
\end{align*}

On a un diagramme commutatif :
\[ \xymatrix{& & 0 \ar[d]& 0\ar[d] &\\
& & \ker (\tilde{\Xi}_{V,n}^{\eps}) \ar[r]^{(2)}\ar[d] & \dcris (V)^{\varphi=p^{-1}} \ar[r]
\ar[d]&0\\
0\ar[r]&\frac{\dcris (V)}{(1-\varphi)\dcris(V)}\ar[ur]^{(1)}\ar[r]^-{(3)}
&\left (\calD(V)^{\Delta_0=0}_{\Gamma_n}\right )_{\eta_0} \ar[r]^{(4)}
\ar[d]
&(1-\varphi)\dcris (V)\ar[r]
\ar[dl]^{(5)}
& 0\\
& & \frac{\dcris(V)}{\dcris (V)^{\varphi=1}}\ar[d]^{(6)} & &\\
& & \frac{\dcris(V)}{(1-p^{-1}{\varphi}^{-1})\dcris (V)}\ar[d] & &\\
& & 0 & &} \]
dont les fl\`eches sont donn\'ees par les formules suivantes :

\begin{itemize}
  \item[(1) et (3)] : $d\mapsto d\otimes (\gamma_n-1)(1+X)$;
  \item[(2) et (4)] : $\alpha (X)\mapsto \alpha (0)$;
  \item[(5)] : $d \mapsto [K_n:K]^{-1} (1-p^{-1}\varphi^{-1})(1-\varphi)^{-1}(d)
\mod{\dcris(V)^{\varphi=1}}$;
  \item[(6)]  : $d\mapsto [K_n:K](1-\varphi)(d)$.
\end{itemize}

On en d\'eduit que la $\eta_0$-composante de $\kappa_{V,n}$ s'\'ecrit comme le compos\'e des isomorphismes suivants; par les suites (eq3) et (eq4), $\det_{\Qp}\left ( \calD(V)^{\Delta_0=0}_{\Gamma_n}\right )_{\eta_0} \otimes \det_{\Qp}^{-1} \dcris (V)$ est isomorphe \`a :
\begin{multline*} 
\left (\det_{\Qp} \left ( \frac{\dcris (V)}{(1-\varphi)\dcris (V)}\right ) \otimes
\det_{\Qp} ((1-\varphi)\dcris(V)) \right ) \\ 
\otimes 
\left (\det_{\Qp}^{-1} \left (\frac{\dcris (V)}{\dcris (V)^{\varphi=1}}\right ) \otimes
\det_{\Qp}^{-1} \dcris (V)^{\varphi=1}\right ), 
\end{multline*}
qui est isomorphe \`a :
\begin{multline*} 
\left ( \det_{\Qp} ((1-\varphi)\dcris (V)) \otimes
\det_{\Qp}^{-1} \left( \frac{\dcris (V)}{\dcris (V)^{\varphi=1}}\right) \right) \\
\otimes
\left( \det_{\Qp}^{-1} \dcris (V)^{\varphi=1} \otimes
\det_{\Qp} \left( \frac{\dcris (V)}{(1-\varphi)\dcris (V)}\right)\right),
\end{multline*}
qui, par les suites (eq5) et (eq6), est lui-m\^eme isomorphe \`a :
\begin{multline*} 
\left (\det_{\Qp}\dcris (V)^{\varphi=p^{-1}} \otimes
\det_{\Qp}^{-1} \left( \frac{\dcris (V)}{(1-p^{-1}\varphi^{-1})\dcris (V)} \right) \right) \\
\otimes \left (\det^{-1}_{\Qp}\dcris (V)^{\varphi=1} \otimes
\det_{\Qp} \left( \frac{\dcris (V)}{(1-\varphi)\dcris (V)}\right)\right), 
\end{multline*}
qui est $\simeq \Qp$.

On a $(\calD(V)^{\Delta_0=0})_{\eta_0} \simeq (1-\varphi )\dcris (V) \oplus D'$, o\`u $D'= \Qp \otimes_{\Zp} N'$ et la deuxi\`eme ligne du gros diagramme ci-dessus s'\'ecrit donc :
\[ 0\to\frac{\dcris (V)}{(1-\varphi) \dcris (V)} \overset{d \mapsto (0,p^{n-1}d)}{\longrightarrow} 
(1-\varphi)\dcris (V)\oplus D' \overset{(a,b)\mapsto a}{\longrightarrow} 
(1-\varphi)\dcris (V) \to 0. \tag{eq8} \]

Posons \`a pr\'esent $\tilde{n}=\wedge_i n_i \in \det_{\Qp}(1-\varphi)
\dcris (V)$ et $\tilde{n}'=\wedge_j n_j' \in \det_{\Qp} D'$. Fixons $\tilde{m}_1 \in \det_{\Qp} M^{\varphi=1}$ et $\tilde{m}_2 \in \det_{\Qp} \calD_M(V)^{\Delta_0=0}_{\Gamma_n} (M/M^{\varphi=1})$, v\'erifiant $\tilde{m}_1 \otimes \tilde{m}_2 \simeq \tilde{n} \otimes \tilde{n}'$, et posons  $\tilde{b}_i = \tilde{m}_i \otimes x_n$. En utilisant le diagramme tautologique :
\[ \xymatrix{& 0 &0 & \\
0\to \dcris (V)^{\varphi=p^{-1}}\ar[r]&
(1-\varphi)\dcris (V)\ar[u]\ar[r]^{(5)} &\frac{\dcris (V)}{\dcris (V)^{\varphi=1}}\ar[u]\ar[r]
&\frac{\dcris (V)}{(1-p^{-1}\varphi^{-1})\dcris (V)}\to 0\\
0\to\dcris (V)^{\varphi=p^{-1}}\ar[u]^{p^{-n}}\ar[r]&
\dcris (V)\ar[u]^{[K_n:K](1-\varphi)}\ar[r]^{1-p^{-1}\varphi^{-1}} &\dcris (V)
\ar[u]^{\mathrm{pr}}\ar[r]
&\frac{\dcris (V)}{(1-p^{-1}\varphi^{-1})\dcris (V)}\ar[u]^{p^{-n}}\to 0\\
 & \dcris (V)^{\varphi=1} \ar[u] \ar[r]^{1-\frac{1}{p}}&\dcris (V)^{\varphi=1}\ar[u] & \\
 & 0\ar[u]&0\ar[u]&  } \]
on montre que l'image de $e_{\eta_0}(\tilde{\alpha}) \otimes e_{\eta_0}^{-1}(\tilde{b}_2)$ par l'application :
\begin{multline*}
\det_{\Qp}(1-\varphi )\dcris (V) \otimes 
\det_{\Qp}^{-1}\left (\frac{\dcris (V)}{\dcris (V)^{\varphi=1}}\right )\simeq \\
\det_{\Qp}\dcris (V)^{\varphi =p^{-1}} \otimes 
\det_{\Qp}^{-1}\frac{\dcris (V)}{(1-p^{-1}\varphi^{-1})\dcris (V)}
\simeq \\
\det_{\Qp}^{-1}\dcris (V) \otimes \det_{\Qp}\dcris (V) \overset{\mathrm{id}}{\simeq} \Qp
\end{multline*}
est \'egale \`a :
\[ (-1)^{\dim_{\Qp}\dcris (V)}
\left (\frac{p}{p-1}\right )^{\dim_{\Qp} \dcris (V)^{\varphi=1}}
p^{(1-n)\dim_{\Qp} (\dcris (V)/\dcris (V)^{\varphi=1})} [(1-\varphi)\tilde{m}_2 :\tilde{n}]. \]
D'autre part, en utilisant (eq8), on montre que l'application compos\'ee : 
\begin{multline*}
\det_{\Qp}^{-1}\dcris (V)^{\varphi=1} \otimes \det_{\Qp} D' 
\simeq 
\det_{\Qp}^{-1}\dcris (V)^{\varphi=1}  \otimes 
\det_{\Qp}\left (\frac{\dcris (V)}{(1-\varphi)\dcris (V)}\right )
\\ \simeq  \det_{\Qp}^{-1}\dcris (V) \otimes \det_{\Qp}\dcris (V)
\overset{\mathrm{id}}{\simeq} \Qp
\end{multline*}
envoie $e_{\eta_0}^{-1}(\tilde{b}_1) \otimes e_{\eta_0}(\tilde{\alpha}')$ sur : 
\[ \left (\frac{1-p}{p^{n-1}}\right )^{\dim_{\Qp} \dcris (V)^{\varphi=1}} [(1-\varphi)\tilde{m}_2 :\tilde{n}]^{-1}. \]

On en conclut finalement que l'application $\kappa_{V,n}$ envoie $e_{\eta_0}(\tilde{a}) \otimes e_{\eta_0}(\tilde{b})^{-1}$
sur :
\[ (-1)^{\dim_{\Qp}\dcris (V)} p^{(1-n)\dim_{\Qp} \dcris (V) + 
\dim_{\Qp} \dcris (V)^{\varphi=1}} e_{\eta_0}, \]
ce qui termine enfin la d\'emonstration.
\end{proof}

\begin{prop}\label{old349}
L'image de $\det_{\Zp[G_n]} \calD_M(V)_{\Gamma_n}
\otimes \det_{\Zp[G_n]}^{-1} M_n$ dans $\Qp[G_n]$ est engendr\'ee par :
\[ p^{-nd} \sum_{\eta \neq \eta_0} \det ( \varphi  \mid \dcris (V))^{-a (\eta)} 
e_\eta + (-1)^{d} p^{(1-n)d + 
n \dim_{\Qp}\dcris (V)^{\varphi=1}} e_{\eta_0}. \]
\end{prop}

\begin{proof} 
Il suffit de calculer l'image de : 
\[ \det_{\Zp[G_n]} \calD_M(V)^{\Delta_0=0}_{\Gamma_n}
\otimes \det^{-1}_{\Zp[G_n]} \calD_M(V)_{\Gamma_n} \]
dans $\Qp[G_n]$ et d'utiliser la proposition \ref{old344}. 

Si $\eta \in X(\Gamma_1)$ est un caract\`ere non-trivial, alors on a 
$(\calD(V)^{\Delta_0=0}_{\Gamma_n})_{\eta}
\simeq (\calD_M(V)^{\Delta_0=0}_{\Gamma_n})_{\eta}$
et la formule (eq7) montre que l'image de $\det_{\OO_\eta} (\calD_M(V)^{\Delta_0=0}_{\Gamma_n})_{\eta}
\otimes \det^{-1}_{\OO_\eta} (\calD_M(V)_{\Gamma_n})_{\eta}$
est engendr\'ee par $(\eta (\gamma_1)-1)^{\dim \dcris (V)^{\varphi =1}} 
e_\eta$.

Supposons maintenant que $\eta=\eta_0$. La suite exacte
$0\to \Lambda_{\Gamma_1} \overset{T}{\to} \Lambda_{\Gamma_1}\to \Zp 
\to 0$ induit une suite exacte
\[ 0 \to \Zp \overset{s}{\to} \Zp[\Gamma_1/\Gamma_{m}] \overset{T}{\to} \Zp[\Gamma_1/\Gamma_{m}]\to\Zp \to 0, \]
o\`u $s$  est la multiplication par $((1+T)^{p^{m}}-1)/T$. On en d\'eduit que l'application :
\begin{multline*}
\det_{\Qp[\Gamma_1/\Gamma_{m}]}(\Qp[\Gamma_1/\Gamma_{m}])
\otimes
\det_{\Qp[\Gamma_1/\Gamma_{m}]}^{-1}
(\Qp[\Gamma_1/\Gamma_{m}])
\simeq \\
\det_{\Qp[\Gamma_1/\Gamma_{m}]}(\Qp)
\otimes
\det_{\Qp[\Gamma_1/\Gamma_{m}]}^{-1}(\Qp) 
\overset{\mathrm{id}}{\simeq}
\Qp[\Gamma_1/\Gamma_{m}]
\end{multline*}
envoie $\det_{\Zp} (\Zp e_{\eta_0})\otimes
\det_{\Zp}^{-1} (\Zp e_{\eta_0})$ sur $p^{-m}\Zp$. En posant $m=n-1$ on obtient que l'image de $\det_{\Zp}  (\calD_M(V)^{\Delta_0=0}_{\Gamma_n})_{\eta_0} \otimes
\det^{-1}_{\Zp} (\calD_M(V)_{\Gamma_n})_{\eta_0}$ est engendr\'ee par $p^{(1-n) \dim_{\Qp}\dcris (V)^{\varphi=1}}$, d'o\`u la proposition.
\end{proof}

Ceci termine notre \'etude de $\Xi^{\eps}_{V,n}$. 

\Subsection{\'Equivalence de $C_{\mathrm{Iw}}$ et de 
$C_{\mathrm{EP}}$ : \'etude de $\Exp^{\eps}_{V,h,n}$}\label{same2}

Nous passons maintenant \`a l'\'etude de l'application :
\[ \Exp^{\eps}_{V,h,n} : (\calD(V)^{\Delta=0})_{\Gamma_n}\to
\Qp \otimes_{\Zp} (H_{\mathrm{Iw}}^1(K,T)/T^{H_K})_{\Gamma_n} \]
d\'eduite de $\Exp^{\eps}_{V,h}$. 

\begin{lemm}\label{old352} 
Si $V$ est une repr\'esentation cristalline, alors : 
\begin{enumerate}
  \item L'application naturelle de $V^{G_K}$ dans $H^1(\Gamma_n, V^{H_K})$ est un isomorphisme;  
  \item L'application compos\'ee :
\[ V^{G_K} \to H^1(K,V) \overset{\exp^*_{V,K}}{\longrightarrow} \ddR^{K_n}(V) \]
co\"{\i}ncide avec l'injection $V^{G_K} \overset{\log^{-1}\chi (\gamma_n)}{\longrightarrow} \Fil^0 \ddR^{K_n}(V)$;
  \item On a $V^{G_K} \cap H^1_g(K_n,V) = \{0\}$.
\end{enumerate}
\end{lemm}

\begin{proof} 
Comme $V$ est cristalline, on a un isomorphisme $V^{H_K} \simeq \oplus_{i \in \ZZ} \Qp(i)^{d_i}$ (voir \cite[lemme 3.4.3]{PR1}). La premi\`ere assertion
s'en d\'eduit. 

Pour montrer la deuxi\`eme, on remarque que l'application $\cup \log \chi : \Fil^0 \ddR^{K_n}(V) \to H^1(K,\Fil^0\bdR \otimes V)$ est un isomorphisme et que $\exp_{V,K_n}^*$ co\"{\i}ncide avec l'application compos\'ee : 
\[ H^1(K_n,V) \to H^1(K_n,\Fil^0\bdR\otimes V) \iso \Fil^0 \ddR^{K_n}(V) \]
(c'est la formule de Kato, voir \cite[\S 1.2-1.4]{K1}). 

Enfin, comme $\ker (\exp^*_{V,K_n}) = H^1_g(K_n,V)$, on en d\'eduit le (3). 
\end{proof}

Comme $V^*(1)^{G_K}$ est isomorphe \`a $\Fil^0 \dcris(V^*(1))^{\varphi=1}$, on a un isomorphisme :
\[ \frac{\dcris (V)}{\Fil^0 \dcris(V)+(1-p^{-1}\varphi^{-1})\dcris (V)}
\simeq (V^*(1)^{G_K})^*, \]
d'o\`u la suite exacte courte :
\[  0\to\frac{\calD(V)^{\Delta=0}_{\Gamma_n}}{\ker (\Exp^{\eps}_{V,h,n})} \overset{\tilde{\Xi}^{\eps}_{V,n}}{\longrightarrow}
\frac{t_V(K_n)}{\dcris (V)^{\varphi=1}/V^{G_K}} \overset{(1-\varphi)\mathrm{Tr}_{K_n/K}}{\longrightarrow} (V^*(1)^{G_K})^* \to 0. \]
D'autre part, pour toute repr\'esentation $p$-adique, on a une suite exacte :
\[ 0 \to \left (\frac{H^1_{\mathrm{Iw}}(K,V)}{V^{H_K}} 
\right )_{\Gamma_n} \to \frac{H^1(K_n,V)}{ H^1(\Gamma_n, V^{H_K})} \to
H^2_{\mathrm{Iw}}(K,V)^{\Gamma_n} \to 0 \]
(voir proposition \ref{old125} ou bien \cite[prop 3.2.1]{PR1}). On a 
$V^{G_K}\simeq H^1(\Gamma_n, V^{H_K})$ et
$H^2_{\mathrm{Iw}}(K,V)^{\Gamma_n} \simeq 
((V^*(1)^{H_K})^*)^{\Gamma_n} \simeq (V^*(1)^{G_K})^*$. La fl\`eche $H^1(K_n,V) \to H^2_{\mathrm{Iw}}(K,V)^{\Gamma_n}$ est duale de l'application d'inflation $V^*(1)^{G_K} \to H^1(K_n,V^*(1))$.

\begin{prop}\label{old354}
On a un diagramme commutatif dont les fl\`eches horizontales sont des isomorphismes :
\[ \xymatrix{
0 \ar[d]&&0\ar[d]\\
\frac{\calD(V)^{\Delta=0}_{\Gamma_n}}{\ker (\Exp^{\eps}_{V,h,n})} 
\ar[d]^{\tilde{\Xi}^{\eps}_{V,n}}
\ar[rr]^{\Exp^{\eps}_{V,h,n}}
\ar @{}[drr] |{(1)}
&&\mathrm{Im}(\Exp^{\eps}_{V,h,n})\ar[d]\\
\frac{t_V(K_n)}{(\dcris (V)^{\varphi =1}/V^{G_K})}
\ar[d]^{\mathrm{Tr}_{K_n/K}} 
\ar[rr]^-{(h-1)!\exp_{V,K_n}}
\ar @{}[drr] |{(2)}
&&\frac{H^1_e(K_n,V)+V^{G_K}}{V^{G_K}} \ar[d]\\
(V^*(1)^{G_K})^* 
\ar[rr]^{(h-1)! \log^{-1} \chi (\gamma_n)} 
\ar[d]  
&&H^2_{\mathrm{Iw}}(K,V)^{\Gamma_n}\ar[d]\\
0&&0} \]
\end{prop}

\begin{proof} 
Il r\'esulte du th\'eor\`eme \ref{old313} que le carr\'e (1) du diagramme est commutatif. D'apr\`es le lemme \ref{old352}, le diagramme :
\[ \xymatrix{
V^*(1)^{G_K} \times \ddR^{K_n}(V) 
\ar[d]^{(\log^{-1} \chi (\gamma_n),\mathrm{id})}
\ar[rrr]^{(\mathrm{inf}_{K_n/K}, \exp_{V,K_n})}
&&& H^1(K_n, V^*(1))\times H^1 (K_n,V)
\ar[d]\\
\ddR^{K_n}(V^*(1))\times \ddR^{K_n}(V)
\ar[rrr]^-{\mathrm{Tr}_{K_n/\Qp}[\cdot,\cdot]}
&&&\Qp } \]
est commutatif, ce qui fait que l'application compos\'ee :
\[ \ddR^{K_n}(V) \overset{\exp_{V,K_n}}{\longrightarrow} H^1(K_n,V) \to (V^*(1)^{G_K})^* \]
co\"{\i}ncide avec l'application :
\[ \ddR^{K_n}(V) \overset{\log^{-1} \chi (\gamma_n)  \mathrm{Tr}_{K_n/K}}{\longrightarrow}
\dcris (V) \to (V^*(1)^{G_K})^*. \]

On en d\'eduit que le carr\'e (2) du diagramme commute. Il est clair que toutes les fl\`eches horizontales sont des isomorphismes. Comme la colonne de gauche est exacte, la colonne de droite l'est aussi.
\end{proof}

\begin{coro}\label{old355}
On a un isomorphisme canonique :
\[ \mathrm{coker} (\Exp^{\eps}_{V,h,n}) \simeq \frac{H^1(K_n,V)}{H^1_e(K_n,V)+V^{G_K}}. \]
\end{coro}

\begin{proof} 
Il suffit d'appliquer le lemme du serpent au diagramme :
\[ \xymatrix{
0 \ar[r] &\mathrm{Im}(\Exp^{\eps}_{V,h,n})\ar[r]\ar[d]
&\frac{H^1_e(K_n,V)+V^{G_K}}{V^{G_K}}
\ar[r] \ar[d]
&H^2_{\mathrm{Iw}}(K,V)^{\Gamma_n} \ar[r]\ar @{=}[d]
&0\\
0\ar[r]
&\frac{H^1_{\mathrm{Iw}}(K,V)_{\Gamma_n}}{V^{G_K}}
\ar[r]
&\frac{H^1(K_n,V)}{V^{G_K}} 
\ar[r]
&H^2_{\mathrm{Iw}}(K,V)^{\Gamma_n}
\ar[r] &0.} \]
\end{proof}

On note $B_{V,h,n}^\eps : \ker (\Exp^{\eps}_{V,h,n}) \to \mathrm{coker}(\Exp^{\eps}_{V,h,n})$ l'application d\'eduite de $\Exp^{\eps}_{V,h,n}$.

\begin{prop}\label{old357} 
On a un diagramme commutatif dont les fl\`eches horizontales sont des isomorphismes :
\[ \xymatrix{
\ker \tilde{\Xi}^{\eps}_{V,n}
\ar[d]\ar[rrr]
\ar @{}[drrr] |{(1)}
&&& \frac{H^1_g(K_n,V)}{H^1_e(K_n,V)}
\ar[d]\\ 
\ker (\Exp^{\eps}_{V,h,n}) \ar[rrr]^{B^\eps_{V,h,n}} \ar[d]^{\tilde{\Xi}^{\eps}_{V,n}}
\ar @{}[drrr] |{(2)}
&&&\mathrm{coker}(\Exp^{\eps}_{V,h,n})\ar[d]^{\exp^*_{K_n,V}}\\
\frac{\Fil^0 \ddR^{K_n}(V)}{V^{G_K}}
\ar[rrr]^{(h-1)!\log^{-1}\chi (\gamma_n)}
&&&\frac{\Fil^0 \ddR^{K_n}(V)}{V^{G_K}}} \]
\end{prop}

\begin{proof} 
La commutativit\'e du deuxi\`eme carr\'e est d\'emontr\'ee dans \cite[lemme 3.5.9]{PR1} en utilisant la loi de r\'eciprocit\'e explicite. Le reste est une cons\'equence imm\'ediate de la proposition suivante. 
\end{proof}

\begin{prop}\label{old358} 
Le diagramme :
\[ \xymatrix{
\frac{\dcris (V)}{(1-\varphi)\dcris (V)}
\ar[rr]^{e_{V,f/e,n}}
\ar[d]
&&\frac{H^1_f(K_n,V)}{H^1_e(K_n,V)}
\ar[d]\\
\ker (\tilde{\Xi}^{\eps}_{V,n})
\ar[rr]^{B^\eps_{V,h,n}}
\ar[d]
&&\frac{H^1_g(K_n,V)}{H^1_e(K_n,V)}
\ar[d]\\
\dcris (V)^{\varphi=p^{-1}} 
\ar[rr]^{e_{V,g/f,n}}
&&\frac{H^1_g(K_n,V)}{H^1_f(K_n,V)},} \]
o\`u  
\begin{align*}
e_{V,f/e,n}(a)&=(h-1)!p^{-n}\exp_{V,f/e}( a),\\
e_{V,g/f,n}(b)&=(h-1)!\log^{-1}\chi (\gamma_n)(\exp^*_{V,g/f})^{-1}(b),
\end{align*}
est commutatif.
\end{prop}

\begin{proof} 
Si $a\in \dcris(V)/(1-\varphi)\dcris (V)$, on choisit un \'el\'ement
$f(X)\in \calD(V)$ v\'erifiant $\Delta f(X)=a$ et l'on pose
$g(X)= (\gamma_n-1)f(X)$ ce qui fait que $g(X)$ est l'image de $a$ dans
$\ker (\tilde{\Xi}_{V,n}^\eps) \subset \calD(V)^{\Delta=0}_{\Gamma_n}$.
Soit $F(X)\in \calH(V)$ un \'el\'ement v\'erifiant l'\'equation
$(1-\varphi)F(X)=f(X)-a$. Si on pose $\alpha (X)= (\partial^h\otimes e_h) f(X)$
et $A(X)= (\partial^h\otimes e_h) F(X)$, alors $A(X)$ v\'erifie $(1-\varphi)A(X)=\alpha (X)$. 

Soit $h$ un entier tel que $\Fil^0 \dcris (V(-h))=\dcris (V(-h))$ et :
\[ \Sigma^{\eps}_{V(-h),h+k,,m} : \calH (V(-h)) \to H^1(K_m, V(k)) \]
le syst\`eme d'applications construit dans \cite[\S4.2-4.3]{Ben} et dont la construction a \'et\'e rappel\'ee au paragraphe \ref{altexpo}. Posons :
\[ z_{k,m}=(-1)^h \Sigma^{\eps}_{V(-h),h+k,m} ((\sigma \otimes \varphi)^{-m}A(X)), \]
et notons $\overline{z}_{k,m}$ son image dans $H^1(K_m,V(k)) / H^1(\Gamma_m, V(k)^{H_K})$. Le th\'eor\`eme 4.3 de \cite{Ben} montre que
$\mathrm{cor}_{K_{n+1}/K_n} (\overline{z}_{k,n+1})=\overline{z}_{k,n}$ pour tout $n \geq 1$ et qu'il existe $s\geq 0$ tel que la suite :
\[ p^{(s-j)m} \sum_{k=0}^j (-1)^k \binom{j}{k} \mathrm{Tw}^{\eps}_{-k}\circ 
\mathrm{res}_{K_\infty/K_n}(z_{m,k}) \]
converge vers $0$ quand $m \to \infty$. On v\'erifie  que
$\overline{z}_{k,m} \in ( H^1_{\mathrm{Iw}}(K,V(k))/V(k)^{H_K})_{\Gamma_m}$
(par exemple, on peut utiliser les arguments de \cite{Ben}, pour montrer
que le cup-produit de $\overline{z}_{k,m}$ avec les \'el\'ements de
$V^*(1-k)^{G_K}$ est nul) et il existe donc un unique \'el\'ement 
$z \in \calH(\Gamma) \otimes_{\Lambda} (H^1_{\mathrm{Iw}}(K,V(k))/ V(k)^{H_K})$ tel que $\mathrm{pr}_{V(k),m} (\mathrm{Tw}^{\eps}_{k}(z))= \overline{z}_{k,m}$ pour tous $k \in \ZZ$ et $m \geq 1$.

L'\'el\'ement $B(X) =(\gamma_n-1)A(X)$ v\'erifie $(1-\varphi) B(X)= (\partial^h \otimes e_h)g(X)$ et on a :
\[ \mathrm{pr}_{V(k),m} (\mathrm{Tw}^{\eps}_{V,k} \circ \Exp^{\eps}_{V,h}(g))=
(-1)^h \Sigma^{\eps}_{V(-h),h+k,m}((\sigma \otimes \varphi)^{-m}B(X)), \]
et donc $(\gamma_n-1) z= \Exp^{\eps}_{V,h}(g)$ et $B^{\eps}_{V,h,n} (a) = z_{0,n} \mod {H^1_e (K_n,V)}$. D'autre part, le m\^eme argument que dans \cite[\S 4.4.5]{Ben} montre que :
\[ z_{0,n} = (h-1)! p^{-n}\exp_{V,K_n}(-\varphi^{-n}(a), 
(\sigma\otimes \varphi)^{-n}F(\zn -1)), \]
et on en d\'eduit la commutativite du premier carr\'e du diagramme.

D\'emontrons la commutativit\'e du deuxi\`eme carr\'e.  Fixons un entier $k \geq 1$ sup\'erieur \`a la longueur de la filtration de Hodge de $V^*(1)$ et tel que $\Fil^{-k}\dcris(V^*(1)) = \dcris (V^*(1))$. La loi de r\'eciprocit\'e s'\'ecrit :
\[ \pscal{\Exp^{\eps}_{V,h}(f),\Exp^{\eps^{-1}}_{V^*(1),k}(g^{\iota})}  (1+X) =
-\prod_{i=0}^{h-1}\ell_i  \prod_{j=0}^{k-1}\ell_i^{\iota} (f \star_{\calD} g). \]

Si $f\in \ker (\tilde{\Xi}_{V,n}^{\eps})$, alors $\Exp^{\eps}_{V,h}(f)=
(\gamma_n-1)x$ o\`u $x\in \calH(\Gamma)\otimes H^1_{\mathrm{Iw}}(K,V)$
et $B^{\eps}_{V,h,n}(f)=\mathrm{pr}_{V,n}(x) \mod{H^1_e(K_n,V)}$. Soient
$b\in \dcris (V^*(1))/(1-\varphi)\dcris (V^*(1))$ et $\beta (X) \in 
\calD(V^*(1))$ un \'el\'ement v\'erifiant $\Delta \beta (X)=b$. Posons $g(X)= (\gamma_n-1)\beta(X)$. On a alors :
\begin{align*}
\Exp^{\eps^{-1}}_{V^*(1),k}(g) & = (\gamma_n-1)y ,\\
\pscal{\Exp^{\eps}_{V,h}(f),\Exp^{\eps^{-1}}_{V^*(1),k}(g^{\iota})} & =
(\gamma_n-1)^2 \pscal{x,y^{\iota}}.
\end{align*}
 
On a $\mathrm{pr}_{V^*(1),n}(y)\in H^1_f(K_n,V^*(1))$ et
$B^\eps_{V^*(1),h,n}(b)=\mathrm{pr}_{V^*(1),n}(y) \mod{H^1_e(K_n,V^*(1))}$. 
Il est facile de voir que $g^{\iota}$ v\'erifie aussi $\Delta (g^{\iota})=b$
et que $B^\eps_{V,h,n}(b)$ ne d\'epend pas du choix de $\eps$ d'o\`u :
\[ B_{V^*(1),h,n}^\eps(b)=\mathrm{pr}_{V^*(1),n}(y^{\iota}) \mod{H^1_e(K_n,V^*(1))}. \]

Comme $ \ell_0 \equiv -\frac{\gamma_n-1}{\log \chi (\gamma_n)} \mod{(\gamma_n-1)^2}$ et comme le coefficient de $(1+X)$ dans le polyn\^ome d'interpolation de $f\star_{\calD}\beta$ modulo $(1+X)^{p^n}-1$ est \'egal \`a :
\[ \frac{1}{p^n}\mathrm{Tr}_{K/\Qp} \sum_{\zeta \in \mu_{p^n}}
[f(\zeta-1),\beta^{\iota}(\zeta^{-1}-1)]_V \]
(voir, par exemple, \cite[prop 4.3.2]{PR1}), on d\'eduit de la loi de r\'eciprocit\'e la formule suivante :
\[ (B^\eps_{V,h,n}(f),B^\eps_{V^*(1),h,n}(b))_{V,K_n}=\frac{(h-1)!(k-1)!}{p^n\log 
\chi (\gamma_n)}\mathrm{Tr}_{K/\Qp} \sum_{\zeta \in \mu_{p^n}}
[f(\zeta-1),\beta^{\iota}(\zeta^{-1}-1)]_V. \]
 
Soit $F(X)$ un \'el\'ement tel que $(1-\varphi)F(X)=f(X)$. Comme $\tilde{\Xi}^{\eps}_{V,n}(f)=0$, on a $F(\zn-1)\in \dcris (V)^{\varphi=1}$. On peut modifier $F(X)$ par cet  \'el\'ement et on a alors $F(\zn-1)=0$. Comme $F(X)$ v\'erifie l'\'equation $\sum_{\zeta^p=1} F(\zeta (1+X)-1) = p F^{\varphi}(X)$,
on a $F(\zeta_{p^m}-1)=0$ pour tout $1 \leq m \leq n$. On en d\'eduit que $f(\zeta_{p^m}-1)=0$  si $2\leq m \leq n$ et $f(\zeta_p-1)=-F^{\sigma}(0)$ ce qui fait que :
\[ \sum_{\zeta \in \mu_{p^n}}
[f(\zeta-1),\beta^{\iota}(\zeta^{-1}-1)]_V=
\sum_{\substack{\zeta^p=1\\ \zeta \neq 1}} 
[-F^{\sigma}(0),\beta^{\iota}(\zeta^{-1}-1)]_V +[f(0),b]_V=
[f(0),b]_V, \]
d'o\`u :
\[ (B^\eps_{V,h,n}(f),B^\eps_{V^*(1),h,n}(b))_{V,K_n}=\frac{(h-1)!(k-1)!}{p^{2n}\log \chi (\gamma_n)}\mathrm{Tr}_{K/\Qp}[\varphi^{-n} f(0), \varphi^{-n}(b)]_V. \]

Comme $B^\eps_{V^*(1),h,n}(b)=e_{V^*(1),f/e,n}(b)$, la commutativit\'e
du deuxi\`eme carr\'e r\'esulte de la d\'efinition de l'application $\exp^*_{V,g/f}$.
\end{proof}

\begin{proof}[D\'emonstration du th\'eor\`eme \ref{old326}]
Soit $h$ un entier sup\'erieur \`a la longueur de la filtration de Hodge de
$V$. L'application $\Exp^{\eps}_{V,h}$ est semi-simple et les r\'esultats pr\'ec\'edents fournissent un diagramme commutatif :
\[ \xymatrix{
\det_{\Lambda_{\Qp}} \calD(V) \otimes 
\det_{\Lambda_{\Qp}} \RR\Gamma_{\mathrm{Iw}} (K,V)
\ar[d]\ar[r] & \calH(\Gamma)
\ar[d]\\
\det_{\Qp [G_n]} \calD(V)_{\Gamma_n} \otimes 
\det_{\Qp[G_n]} \RR\Gamma (K_n,V)
\ar[d]^{\kappa_{V,n}}\ar[r] & \Qp[G_n]
\ar[d]^{\mathrm{id}}\\
\det_{\Qp [G_n]} \ddR^{K_n}(V) \otimes 
\det_{\Qp[G_n]} \RR\Gamma (K_n,V)
\ar[r] & \Qp[G_n]. } \tag{eq9} \]

La deuxi\`eme ligne du diagramme est induite par l'application
$\Exp^{\eps}_{V,h,n}$. Le th\'eor\`eme \ref{old324} entra\^{\i}ne que l'image de $\det_\Lambda \calD_M(V) \otimes \det_\Lambda \RR\Gamma_{\mathrm{Iw}} (K,T)$ dans $\calH(\Gamma)$ s'\'ecrit sous la forme $\mathbf{\Gamma}_h(V)\sum_{i=0}^{p-2} a_i \Lambda_i$, o\`u $a_i \in \Qp$ et $\Lambda_i = \Lambda_{\Gamma_1} \delta_i$.

Les propositions pr\'ec\'edentes et le lemme \ref{old352} montrent que la troisi\`eme ligne co\"{\i}ncide avec l'application :
\[ \mathbf{\Gamma}_h^*(V) \Gamma^*(V)(\mathrm{id}+ (p^{-n\dim \dcris (V)^{\varphi=1}}-1)e_{\eta_0}) \delta'_{V,K_n/K}. \]

La conjecture $C_{\mathrm{Iw}}(K_\infty/K,V)$ est vraie si et seulement si $a_i \in \Zp^\times$ c'est-\`a-dire si et seulement si 
l'image de $\det_{\Zp[G_n]} \calD_M(V)_{\Gamma_n} \otimes 
\det_{\Zp[G_n]} \RR\Gamma (K_n,T)$ dans $\Qp[G_n]$ est engendr\'ee par $\mathbf{\Gamma}^*_h(V)$ (voir lemme \ref{old3212}). Les propositions \ref{old333} et \ref{old349} entra\^{\i}nent que 
cela \'equivaut \`a dire que l'application $\delta'_{V,K_n/K}$ envoie $\det_{\Zp[G_n]} M_n \otimes \det_{\Zp[G_n]} \RR\Gamma (K_n,T)$
sur $\beta_{V,K_n/K}(M,T)^{-1}$ d'o\`u le th\'eor\`eme.
\end{proof}

\Subsection{R\'esultats principaux}\label{main}

Dans ce paragraphe, on d\'emontre la conjecture $C_{\mathrm{Iw}}(K_\infty/K,V)$.

\begin{theo}\label{old431} 
Si $V$ est une repr\'esentation cristalline dont les oppos\'es des poids
de Hodge-Tate sont $0=r_1\leq  r_2 \leq \cdots \leq r_d=h$, et qui n'a pas de sous-quotient isomorphe \`a $\Qp(m)$, alors :
\begin{enumerate}
\item l'application $\varphi^{-1}$ induit un isomorphisme :
\[ i_V : \frac{\dfont(V)^{\psi=1}}{{(\varphi^*\nwach(V))}^{\psi=1}} \iso
\oplus_{k=1}^h (K_1t^{-k} \otimes_K  \Fil^k  \dcris (V)). \]
\item On a :
\[ \det_\Lambda \left ( \frac{\dfont(T)^{\psi=1}}{(\varphi^*\nwach(T))^{\psi=1}} \right )
= \prod_{k=1}^h  (\det_\Lambda (\Zp[\Delta]\otimes \Zp(-k)))^{\dim_{\Qp} \Fil^k \dcris (V)}. \]
\end{enumerate}
\end{theo}

Nous montrons ce th\'eor\`eme un peu plus bas. Si $f,g \in \Lambda_{\Gamma_1}$, on \'ecrit $f\sim g$ si $f$ et $g$ sont associ\'ees, c'est-\`a-dire s'il existe $u\in \Lambda_{\Gamma_1}^\times$
tel que $f=gu$. De m\^eme, si $a,b\in\ZZ$ on \'ecrit $a \sim_p b$ si $a$ et $b$ sont associ\'es dans $\Zp$. Si $M$ est un $\Lambda$-module de torsion et de type fini, on note 
$\mathrm{car}_\Lambda(M)$ son polyn\^ome caract\'eristique. Rappelons que $\det_\Lambda(M) = \mathrm{car}_\Lambda(M)^{-1}\Lambda$. 

\begin{prop}\label{old432}
Si $V$ est une repr\'esentation cristalline v\'erifiant les conditions
du th\'eor\`eme \ref{old431} ci-dessus, alors :
\begin{enumerate}
\item l'application $\varphi^{-1}$ induit une injection :
\[ i_V : \frac{\dfont(V)^{\psi=1}}{{(\varphi^*\nwach(V))}^{\psi=1}} \hookrightarrow
\oplus_{k=1}^h (K_1t^{-k} \otimes_K  \Fil^k \dcris (V)). \]
\item On a :
\[ \mathrm{car}_\Lambda \left ( \frac{\dfont(T)^{\psi=1}}{(\varphi^*\nwach(T))^{\psi=1}} \right )
 \Bigg| \prod_{k=1}^h  (\mathrm{car}_\Lambda (\Zp[\Delta]\otimes \Zp(-k)))^{\dim_{\Qp} \Fil^k \dcris (V)}. \]
\end{enumerate}
\end{prop}

\begin{proof}
La proposition \ref{old414} nous dit que : 
\[ [\bhol \otimes_{\bplus_K} \nwach(V) : \bhol \otimes_K \dcris (V) = [(t/X)^{r_1};\hdots ; (t/X)^{r_d}], \]
et on en d\'eduit que : 
\[ \nwach(V) \subset \left (\frac{X}{t} \right )^h  \bhol \otimes_K \dcris (V), \] 
et que le plongement de $\bhol$ dans $K[[t]]$ donne un isomorphisme 
$K[[t]]  \otimes_{\bplus_K} \nwach(V) \simeq K[[t]] \otimes_K \dcris(V)$.
Comme $q^h \nwach(T) \subset \varphi^*(\nwach(T))$, on a $X^{-h} \nwach(T) \subset \varphi (X)^{-h}\varphi^*(\nwach(T))$ d'o\`u :
\[ \dfont(T)^{\psi=1}=\left (\frac{1}{X^{h}} \nwach(T)
\right )^{\psi=1}=
\left (\frac{1}{\varphi (X)^{h}}\varphi^*(\nwach(T))\right )^{\psi=1}. \]
Comme $\varphi^{-1} (1/X) \in \bdR^+$, l'application $\varphi^{-1}$ induit une injection :
\[ \varphi^{-1} : \dfont(T)^{\psi=1} \to \Fil^0 (K_1((t))\otimes_K \dcris(V)). \]
Pour all\'eger les notations, on pose :
\[ D = \Fil^0 (K_1((t))\otimes_K \dcris(V)) = 
\oplus_{m=-\infty}^h K_1t^{-m} \otimes_K \Fil^m \dcris (V), \]
et on d\'efinit une filtration croissante de $D$ par des sous-espaces $D_k$ : 
\[ D_k= \oplus_{m=-\infty}^k K_1t^{-m} \otimes_K \Fil^m \dcris (V). \]
On a $D_h=D$ et $D_k/D_{k-1} \simeq \Fil^k  K_1 t^{-k} \otimes_K \dcris(V)$.
On d\'efinit aussi une filtration sur $(X^{-h} \nwach(V))^{\psi=1}$ par des sous-espaces $\nwach_k(V)$ en posant :
\[ \nwach_k(V) = \left (\frac{1}{\varphi (X)^k}\varphi^*(\nwach(V))\right )^{\psi=1}, \]
et  $\varphi^{-1}$ induit alors une injection $\nwach_k(V) \hookrightarrow D_k$. 
Pour montrer que l'application $\nwach_h(V)/\nwach_0(V)\to D_h/D_0$ est injective,
il suffit de montrer que les applications $i_{V,k}  : \nwach_k(V)/\nwach_{k-1}(V) \to K_1t^k \otimes_K \Fil^k \dcris (V)$ sont injectives pour $k=1,\hdots ,h$. Pour cela, soit $n_1,\hdots ,n_d$ une base de $\nwach(V)$ et soit :
\[ x  = \frac{1}{\varphi (X)^k}(a_1\varphi (n_1)+\cdots +
a_d\varphi (n_d)) \]
un \'el\'ement de $\ker i_{V,k}$. On a alors :
\[ \frac{\varphi^{-1}(a_1) n_1+\cdots +
\varphi^{-1}(a_d)n_d}{X^k}\in  t^{-k+1}K_1[[t]] \otimes_K \dcris(V), \]
et comme  $n_1,\hdots ,n_d$ forment une base de $K[[t]]  \otimes_K \dcris(V)$, on a $a_i \equiv 0\mod q$. Si l'on \'ecrit $a_i= qb_i$ avec $b_i = \sum_{j=0}^\infty b_{ij}X^j$, alors la condition $\psi (x)=x$ s'\'ecrit :
\[ \left (\frac{b_1}{X} \right ) \varphi (n_1)+ \cdots + \left (\frac{b_d}{X} \right ) \varphi (n_d)=
q^{k-1}\left (\psi \left (\frac{b_1}{X} \right ) n_1 + \cdots +\psi \left (\frac{b_d}{X} \right ) 
n_d \right ). \]
On en d\'eduit que $y=\sum_{i=1}^d \varphi^{-1} (b_{i0})n_i$ appartient \`a
$\dcris (V)^{\varphi=p^{k-1}}$. Mais il appartient aussi \`a $\Fil^{k-1}\dcris (V)$ par construction de $i_{V,k}$, donc \`a $V(k-1)^{G_{K}}(1-k)=0$. Ceci montre que $i_{V,k}$ est injective et 
la deuxi\`eme assertion s'en d\'eduit.
\end{proof}

\begin{proof}[D\'emonstration du th\'eor\`eme \ref{old431}] 
Nous d\'emontrons d'abord le (2). Posons :
\begin{align*}
\sum_{i=0}^{p-2} \delta_i \otimes f_{T,i}(\gamma_1-1) & = \mathrm{car}_\Lambda \left (\frac{\dfont(T)^{\psi=1}}{(\varphi^*\nwach(T))^{\psi=1}}\right ),\\
\sum_{i=0}^{p-2} \delta_i \otimes g_{T,i}(\gamma_1-1) & = \prod_{m=1}^h \mathrm{car}_\Lambda \left (\Zp[\Delta_K]\otimes \Zp(-m) \right )^
{\dim_{\Qp}\Fil^m \dcris(V)},
\end{align*}
o\`u $\delta_i = \sum_{g\in \Delta_K} \chi^{-i}(g)g$.
Posons aussi $f_T(\gamma_1-1)= \prod_{i=0}^{p-2} 
f_{T,i}(\gamma_1-1)$ et $g_T(\gamma_1-1) = \prod_{i=0}^{p-2} 
g_{T,i}(\gamma_1-1)$. On voit que pour $i=0,\hdots,p-2$, on a :
\[ g_{T,i}(\gamma_1-1)=\prod_{m=1}^h (\gamma_1-\chi (\gamma_1)^{-m})^{\dim_{\Qp}\Fil^m \dcris(V)}, \]
d'o\`u l'on d\'eduit que :
\[ g_{T,i}(\chi(\gamma_1)^{-k}-1) \sim_p 
p^{[K:\Qp]t_H(V)} \prod_{m=1}^h (k-m)^{\dim_{\Qp}\Fil^m \dcris(V)}. \]
Comme $t_H(V)+t_H(V^*(-h))=h\dim_{\Qp}V$ et comme : 
\[ \dim_{\Qp}\Fil^m\dcris (V)+\dim_{\Qp}\Fil^{1-m}\dcris (V^*)=[K:\Qp]\dim_{\Qp}V,  \]
on voit que :
\[ g_{T}(\chi(\gamma_1)^{-k}-1)g_{T^*(-h)}(\chi(\gamma_1)^{k-h-1}-1)
\sim_p  \left ({p^h\frac{\Gamma^*(k)}{\Gamma^*(k-h)}} \right )^{[K_1:\Qp]\dim (V)}. \tag{eq12} \]
D'autre part, comme : 
\[ \left (\frac{\dfont(T(k))^{\psi=1}}{(\varphi^*\nwach(T)(k))^{\psi=1}}\right )^{\Gamma_1}=0, \] si $k\notin [1,h]$, on a :
\[ f_{T}(\chi(\gamma_1)^{-k}-1)=[\dfont(T(k))^{\psi=1}_{\Gamma_1}:
(\varphi^*\nwach(T)(k))^{\psi=1}_{\Gamma_1}], \]
et le corollaire \ref{old4211} donne :
\[ f_{T}(\chi(\gamma_1)^{-k}-1)f_{T^*(-h)}(\chi(\gamma_1)^{k-h-1}-1)
\sim_p \left (
{p^h\frac{\Gamma^*(k)}{\Gamma^*(k-h)}} \right )^{[K_1:\Qp]\dim (V)}.
\tag{eq13} \]

Comme la divisibilit\'e $f_{T,i}(\gamma_1-1)  \mid g_{T,i}(\gamma_1-1)$ est d\'emontr\'ee dans la proposition \ref{old432}, les formules (eq12) et (eq13) entra\^{\i}nent :
\[ f_{T}(\chi(\gamma_1)^{-k}-1) \sim_p g_{T}(\chi(\gamma_1)^{-k}-1) \]
pour tout $k\notin [1,h]$ et donc $f_{T}(\gamma_1-1)  \sim g_{T}(\gamma_1-1)$
et l'assertion (2) est d\'emontr\'ee.

Montrons maintenant le (1). Si on appelle $Y$ le conoyau de l'injection :
\[ \frac{\dfont(V)^{\psi=1}}{{(\varphi^*\nwach(V))}^{\psi=1}}\hookrightarrow
\oplus_{k=1}^h (K_1t^{-k} \otimes_K \Fil^k \dcris (V)), \]
alors $\det_{\Lambda_{\Qp}}(Y)= \Lambda_{\Qp}$ par le (2), d'o\`u $Y=0$ car  $\Lambda_{\Qp}$ est un anneau principal, ce qui montre le (1). 
\end{proof}

\begin{coro}\label{old434} 
Si $k \notin [1,h]$, alors le conoyau de l'application :
\[ \Omega_{T,k,1}^\eps : \calD(T(k))_{\Gamma_1} \to
H_{\mathrm{Iw}}^1(K_1,T(k))_{\Gamma_1} \]
est isomorphe \`a :
\[ \oplus_{m=1}^h ( (\ZZ/(k-m)p\ZZ) [\Delta_K] )^{\dim_{\Qp} \Fil^m \dcris(V)}. \]
En particulier,
\[ \det_{\Zp[\Delta_K]}(\Omega_{T,k,1}^\eps)=
p^{[K:\Qp]t_H(V)}
\prod_{m=1}^h (k-m)^{\dim_{\Qp}\Fil^m \dcris(V)}\Zp[\Delta_K]. \]
\end{coro}

Nous pouvons maintenant d\'emontrer le r\'esultat principal de cet article.

\begin{theo}\label{new345}
Si $V$ est une repr\'esentation cristalline de $G_K$, alors :
\begin{enumerate}
\item la conjecture $C_{\mathrm{Iw}}(K_\infty/K,V)$ est vraie; 
\item la conjecture $C_{\mathrm{EP}}(L/K,V)$ est vraie pour toute extension $L/K$ contenue dans $K_\infty$.
\end{enumerate}
\end{theo}

\begin{proof} 
Comme la conjecture $C_{\mathrm{Iw}}(K_\infty/K,V) $ est stable par suites exactes et comme le cas $V=\Qp(m)$ a \'et\'e d\'emontr\'e dans \cite[page 143]{PR1}, on peut supposer que 
$V$ n'a pas de sous-quotient fix\'e par $H_K$. Sous ces hypoth\`eses, $H_{\mathrm{Iw}}^1(K,T(k))^{\Gamma_1}=0$, le groupe $H_{\mathrm{Iw}}^2(K,T(k))$ est fini et pour $n=1$ le diagramme (eq9) du paragraphe \ref{same2} s'\'ecrit :
\[ \xymatrix{  
\det_{\Lambda_{\Qp}}\calD(V(k))\otimes 
\det_{\Lambda_{\Qp}}^{-1}H_{\mathrm{Iw}}^1(K,V(k))\ar[r] \ar[d] & \calH(\Gamma)
\ar[d]\\
\det_{\Qp[\Delta_K]}\calD(V(k))_{\Gamma_1}\otimes 
\det_{\Qp[\Delta_K]}^{-1}H_{\mathrm{Iw}}^1(K,V(k))_{\Gamma_1}\ar[r]  &
\Qp[\Delta_K].} 
\tag{eq14} \]
La premi\`ere ligne de ce diagramme est induite par l'application
$\Exp^{\eps}_{V(k),k}=(-1)^k\mathrm{Tw}_{V,k}^\eps \circ \Exp^{\eps}_{V,0}
\circ (\partial^k \otimes e_k)$ et le th\'eor\`eme \ref{old324} (la conjecture
$\delta_{\Qp}(V)$ de Perrin-Riou) entra\^{\i}ne que l'image de
$\det_{\Lambda} \calD(T(k))\otimes \det_{\Lambda} \RR\Gamma(K,T(k))$ 
dans $\calH(\Gamma)$ s'\'ecrit sous la forme $\prod_{j=1}^h (\ell_{k-j})^{{\mathrm{dim}}_
{\Qp}\Fil^j \dcris (V)} \sum_{i=0}^{p-2} a_i\Lambda_i$, o\`u $a_i\in \Qp$ et $\Lambda_i=\Lambda_{\Gamma_1}\delta_i$. 
Pour conclure, il suffit de montrer que $a_i\in \Zp^*$ pour $i=0,\hdots,p-2$. Le (2) de la proposition \ref{old423} montre que la deuxi\`eme ligne de (eq14)
est induite par $(-1)^k \Omega_{T,k}^{\eps}((\sigma \otimes \varphi)^{-1}
\circ (\partial^k \otimes e_k))$ et on en d\'eduit que l'image de
$\det_{\Zp[\Delta_K]}\calD(T(k))_{\Gamma_1}\otimes 
\det_{\Zp[\Delta_K]}^{-1} H_{\mathrm{Iw}}^1(K,T(k))_{\Gamma_1}$ 
dans $\Qp[\Delta_K]$ est \'egale \`a $\prod_{j=1}^h (k-j)^{{\mathrm{dim}}_
{\Qp}\Fil^j \dcris (V)}\sum_{i=0}^{p-2} a_i\delta_i\Zp$.

D'autre part, comme le $\varphi$-module filtr\'e $\dcris(V)$ est admissible, on a
$[M : \varphi (M)] = p^{[K:\Qp]t_H(V)}$ et en comparant avec la formule du corollaire \ref{old434}, on obtient que $a_i\in \Zp^*$ et la conjecture $C_{\mathrm{Iw}}(K_\infty/K,V)$ est d\'emontr\'ee. 

Le th\'eor\`eme \ref{old325} et la proposition \ref{old242} montrent enfin
que la conjecture $C_{\mathrm{EP}}(L/K,V)$ est vraie pour toute extension $L/K$ contenue dans $K_\infty$.
\end{proof}

\begin{coro}\label{old436} 
Si $V$ est une repr\'esentation cristalline de $G_K$, alors :
\begin{enumerate}
\item la conjecture $C_{\mathrm{EP}}(L,V)$ est vraie pour toute extension
$L/K$ contenue dans $\Qp^{\mathrm{ab}}$. 
\item la conjecture $C_{\mathrm{EP}}(K,V(\eta))$ est vraie pour tout caract\`ere
de Dirichlet $\eta$ de $\Gamma$.
\end{enumerate}
\end{coro}

\begin{proof} 
C'est une cons\'equence imm\'ediate de la proposition \ref{old242}.
\end{proof}

\end{document}